\documentclass{article}

\usepackage[english]{babel}

\usepackage[scale=0.85]{geometry}

\usepackage{bm}
\usepackage{graphicx}
\usepackage{amsthm}
\usepackage{thmtools}
\usepackage{amssymb}
\usepackage{mathtools}
\usepackage{subcaption}
\usepackage{float}
\usepackage{multirow}
\usepackage[colorlinks=true, allcolors=blue]{hyperref}
\usepackage[capitalise,nameinlink,noabbrev]{cleveref}

\usepackage{siunitx}
\makeatletter
\newcommand{\printslopeinv}[4]{
    \tikzset{fixed point arithmetic}
    \def\nero@printslope@orderlist{#1}
    \edef\nero@printslope@xpos{#2}
    \edef\nero@printslope@ypos{#3}
    \edef\nero@printslope@width{#4}
    \pgfmathparse{\nero@printslope@xpos+\nero@printslope@width}
    \edef\nero@printslope@px{\pgfmathresult}
    \edef\nero@printslope@py{\nero@printslope@ypos}
    \edef\nero@printslope@qx{\pgfmathresult}
    \edef\nero@printslope@ry{\nero@printslope@ypos}
    \foreach \nero@printslope@order in {#1}{
        \pgfmathparse{
        ((\nero@printslope@px/\nero@printslope@xpos)^(\nero@printslope@order))*\nero@printslope@ypos}
        \edef\nero@printslope@qy{\pgfmathresult}
            \edef\nero@aux1{\noexpand\draw[line width=0.6pt]
            (axis cs:\nero@printslope@xpos,\nero@printslope@ypos)
            -- (axis cs:\nero@printslope@qx,\nero@printslope@qy)
            -- (axis cs:\nero@printslope@px,\nero@printslope@py);}
        \nero@aux1
        \pgfmathparse{10^((ln(\nero@printslope@ry)+ln(\nero@printslope@qy))/(ln(10)*2))}
        \edef\nero@printslope@labelpos{\pgfmathresult}
        \edef\nero@aux2{\noexpand\node[anchor=west] at
            (axis cs:\nero@printslope@qx,\nero@printslope@labelpos)
            {\noexpand\tiny \nero@printslope@order};}
        \nero@aux2
        \global\edef\nero@printslope@ry{\nero@printslope@qy}
    }
    \draw[line width=0.6pt] (axis cs:\nero@printslope@xpos,\nero@printslope@ypos)
        |- (axis cs:\nero@printslope@px,\nero@printslope@py);
    \pgfmathparse{10^((ln(\nero@printslope@px)+ln(\nero@printslope@xpos))/(ln(10)*2))}
    \edef\nero@printslope@labelpos{\pgfmathresult}
    \node[anchor=north] at (axis cs:\nero@printslope@labelpos,\nero@printslope@ypos) {\tiny 1};
}
\makeatother

\usepackage{amssymb}
\usepackage{mathtools}
\usepackage{pgfplots}
\usepackage{pgfplotstable}
\usepackage{datatool}
\usepackage{fp}

\pgfplotsset{
    compat=newest,
    smaller labels/.style={
        label style={font=\footnotesize},
        tick label style={font=\footnotesize}
    }
}

\tikzset{font=\small}
\usetikzlibrary{
    fpu,
    fixedpointarithmetic,
    babel,
    external,
    arrows.meta,
    plotmarks,
    positioning,
    angles,
    quotes,
    intersections,
    calc,
    spy,
    decorations.pathreplacing,
    matrix,
    fit,
}
\usepgfplotslibrary{fillbetween}

\gdef\iterator{0}

\definecolor{femcolor}{RGB}{51, 138, 55} 
\definecolor{addcolor}{RGB}{217,95,2} 
\definecolor{add2color}{RGB}{199,39,34} 
\definecolor{add3color}{RGB}{228, 155, 15} 
\definecolor{multcolor3}{RGB}{117,112,179} 
\definecolor{multcolor100}{RGB}{0,0,0} 
\definecolor{multcolor0weak}{RGB}{49, 73, 181} 
\definecolor{multcolor0strong}{RGB}{49, 181, 161} 

\newlength{\plotwidth}
\newlength{\plotheight}
\newenvironment{cvgh}[6][]{
    \begin{tikzpicture}[every node/.style={inner sep=0.5,outer sep=0.5}]
        \edef\filename{#2}
        \edef\legendcolumns{#3}
        \edef\slopes{#4}
        \edef\ypos{#5}
        \edef\normtype{#6} 
        \edef\extraoptions{#1}

        \pgfplotstableread[col sep=comma]{\filename}\datatable

        \pgfmathtruncatemacro{\secondrow}{1} 
        \pgfplotstablegetelem{\secondrow}{h}\of\datatable
        \pgfmathsetmacro{\second}{\pgfplotsretval} 

        \pgfmathtruncatemacro{\firstrow}{0} 
        \pgfplotstablegetelem{\firstrow}{h}\of\datatable
        \pgfmathsetmacro{\first}{\pgfplotsretval} 

        \pgfmathsetmacro{\diff}{\first - \second}

        \pgfmathtruncatemacro{\iterator}{\iterator+1}

        \begin{loglogaxis}[
            smaller labels,
            name = left_plot,
            axis lines = left,
            enlarge x limits={abs=10pt},
            enlarge y limits={abs=10pt},
			xmode=log,
            xlabel = {$h$},
            ylabel = {\rotatebox{270}{\normtype}},
            xlabel style={at={(ticklabel* cs:1.01)},anchor=south west},
            ylabel style={at={(ticklabel* cs:1.01)},anchor=west},
            xtick=data,
            xticklabels from table={\datatable}{h},
            width=\plotwidth, height=\plotheight,
            mark options={solid, scale=1},
            grid = major,
            legend columns=\legendcolumns,
            legend to name=leg:legendFEMCORR_\iterator,
            legend image post style={mark options={solid, scale=1}},
            \extraoptions 
        ]
        \expandafter\printslopeinv\expandafter{\slopes}{\second}{\ypos}{\diff}
    }
    {
        \end{loglogaxis}
        \node[yshift=-20pt] at (left_plot.outer south) {\pgfplotslegendfromname{leg:legendFEMCORR_\iterator}};

    \end{tikzpicture}
}

\newcommand{\cvgFEMCorrAlldegGeneral}[5]{
    \edef\fem{#1}
    \edef\add{#2}

    \begin{cvgh}{\fem}{3}{#3}{#4}{#5}
        \addlegendentry{\,FEM $\mathbb{P}_1$\;}
        \addlegendentry{\,FEM $\mathbb{P}_2$\;}
        \addlegendentry{\,FEM $\mathbb{P}_3$\;}
        \addlegendentry{\,Add $\mathbb{P}_1$\;}
        \addlegendentry{\,Add $\mathbb{P}_2$\;}
        \addlegendentry{\,Add $\mathbb{P}_3$\;}

        \addplot [style={solid}, mark=square*, mark size=2, color=femcolor, line width=0.8pt ]
        table [x=h, y=P1, col sep=comma]
            {\fem};
        
        \addplot [style={solid}, mark=*, mark size=2, color=femcolor, line width=0.8pt ]
        table [x=h, y=P2, col sep=comma]
            {\fem};
        
        \addplot [style={solid}, mark=triangle*, mark size=2, color=femcolor, line width=0.8pt ]
        table [x=h, y=P3, col sep=comma]
            {\fem};

        \addplot [style={dashed}, mark=square*, mark size=2, color=addcolor, line width=0.8pt ]
        table [x=h, y=P1, col sep=comma]
            {\add};

        \addplot [style={dashed}, mark=*, mark size=2, color=addcolor, line width=0.8pt ]
        table [x=h, y=P2, col sep=comma]
            {\add};

        \addplot [style={dashed}, mark=triangle*, mark size=2, color=addcolor, line width=0.8pt ]
        table [x=h, y=P3, col sep=comma]
            {\add};

    \end{cvgh}
}

\newcommand{\cvgFEMCorrAlldeg}[3]{
    \cvgFEMCorrAlldegGeneral{#1}{#2}{2,3,4}{#3}{$L^2$}
}

\newcommand{\cvgFEMCorrAlldegHun}[3]{
    \cvgFEMCorrAlldegGeneral{#1}{#2}{1,2,3}{#3}{Semi-$H^1$}
}

\newcommand{\cvgFEMCorrMultOnedegGeneral}[8]{
    \edef\fem{#1}
    \edef\femsec{#2}
    \edef\add{#3}
    \edef\mult{#4}
    \edef\multHundred{#5}

    \begin{cvgh}{\fem}{3}{#6}{#7}{#8}
        \addlegendentry{\,FEM $\mathbb{P}_1$\;}
        \addlegendentry{\,Mult $\mathbb{P}_1$ (M=3)\;}
        \addlegendentry{\,Add $\mathbb{P}_1$\;}
        \addlegendentry{\,FEM $\mathbb{P}_2$\;}
        \addlegendentry{\,Mult $\mathbb{P}_1$ (M=100)\;}

        \addplot [style={solid}, mark=square*, mark size=2, color=femcolor, line width=0.8pt ]
        table [x=h, y=err, col sep=comma]
            {\fem};

        \addplot [style={dotted}, mark=square*, mark size=2, color=multcolor3, line width=1.0pt ]
        table [x=h, y=err, col sep=comma]
            {\mult};

        \addplot [style={dashed}, mark=square*, mark size=2, color=addcolor, line width=0.8pt ]
        table [x=h, y=err, col sep=comma]
            {\add};

        \addplot [style={solid}, mark=*, mark size=2, color=femcolor, line width=0.8pt ]
        table [x=h, y=err, col sep=comma]
            {\femsec};

        \addplot [style={dotted}, mark=square, mark size=2, color=multcolor100, line width=1.0pt ]
        table [x=h, y=err, col sep=comma]
            {\multHundred};
    \end{cvgh}
}

\newcommand{\cvgFEMCorrMultOnedeg}[6]{
    \cvgFEMCorrMultOnedegGeneral{#1}{#2}{#3}{#4}{#5}{2,3}{#6}{$L^2$}
}

\newcommand{\cvgFEMCorrMultOnedegHun}[6]{
    \cvgFEMCorrMultOnedegGeneral{#1}{#2}{#3}{#4}{#5}{1,2}{#6}{Semi-$H^1$}
}

\newcommand{\cvgFEMCorrMultSWOnedeg}[6]{
    \edef\fem{#1}
    \edef\femsec{#2}
    \edef\add{#3}
    \edef\mults{#4}
    \edef\multw{#5}

    \begin{cvgh}{\fem}{3}{2,3}{#6}{$L^2$}
        \addlegendentry{\,FEM $\mathbb{P}_1$\;}
        \addlegendentry{\,Mult strong $\mathbb{P}_1$\;}
        \addlegendentry{\,Add $\mathbb{P}_1$\;}
        \addlegendentry{\,FEM $\mathbb{P}_2$\;}
        \addlegendentry{\,Mult weak $\mathbb{P}_1$\;}

        \addplot [style={solid}, mark=square*, mark size=2, color=femcolor, line width=0.8pt ]
        table [x=h, y=err, col sep=comma]
            {\fem};

        \addplot [style={dotted}, mark=square*, mark size=2, color=multcolor0strong, line width=0.8pt ]
        table [x=h, y=err, col sep=comma]
            {\mults};

        \addplot [style={dashed}, mark=square*, mark size=2, color=addcolor, line width=0.8pt ]
        table [x=h, y=err, col sep=comma]
            {\add};

        \addplot [style={solid}, mark=*, mark size=2, color=femcolor, line width=0.8pt ]
        table [x=h, y=err, col sep=comma]
            {\femsec};

        \addplot [style={dotted}, mark=square*, mark size=2, color=multcolor0weak, line width=0.8pt ]
        table [x=h, y=err, col sep=comma]
            {\multw};
    \end{cvgh}
}

\newcommand{\cvgFEMCorrTwoPriors}[6]{
    \edef\degree{#1}
    \edef\fem{#2}
    \edef\add{#3}
    \edef\addDeux{#4}
    \edef\slope{\number\numexpr\degree + 1\relax}
    \ifnum\degree=1
        \def\marker{square*}
    \else\ifnum\degree=2
        \def\marker{*}
    \else\ifnum\degree=3
        \def\marker{triangle*}
    \fi\fi\fi

    \begin{cvgh}{\fem}{3}{\slope}{#5}{$L^2$}
        \addlegendentry{\,FEM $\mathbb{P}_\degree$\;}
        \addlegendentry{\,Add $\mathbb{P}_\degree$ ($u_\theta$)\;}
        \addlegendentry{\,Add $\mathbb{P}_\degree$ (#6)\;}

        \addplot [style={solid}, mark=\marker, mark size=2, color=femcolor, line width=0.8pt ]
            table [x=h, y=P\degree, col sep=comma]
                {\fem};

        \addplot [style={dashed}, mark=\marker, mark size=2, color=addcolor, line width=0.8pt ]
            table [x=h, y=P\degree, col sep=comma]
                {\add};

        \addplot [style={dashed}, mark=\marker, mark size=2, color=add2color, line width=0.8pt ]
            table [x=h, y=P\degree, col sep=comma]
                {\addDeux};
    \end{cvgh}
}

\newcommand{\cvgFEMCorrThreePriors}[8]{
    \edef\degree{#1}
    \edef\fem{#2}
    \edef\add{#3}
    \edef\addDeux{#4}
    \edef\addTrois{#5}
    \edef\slope{\number\numexpr\degree + 1\relax}
    \ifnum\degree=1
        \def\marker{square*}
    \else\ifnum\degree=2
        \def\marker{*}
    \else\ifnum\degree=3
        \def\marker{triangle*}
    \fi\fi\fi

    \begin{cvgh}{\fem}{2}{\slope}{#6}{$L^2$}
        \addlegendentry{\,FEM $\mathbb{P}_\degree$\;}
        \addlegendentry{\,Add $\mathbb{P}_\degree$ ($u_\theta$)\;}
        \addlegendentry{\,Add $\mathbb{P}_\degree$ (#7)\;}
        \addlegendentry{\,Add $\mathbb{P}_\degree$ (#8)\;}

        \addplot [style={solid}, mark=\marker, mark size=2, color=femcolor, line width=0.8pt ]
            table [x=h, y=P\degree, col sep=comma]
                {\fem};

        \addplot [style={dashed}, mark=\marker, mark size=2, color=addcolor, line width=0.8pt ]
            table [x=h, y=P\degree, col sep=comma]
                {\add};

        \addplot [style={dashed}, mark=\marker, mark size=2, color=add3color, line width=0.8pt ]
            table [x=h, y=P\degree, col sep=comma]
                {\addDeux};

        \addplot [style={dashed}, mark=\marker, mark size=2, color=add2color, line width=0.8pt ]
            table [x=h, y=P\degree, col sep=comma]
                {\addTrois};
    \end{cvgh}
}

\newcommand{\cvgFEMCorrOnedeg}[3]{
    \edef\fem{#1}
    \edef\add{#2}

    \begin{cvgh}[xticklabel style={rotate=20, anchor=east, yshift=-1pt}]{\fem}{3}{2}{#3}{$L^2$}
        \addlegendentry{\,FEM $\mathbb{P}_1$\;}
        \addlegendentry{\,Add $\mathbb{P}_1$\;}

        \addplot [style={solid}, mark=square*, mark size=2, color=femcolor, line width=0.8pt ]
        table [x=h, y=err, col sep=comma]
            {\fem};
        
        \addplot [style={dashed}, mark=square*, mark size=2, color=addcolor, line width=0.8pt ]
        table [x=h, y=err, col sep=comma]
            {\add};
    \end{cvgh}
}
\usepackage{booktabs}

\newcommand{\GainsTableOnedeg}[1]{
    \pgfplotstabletypeset[
        col sep=comma,
        every head row/.style={
        before row={\toprule[1.pt]
        & & \multicolumn{4}{c}{\textbf{Gains in $L^2$ rel error}} &
		\multicolumn{4}{c}{\textbf{Gains in $L^2$ rel error}} \\
		& & \multicolumn{4}{c}{\textbf{of our method w.r.t. PINN}} &
		\multicolumn{4}{c}{\textbf{of our method w.r.t. FEM}} \\
		\cmidrule(lr){3-6} \cmidrule(lr){7-10}
        },
        after row=\cmidrule(lr){1-1} \cmidrule(lr){2-2} \cmidrule(lr){3-6} \cmidrule(lr){7-10}},
        every last row/.style={after row=\bottomrule[1.pt]},
        every nth row={2}{before row=\cmidrule(lr){1-1} \cmidrule(lr){2-2} \cmidrule(lr){3-6} \cmidrule(lr){7-10}},
		columns/method/.style={column name=\textbf{method},string type},
        columns/N/.style={column name=\textbf{N}},
		columns/min_PINNs/.style={column name=\textbf{min},fixed},
        columns/max_PINNs/.style={column name=\textbf{max},fixed},
        columns/mean_PINNs/.style={column name=\textbf{mean},fixed},
		columns/std_PINNs/.style={column name=\textbf{std},fixed},
        columns/min_FEM/.style={column name=\textbf{min},fixed},
        columns/max_FEM/.style={column name=\textbf{max},fixed},
        columns/mean_FEM/.style={column name=\textbf{mean},fixed},
		columns/std_FEM/.style={column name=\textbf{std},fixed},
        columns={method,N,min_PINNs,max_PINNs,mean_PINNs,std_PINNs,min_FEM,max_FEM,mean_FEM,std_FEM},
        precision=2
    ]{#1}
}

\newcommand{\GainsTableOnedegData}[1]{
    \pgfplotstabletypeset[
        col sep=comma,
        every head row/.style={
        before row={\toprule[1.pt]
        & & \multicolumn{4}{c}{\textbf{Gains in $L^2$ rel error}} &
		\multicolumn{4}{c}{\textbf{Gains in $L^2$ rel error}} \\
		& & \multicolumn{4}{c}{\textbf{of our method w.r.t. Data Network}} &
		\multicolumn{4}{c}{\textbf{of our method w.r.t. FEM}} \\
		\cmidrule(lr){3-6} \cmidrule(lr){7-10}
        },
        after row=\cmidrule(lr){1-1} \cmidrule(lr){2-2} \cmidrule(lr){3-6} \cmidrule(lr){7-10}},
        every last row/.style={after row=\bottomrule[1.pt]},
        every nth row={2}{before row=\cmidrule(lr){1-1} \cmidrule(lr){2-2} \cmidrule(lr){3-6} \cmidrule(lr){7-10}},
		columns/method/.style={column name=\textbf{method},string type},
        columns/N/.style={column name=\textbf{N}},
		columns/min_PINNs/.style={column name=\textbf{min},fixed},
        columns/max_PINNs/.style={column name=\textbf{max},fixed},
        columns/mean_PINNs/.style={column name=\textbf{mean},fixed},
		columns/std_PINNs/.style={column name=\textbf{std},fixed},
        columns/min_FEM/.style={column name=\textbf{min},fixed},
        columns/max_FEM/.style={column name=\textbf{max},fixed},
        columns/mean_FEM/.style={column name=\textbf{mean},fixed},
		columns/std_FEM/.style={column name=\textbf{std},fixed},
        columns={method,N,min_PINNs,max_PINNs,mean_PINNs,std_PINNs,min_FEM,max_FEM,mean_FEM,std_FEM},
        precision=2
    ]{#1}
}

\newcommand{\GainsTableAlldeg}[1]{
    \pgfplotstabletypeset[
        col sep=comma,
        every head row/.style={
        before row={\toprule[1.pt]
        & & \multicolumn{4}{c}{\textbf{Gains in $L^2$ rel error}} &
		\multicolumn{4}{c}{\textbf{Gains in $L^2$ rel error}} \\
		& & \multicolumn{4}{c}{\textbf{of our method w.r.t. PINN}} &
		\multicolumn{4}{c}{\textbf{of our method w.r.t. FEM}} \\
		\cmidrule(lr){3-6} \cmidrule(lr){7-10}
        },
        after row=\cmidrule(lr){1-1} \cmidrule(lr){2-2} \cmidrule(lr){3-6} \cmidrule(lr){7-10}},
        every last row/.style={after row=\bottomrule[1.pt]},
        every nth row={2}{before row=\cmidrule(lr){1-1} \cmidrule(lr){2-2} \cmidrule(lr){3-6} \cmidrule(lr){7-10}},
		columns/q/.style={column name=\textbf{k}},
        columns/N/.style={column name=\textbf{N}},
		columns/min_PINNs/.style={column name=\textbf{min},fixed},
        columns/max_PINNs/.style={column name=\textbf{max},fixed},
        columns/mean_PINNs/.style={column name=\textbf{mean},fixed},
		columns/std_PINNs/.style={column name=\textbf{std},fixed},
        columns/min_FEM/.style={column name=\textbf{min},fixed},
        columns/max_FEM/.style={column name=\textbf{max},fixed},
        columns/mean_FEM/.style={column name=\textbf{mean},fixed},
		columns/std_FEM/.style={column name=\textbf{std},fixed},
        columns={q,N,min_PINNs,max_PINNs,mean_PINNs,std_PINNs,min_FEM,max_FEM,mean_FEM,std_FEM},
        precision=2
    ]{#1}
}

\newcommand{\GainsTableAlldegh}[1]{
    \pgfplotstabletypeset[
        col sep=comma,
        every head row/.style={
        before row={\toprule[1.pt]
        & & \multicolumn{4}{c}{\textbf{Gains in $L^2$ rel error}} &
		\multicolumn{4}{c}{\textbf{Gains in $L^2$ rel error}} \\
		& & \multicolumn{4}{c}{\textbf{of our method w.r.t. PINN}} &
		\multicolumn{4}{c}{\textbf{of our method w.r.t. FEM}} \\
		\cmidrule(lr){3-6} \cmidrule(lr){7-10}
        },
        after row=\cmidrule(lr){1-1} \cmidrule(lr){2-2} \cmidrule(lr){3-6} \cmidrule(lr){7-10}},
        every last row/.style={after row=\bottomrule[1.pt]},
        every nth row={2}{before row=\cmidrule(lr){1-1} \cmidrule(lr){2-2} \cmidrule(lr){3-6} \cmidrule(lr){7-10}},
		columns/q/.style={column name=\textbf{k}},
        columns/h/.style={column name=\textbf{h},sci},
		columns/min_PINNs/.style={column name=\textbf{min},fixed},
        columns/max_PINNs/.style={column name=\textbf{max},fixed},
        columns/mean_PINNs/.style={column name=\textbf{mean},fixed},
		columns/std_PINNs/.style={column name=\textbf{std},fixed},
        columns/min_FEM/.style={column name=\textbf{min},fixed},
        columns/max_FEM/.style={column name=\textbf{max},fixed},
        columns/mean_FEM/.style={column name=\textbf{mean},fixed},
		columns/std_FEM/.style={column name=\textbf{std},fixed},
        columns={q,h,min_PINNs,max_PINNs,mean_PINNs,std_PINNs,min_FEM,max_FEM,mean_FEM,std_FEM},
        precision=2
    ]{#1}
}
\usepackage{booktabs}

\newcommand{\GainsFixedMu}[3]{
    \edef\whichmu{#1}

	\pgfplotstabletypeset[
        col sep=comma,
        every head row/.style={
        before row={\toprule[1.pt]
        & \textbf{FEM} \\
		\cmidrule(lr){2-2}
        },
        after row=\cmidrule(lr){1-1} \cmidrule(lr){2-2}},
        every last row/.style={after row=\bottomrule[1.pt]},
        every nth row={2}{before row=\cmidrule(lr){1-1} \cmidrule(lr){2-2}},
        columns/N/.style={column name=\textbf{N}},
		columns/err\whichmu/.style={column name=\textbf{error},sci},
        columns={N,err\whichmu},
        precision=2
    ]{#2} \hspace{20pt}
	\pgfplotstabletypeset[
        col sep=comma,
        every head row/.style={
        before row={\toprule[1.pt]
        & & \multicolumn{2}{c}{\textbf{PINN prior $u_{\theta}$}} \\
		\cmidrule(lr){3-4}
        },
        after row=\cmidrule(lr){1-1} \cmidrule(lr){2-2} \cmidrule(lr){3-4}},
        every last row/.style={after row=\bottomrule[1.pt]},
        every nth row={2}{before row=\cmidrule(lr){1-1} \cmidrule(lr){2-2} \cmidrule(lr){3-4}},
		columns/method/.style={column name=\textbf{method},string type},
        columns/N/.style={column name=\textbf{N}},
		columns/err\whichmu/.style={column name=\textbf{error},sci},
        columns/gains\whichmu/.style={column name=\textbf{gain},fixed},
        columns={method,N,err\whichmu,gains\whichmu},
        precision=2
    ]{#3}
}

\newcommand{\GainsFixedMuTwoPriors}[2]{
    \edef\whichmu{1}

	\pgfplotstabletypeset[
        col sep=comma,
        every head row/.style={
        before row={\toprule[1.pt]
        & \textbf{FEM} \\
		\cmidrule(lr){2-2}
        },
        after row=\cmidrule(lr){1-1} \cmidrule(lr){2-2}},
        every last row/.style={after row=\bottomrule[1.pt]},
        every nth row={2}{before row=\cmidrule(lr){1-1} \cmidrule(lr){2-2}},
        columns/N/.style={column name=\textbf{N}},
		columns/err/.style={column name=\textbf{error},sci},
        columns={N,err},
        precision=2
    ]{#1} \hspace{20pt}
	\pgfplotstabletypeset[
        col sep=comma,
        every head row/.style={
        before row={\toprule[1.pt]
        & & \multicolumn{2}{c}{\textbf{PINN prior $u_{\theta}$}} &
		\multicolumn{2}{c}{\textbf{Data prior $u_{\theta}^\text{data}$}} \\
		\cmidrule(lr){3-4} \cmidrule(lr){5-6}
        },
        after row=\cmidrule(lr){1-1} \cmidrule(lr){2-2} \cmidrule(lr){3-4} \cmidrule(lr){5-6}},
        every last row/.style={after row=\bottomrule[1.pt]},
        every nth row={2}{before row=\cmidrule(lr){1-1} \cmidrule(lr){2-2} \cmidrule(lr){3-4} \cmidrule(lr){5-6}},
		columns/method/.style={column name=\textbf{method},string type},
        columns/N/.style={column name=\textbf{N}},
		columns/PINNs_err/.style={column name=\textbf{error},sci},
        columns/PINNs_gains/.style={column name=\textbf{gain},fixed},
        columns/NN_err/.style={column name=\textbf{error},sci},
		columns/NN_gains/.style={column name=\textbf{gain},fixed},
        columns={method,N,PINNs_err,PINNs_gains,NN_err,NN_gains},
        precision=2
    ]{#2}
}
\usepackage{booktabs}

\newcommand{\coststableallq}[1]{
    \pgfplotstabletypeset[
        col sep=comma,
        every head row/.style={
        before row={\toprule[1.pt]
        & & \multicolumn{2}{c}{\textbf{$N$}} &
		\multicolumn{2}{c}{\textbf{$N_\text{dofs}$}} \\
		\cmidrule(lr){3-4} \cmidrule(lr){5-6}
        },
        after row=\cmidrule(lr){1-1} \cmidrule(lr){2-2} \cmidrule(lr){3-4} \cmidrule(lr){5-6}},
        every last row/.style={after row=\bottomrule[1.pt]},
        every nth row={2}{before row=\cmidrule(lr){1-1} \cmidrule(lr){2-2} \cmidrule(lr){3-4} \cmidrule(lr){5-6}},
		columns/q/.style={column name=\textbf{k}},
        columns/e/.style={column name=\textbf{e},sci},
		columns/FEM_N/.style={column name=\textbf{FEM},fixed},
        columns/Add_N/.style={column name=\textbf{Add},fixed},
        columns/FEM_dofs/.style={column name=\textbf{FEM},fixed},
		columns/Add_dofs/.style={column name=\textbf{Add},fixed},
        columns={q,e,FEM_N,Add_N,FEM_dofs,Add_dofs},
        precision=2
    ]{#1}
}

\newcommand{\coststableallqhundred}[1]{
    \pgfplotstabletypeset[
        col sep=comma,
        every head row/.style={
        before row={\toprule[1.pt]
        & & \multicolumn{2}{c}{\textbf{$n_p=1$}} &
		\multicolumn{2}{c}{\textbf{$n_p=100$}} \\
		\cmidrule(lr){3-4} \cmidrule(lr){5-6}
        },
        after row=\cmidrule(lr){1-1} \cmidrule(lr){2-2} \cmidrule(lr){3-4} \cmidrule(lr){5-6}},
        every last row/.style={after row=\bottomrule[1.pt]},
        every nth row={2}{before row=\cmidrule(lr){1-1} \cmidrule(lr){2-2} \cmidrule(lr){3-4} \cmidrule(lr){5-6}},
		columns/q/.style={column name=\textbf{k}},
        columns/e/.style={column name=\textbf{e},sci},
		columns/FEM1/.style={column name=\textbf{FEM},fixed},
        columns/Corr1/.style={column name=\textbf{Add},fixed},
        columns/FEM100/.style={column name=\textbf{FEM},fixed},
		columns/Corr100/.style={column name=\textbf{Add},fixed},
        columns={q,e,FEM1,Corr1,FEM100,Corr100},
        precision=2
    ]{#1}
}
\usepackage{booktabs}
\usepackage{amssymb}
\usepackage{mathtools}
\usepackage{pgfplots}
\usepackage{pgfplotstable}
\usepackage{datatool}
\usepackage{fp}

\pgfplotsset{
    compat=newest,
}
\pgfplotsset{
    smaller labels/.style={
        label style={font=\footnotesize},
        tick label style={font=\footnotesize}
    }
}
\tikzset{font=\small}
\usetikzlibrary{
    fpu,
    fixedpointarithmetic,
    babel,
    external,
    arrows.meta,
    plotmarks,
    positioning,
    angles,
    quotes,
    intersections,
    calc,
    spy,
    decorations.pathreplacing,
    matrix,
    fit,
}
\usepgfplotslibrary{fillbetween}

\definecolor{femcolor}{RGB}{51, 138, 55} 
\definecolor{addcolor}{RGB}{217,95,2} 
\definecolor{addsobcolor}{RGB}{199,39,34} 
\definecolor{multcolor3}{RGB}{117,112,179} 
\definecolor{multcolor100}{RGB}{0,0,0} 
\definecolor{multcolor0weak}{RGB}{49, 73, 181} 
\definecolor{multcolor0strong}{RGB}{49, 181, 161} 

\newcommand{\cvgtimeerror}[2]{
    \begin{tikzpicture}
        \edef\fem{#1}
        \edef\add{#2}

        \begin{loglogaxis}[
            smaller labels,
            name = left_plot,
            axis lines = left,
            enlarge x limits={abs=5pt},
            enlarge y limits={abs=5pt},
            yminorticks=false,
            xminorticks=false,
			xmode=log,
            xlabel = {Time (s)},
            ylabel = {\rotatebox{270}{$L^2$}},
            xlabel style={at={(ticklabel* cs:1.01)},anchor=north},
            ylabel style={at={(ticklabel* cs:1.01)},anchor=west},
            width=\plotwidth, height=\plotheight,
            mark options={solid, scale=1},
            grid = major,
            legend columns=1,
            legend to name=leg:legendFEMCORR_\iterator,
            legend image post style={mark options={solid, scale=1}},
        ]

            \addlegendentry{\,FEM $\mathbb{P}_1$\;}
            \addlegendentry{\,Add $\mathbb{P}_1$\;}

            \addplot [style={solid}, mark=square*, mark size=2, color=femcolor, line width=0.8pt ]
            table [x=time, y=err, col sep=comma]
                {\fem};
            
            \addplot [style={dashed}, mark=square*, mark size=2, color=addcolor, line width=0.8pt ]
            table [x=time, y=err, col sep=comma]
                {\add};

        \end{loglogaxis}
        
        \node[yshift=-10pt,xshift=20pt] at (left_plot.outer north) {\pgfplotslegendfromname{leg:legendFEMCORR_\iterator}};

    \end{tikzpicture}
}

\newcommand{\cvgtimeerrorparam}[1]{
    \begin{tikzpicture}
        \edef\tps{#1}
        \begin{loglogaxis}[
            smaller labels,
            name = left_plot,
            axis lines = left,
            enlarge x limits={abs=5pt},
            enlarge y limits={lower,abs=5pt},
			ymode=normal,
            xmode=normal,
            xlabel = {$n_p$},
            ylabel = {\rotatebox{270}{Time (s)}},
            xlabel style={at={(ticklabel* cs:1.01)},anchor=west},
            ylabel style={at={(ticklabel* cs:1.01)},anchor=west},
            width=\plotwidth, height=\plotheight,
            mark options={solid, scale=1},
            grid = major,
            legend columns=1,
            legend to name=leg:legendFEMCORR_\iterator,
            legend image post style={mark options={solid, scale=1}},
            ymin=0, ymax=1860,
            xmin=0, xmax=50,
        ]

            \addlegendentry{\,FEM $\mathbb{P}_1$\;}
            \addlegendentry{\,Add $\mathbb{P}_1$\;}

            \addplot [style={solid}, mark size=2, color=femcolor, line width=1pt ]
            table [x=np, y=tfem, col sep=comma]
                {\tps};
            
            \addplot [style={dashed}, mark size=2, color=addcolor, line width=1pt ]
            table [x=np, y=tadd, col sep=comma]
                {\tps};
                
            \addlegendimage{
                /pgfplots/legend image code/.code={
                    \draw[femcolor,very thick] (0, -0.1) -- (0, 0.15);
                    \draw[femcolor,very thick] (-0.02, 0.15) -- (0.55, 0.15);
                    \draw[addcolor,very thick] (-0.02, -0.1) -- (0.57, -0.1);
                    \draw[addcolor,very thick] (0.55, -0.1) -- (0.55, 0.17);
                },
            }
            \addlegendentry{Offline cost}

            \addplot [style={dotted}, line width=1.5pt, gray, clip=false]
            coordinates {(19,-200) (19,1900)};
        \end{loglogaxis}
        
        \node[xshift=-22pt, yshift=6pt, gray] at (left_plot.outer south) {\normalsize\bfseries $19$};
        \node[draw=femcolor, thick, rectangle, xshift=16pt, yshift=-36pt, femcolor, inner sep=1.2pt] at (left_plot.outer west) {\normalsize\bfseries $40$};
        \node[draw=addcolor, thick, rectangle, xshift=15pt, yshift=-10pt, addcolor, inner sep=1.2pt] at (left_plot.outer west) {\normalsize\bfseries $708$};
    \end{tikzpicture}
}

\newcommand{\cvgtimeerrorparamD}[1]{
    \begin{tikzpicture}
        \edef\tps{#1}
        \begin{loglogaxis}[
            smaller labels,
            name = left_plot,
            axis lines = left,
            enlarge x limits={abs=5pt},
            enlarge y limits={lower,abs=5pt},
			ymode=normal,
            xmode=normal,
            xlabel = {$n_p$},
            ylabel = {\rotatebox{270}{Time (s)}},
            xlabel style={at={(ticklabel* cs:1.01)},anchor=west},
            ylabel style={at={(ticklabel* cs:1.01)},anchor=west},
            width=\plotwidth, height=\plotheight,
            mark options={solid, scale=1},
            grid = major,
            legend columns=1,
            legend to name=leg:legendFEMCORR_\iterator,
            legend image post style={mark options={solid, scale=1}},
            ymin=0, ymax=6500,
            xmin=0, xmax=50,
        ]

            \addlegendentry{\,FEM $\mathbb{P}_1$\;}
            \addlegendentry{\,Add $\mathbb{P}_1$\;}

            \addplot [style={solid}, mark size=2, color=femcolor, line width=1pt ]
            table [x=np, y=tfem, col sep=comma]
                {\tps};
            
            \addplot [style={dashed}, mark size=2, color=addcolor, line width=1pt ]
            table [x=np, y=tadd, col sep=comma]
                {\tps};
                
            \addlegendimage{
                /pgfplots/legend image code/.code={
                    \draw[femcolor,very thick] (0, -0.1) -- (0, 0.15);
                    \draw[femcolor,very thick] (-0.02, 0.15) -- (0.55, 0.15);
                    \draw[addcolor,very thick] (-0.02, -0.1) -- (0.57, -0.1);
                    \draw[addcolor,very thick] (0.55, -0.1) -- (0.55, 0.17);
                },
            }
            \addlegendentry{Offline cost}

            \addplot [style={dotted}, line width=1.5pt, gray, clip=false]
            coordinates {(5,-200) (5,6500)};
        \end{loglogaxis}
        
        \node[yshift=-10pt, xshift=-10pt] at (left_plot.outer north east) {\pgfplotslegendfromname{leg:legendFEMCORR_\iterator}};

        \node[xshift=-60pt, yshift=6pt, gray] at (left_plot.outer south) {\normalsize\bfseries $5$};
        \node[draw=femcolor, thick, rectangle, xshift=16pt, yshift=-38pt, femcolor, inner sep=1.2pt] at (left_plot.outer west) {\normalsize\bfseries $98$};
        \node[draw=addcolor, thick, rectangle, xshift=15pt, yshift=-27pt, addcolor, inner sep=1.2pt] at (left_plot.outer west) {\normalsize\bfseries $708$};
    \end{tikzpicture}
}

\pgfplotstableset{
    highlight cell/.style={
        every row #1 column \pgfplotstablecol/.style={
            postproc cell content/.append code={
                \pgfkeysalso{@cell content/.add={(}{)}}%
            }
        }
    }
}

\newcommand{\timerequired}[1]{
    \pgfplotstabletypeset[
        col sep=comma,
        every head row/.style={
        before row={\toprule[1.pt]
        & \multicolumn{2}{c}{\textbf{$N$}} &
		\multicolumn{2}{c}{\textbf{$N_\text{dofs}$}} &
		\multicolumn{2}{c}{\textbf{Computation time}} \\
		\cmidrule(lr){2-3} \cmidrule(lr){4-5} \cmidrule(lr){6-7}
        },
        after row=\cmidrule(lr){1-1} \cmidrule(lr){2-3} \cmidrule(lr){4-5} \cmidrule(lr){6-7}},
        every last row/.style={after row=\bottomrule[1.pt]},
        columns/e/.style={column name=\textbf{e},sci},
		columns/FEM_N/.style={column name=\textbf{FEM},fixed},
        columns/ADD_N/.style={column name=\textbf{Add},fixed},
        columns/FEM_dofs/.style={column name=\textbf{FEM},sci},
		columns/ADD_dofs/.style={column name=\textbf{Add},sci},
        columns/FEM_time/.style={column name=\textbf{FEM},sci,highlight cell=0,highlight cell=1,highlight cell=2},
        columns/ADD_time/.style={column name=\textbf{Add},sci,highlight cell=2},
        columns={e,FEM_N,ADD_N,FEM_dofs,ADD_dofs,FEM_time,ADD_time},
        precision=2,
    ]{#1}
}

\newtheorem{thrm}{Theorem}
\newtheorem{rmrk}[thrm]{Remark}
\crefname{paragraph}{Section}{Sections}
\Crefname{paragraph}{Section}{Sections}
\crefname{rmrk}{Remark}{Remarks}
\Crefname{rmrk}{Remark}{Remarks}
\crefname{thrm}{Theorem}{Theorems}
\Crefname{thrm}{Theorem}{Theorems}

\DeclareMathOperator*{\argmin}{argmin}

\newcommand{\boldparagraph}[1]{\medskip\paragraph{\RTwo{\textbf{#1.}}}}

\numberwithin{equation}{section}

\newcommand{\ROne}[1]{#1}
\newcommand{\RTwo}[1]{#1}
\newcommand{\RBoth}[1]{#1}
\newcommand{\ROneN}[1]{#1}
\newcommand{\RTwoN}[1]{#1}
\newcommand{\RBothN}[1]{#1}

\usepackage{authblk}
\newcommand{\address}{\affil}

\title{Enriching continuous Lagrange finite element approximation spaces using neural networks}

\author[1]{Hélène Barucq}
\author[2]{Michel Duprez}
\author[1]{Florian Faucher}
\author[3]{Emmanuel Franck}
\author[2]{Frédérique Lecourtier}
\author[4]{Vanessa Lleras}
\author[3]{Victor Michel-Dansac}
\author[1]{Nicolas Victorion}

\address[1]{Project-Team Makutu, Inria, University of Pau and Pays de l'Adour, TotalEnergies, CNRS UMR 5142, Pau, France}
\address[2]{Université de Strasbourg, CNRS, Inria, ICube, F-67000, Strasbourg, France}
\address[3]{Université de Strasbourg, CNRS, Inria, IRMA, F-67000, Strasbourg, France}
\address[4]{IMAG, University of Montpellier, CNRS UMR 5149, Montpellier, France}

\date{\today}

\begin{document}
    \pgfkeys{/pgf/number format/.cd,1000 sep={\,}}%

    \maketitle
    
    \begin{abstract}
        In this work, we present a study combining two approaches in the context of solving PDEs: the continuous finite element method (FEM) and more recent techniques based on neural networks. In recent years, physics-informed neural networks (PINNs) have become particularly interesting for rapidly solving PDEs, especially in high dimensions.
        However, their lack of accuracy can be a significant drawback in this context, hence the interest in combining them with FEM, for which error estimates are already known. The complete pipeline proposed here consists in modifying the classical FEM approximation spaces by taking information from a prior, chosen as the prediction of a neural network. On the one hand, this combination improves and certifies the prediction of neural networks, to obtain a fast and accurate solution. On the other hand, error estimates are proven, showing that such strategies outperform classical ones by a factor that depends only on the quality of the prior. We validate our approach with numerical results performed on parametric problems with \RTwo{1D, 2D and 3D} geometries. These experiments demonstrate that to achieve a given accuracy, a coarser mesh can be used with our enriched FEM compared to the standard FEM, leading to reduced computational time, particularly for parametric problems.
    \end{abstract}
    
    \tableofcontents
    
    \section*{Introduction}
    \ROne{The finite element method (FEM, e.g.,~\cite{ciarlet2002finite,Ern2004TheoryAP,brenner2008mathematical}) is widely used for computing accurate solutions of complex PDEs for which there is no analytic solution. It begins with a mesh as a subdivision of the computational domain into elements, typically simplexes or quadrilaterals/hexahedra. Then the numerical solution is constructed from degrees of freedom (dofs) that are coefficients of the finite element discretization defined by basis functions related to the elements. Standard FEM performs poorly when the solution has strong local variations (like sharp gradients near cracks, re-entrant corners, material interfaces, etc). Then refining the mesh and/or increasing the order of approximation can help, but the method becomes very (sometimes too) expensive. }

\ROne{This has motivated intensive research on new finite element methods to combine accuracy and lower computational burden. For instance, the Generalized Finite Element Method (GFEM, \cite{strouboulis2001generalized,fries2010extended}) enriches the classical finite element approximation space with special functions that capture local solution features (singularities, oscillations, boundary layers, etc). Basically, GFEM belongs to the family of Partition of Unity Methods \cite{babuvska1997partition}. The Extended Finite Element Method (XFEM) is a special case of GFEM, developed mainly for fracture mechanics. It was popularized by Belytschko and co-workers \cite{sukumar2000extended}, who suggested to enrich the finite element space with discontinuous functions and/or singularity functions, without re-meshing.
\ROneN{Therefore, both GFEM and XFEM rely on the enrichment of the finite element spaces to treat local structures in order to increase the accuracy 
of the solution, \cite{belytschko2009review,babuvska2004generalized}. 
The approach we propose has the same idea, but a different conception as 
it uses a (global) prior of the solution (on the entire domain), and the FEM 
is used to correct it; this prior information further allows the use of a low-order FEM. 
Furthermore, our application focuses on a parametric problem, with a prior learnt from a set of parameter values. In \cite{CanLeyChiGonCueFeuBerHue2016}, the authors deal with parametric problems by enriching finite element spaces additively with functions determined by Proper Generalized Decomposition (PGD). However, PGD requires a separability assumption on the parametric dependence of the solution, which is not necessary in our approach. 
More generally, XFEM/GFEM is sensitive to the choice of functions 
for enrichment, which must be chosen problem-dependent for best 
efficiency, while NN-based priors are more adapted to parametric problems.} Another way to improve the accuracy of FEM solutions is through post-processing
superconvergence (e.g., \cite{wahlbin2006superconvergence} and the references therein).
However, their effectiveness is often problem-specific and depends on particular properties
of the PDE, mesh, and element type. Moreover, they are not particularly suited for
parametric problems.}

\ROne{One approach to lowering the cost further consists of reducing the size of the discrete system. This can be done by decreasing the number of degrees of freedom as in the Trefftz method (e.g.,~\cite{hiptmair2013error,moiola2018space,ImbMoiSto2022})
    or the hybridizable discontinuous Galerkin (HDG) method (e.g.,~\cite{Cockburn2008,hungria2017hdg,Pham2024stabilization}), along with hybrid high-order (HHO) approximations \cite{ern2024convergence}, which solve a global system associated only with the degrees of freedom on the skeleton of the mesh. It is then worth noting that these approaches are much easier to implement in the stationary case.
    In the same vein, we can mention reduced-order methods, which focus more on the number of basis functions used to construct the approximate solution. The idea is to construct a low-dimensional subspace from high-fidelity simulations and then solve the PDE projected onto this subspace. This approach is very interesting when the parametric family or the dominant structure of the solutions is known. Examples include the snapshot method or POD (proper orthogonal decomposition, see~\cite{berkooz1993proper, ravindran2000reduced}) and the reduced basis method \cite{patera2007reduced,prud2002reliable}.
    Avoiding mesh refinement is a significant advantage, making it possible to control often exorbitant computational costs without any special preprocessing.
    This justifies other approaches that do not use meshes. Meshless methods have been investigated in the last decades, and isogeometric analysis has been employed, see e.g.~\cite{bazilevs2006isogeometric, hughes2005isogeometric,Frambati2022practical}.}

In recent years, learning-based alternatives have emerged, such as Physics-Informed Neural Networks
(PINNs, e.g.,~\cite{RAISSI2019686})
or the Deep Ritz method~\cite{e2017deepritzmethoddeep}.
\RTwo{Similar methods were already present in the 1990s, see for instance~\cite{LeeKan1990,MeaFer1994,LagLikFot1998}, but the advances in deep learning and computational power have significantly boosted their applicability in recent years.}
\ROne{The idea is to approximate the solution of the PDE under consideration using a neural network. If this solution lives in some function space $H$, this amounts to projecting it onto a finite-dimensional subset of $H$ (e.g., a submanifold) defined by the parametrization inherent to the neural network. We stress that, conversely, FEM projects the solution onto a finite-dimensional linear subspace of $H$. The neural network is then trained by minimizing a loss function, taking the underlying physics into account.}
Unlike neural networks trained with more conventional data-driven loss functions, these methods share similarities with traditional solvers: they require the same inputs, namely the PDE, physical parameters, boundary, and initial conditions.
In addition, the training phase requires approximating the PDE solution in a discrete set of points of the space domain.
Thus, these approaches have some advantages,
\RTwo{notably their relative dimension-insensitivity and their mesh-free quality.
    Indeed, PINNs do not use a mesh, but rather require sampling points in the domain, as it is done in meshless methods (see for instance \cite{Frambati2022practical}). But on complex geometries, sampling can be challenging; however, it is \RTwoN{very often} easier than meshing.}
\RTwo{Further refinements have been proposed to improve PINNs.
    For instance, we refer to the review paper~\cite{cuomo2022scientific} for some analysis on PINNs, to~\cite{JagKhaKar2020} for a conservative version, to~\cite{KhaZhaKar2021} for PINNs using variational formulas, or to~\cite{DeRMisMol2024} for PINNs using weak formulations to improve discontinuity handling.}
Since they do not require data (in the form of reference solutions), they are particularly well-suited to high-dimensional problems on complex domains.
However, at present, these learning-based techniques are not competitive with classical finite element methods (see~\cite{grossmann2023can}), mainly because network-based methods lack precision and convergence guarantees, see~\cite{sikora2024comparison} for a comparison between PINNs and FEM.
While FEM has a better error/computation time ratio for a single resolution, PINNs are more advantageous for parametric systems where a multitude of resolutions is needed.

This paper aims to propose a new method that combines learning-based and finite element methods \RTwo{applied on coarse meshes.
    The general idea is to assume that some so-called ``prior'' information
    is available about the solution of the PDE.
    It can be thought of as a function that approximates the solution,
    or a family of functions that approximate the solution for different parameters.
    This prior then modifies the finite element approximation space,
    which is finally used to compute the solution of the PDE.}
\RTwo{More precisely, we use a PINN to compute
    either one offline solution in the non-parametric case,
    or a parametric family of offline solutions in the parametric case.
    This is followed by calculating an online solution using coarse finite elements,
    with the PINN solution used to modify the finite element approximation space;
    in the parametric case, this is done for a single parameter.}
The result is a method capable of rapidly predicting a PDE solution while guaranteeing convergence properties, thanks to the FEM framework.
Finite element resolution improves the prediction while remaining cheap as it is performed on a coarse mesh, benefiting from the network prediction.
This paper proposes two ways to enrich the FEM.
In both cases, the finite element error will be exhibited as a function of the network error with respect to the true solution.
These corrections will be called ``additive'' and ``multiplicative'' depending on how the prediction is incorporated in the FEM spaces. \ROneN{The error estimates, resulting from these two enriched approaches, are mainly meant to show that enriching the space lets us get the same convergence orders as the standard approach. In addition, this study helps guide choices about prior learning by showing that the higher the approximation orders, the more high-order derivatives are involved. The manuscript does not aim to study the gain constants obtained in these error estimates, which depend on the (hard to quantify) prediction error associated with neural networks.}

\RTwo{Work has already been undertaken to combine numerical methods and the use of a prior, be it obtained through a neural network prediction or by some other means.
    Starting with the enhancement of FEM-related methods, in the FEM itself, the approximation space can also be enriched to ensure stability, see for instance the introduction of bubble functions in mixed problems (see e.g.~\cite{Ern2004TheoryAP}).
    In $\varphi$-FEM, developed in~\cite{duprez2020phi} (see also~\cite{duprez2023new,cotin2023phi,DupLleLozVui2023,duprez2023phi} for other contexts), the FEM prior is a level-set function used to localize the domain boundary.}
\ROneN{In~\cite{CanLeyChiGonCueFeuBerHue2016}, the authors deal with parametric problems by additively enriching finite element spaces with functions determined by Proper Generalized Decomposition (PGD). However, PGD requires a separability assumption on the parametric dependence of the solution.}
Combining FEM and neural networks, in~\cite{feng_hybrid_2024}, the authors solve a time-dependent PDE by splitting the problem into a spatial resolution (handled by a neural network) and a temporal one (handled by FEM).
Other works have also explored related ideas, such as
\cite{BadLiMar2024} where neural networks are interpolated onto FEM spaces,
\cite{MarJenLesRic2024} where a coarse FEM is augmented with fine information from a neural network,
or \cite{WanLiZha2025} where GFEM is enriched with neural networks.
It is also possible to include the FEM shape functions in PINNs, as proposed in~\cite{skardova_finite_2024,XioLonBorJia2025}.
Concerning other numerical methods, in~\cite{brunet2019physics,AghFraHilMicVig2025}, the authors initialize Newton's algorithm, used when solving nonlinear equations, with an initial guess derived from the prediction of a neural network. Such a prediction can also be used as a prior for discontinuous Galerkin methods (see~\cite{FraMicNav2024}).
In the finite difference context, the authors of~\cite{XiaPenYaoZho2025} propose to replace automatic differentiation with finite difference discretization when possible.

In this paper, we consider general parametric linear elliptic differential equations
defined on a \RTwo{domain $\Omega \subset \mathbb{R}^d$}
with $d$ space dimensions, with \RTwo{a smooth boundary $\partial\Omega$} \ROne{(for instance a Lipschitz polytope)}.
Consider a parameter space $\mathcal{M}=\{\bm{\mu}=(\mu_1,\ldots,\mu_p)\in \mathbb{R}^p\}$.
The typical problem of interest is, for one or several $\bm{\mu}\in \mathcal{M}$, to find $u: \Omega\to \mathbb{R}$ such that
\begin{equation}\label{eq:ob_pde}
    \mathcal{L}\big(u\,;\bm{x},\bm{\mu}\big) = f(\bm{x},\bm{\mu}),
\end{equation}
with $\bm{x}=(x_1,\dots,x_d)\in\Omega$ the space variable, and where $\mathcal{L}$ is the parametric differential operator defined  by
\begin{equation}\label{eq:operatorL}
    \mathcal{L}(\,\cdot\,;\bm{x},\bm{\mu}) : u \mapsto R(\bm{x},\bm{\mu}) u + C(\ROne{\bm{x},}\bm{\mu}) \cdot \nabla u - \frac{1}{\text{Pe}} \nabla \cdot (D(\bm{x},\bm{\mu}) \nabla u),
\end{equation}
with $f(\bm{x},\bm{\mu})\in L^2(\Omega)$ the source term,
$R(\bm{x},\bm{\mu})\in L^{\infty}(\Omega)
$ the reaction coefficient,
\ROne{$C(\bm{x},\bm{\mu}) \in W^{1,\infty}(\Omega)^d$} the convection coefficient,
\smash{$D(\bm{x},\bm{\mu}) \in {(W^{1,\infty}(\Omega))}^{d\times d}$}
the (symmetric and positive definite) diffusion matrix,
and $\text{Pe} \in \mathbb{R}_+^*$ the Péclet number representing the ratio between convection and diffusion.
\ROne{We suppose that $\mathcal{L}$ is coercive.}
The differential operator is considered with Dirichlet, Neumann or Robin boundary conditions, which can also depend on $\bm{\mu}$.

\begin{figure}[!ht]
    \centering
    \includegraphics[scale=0.6]{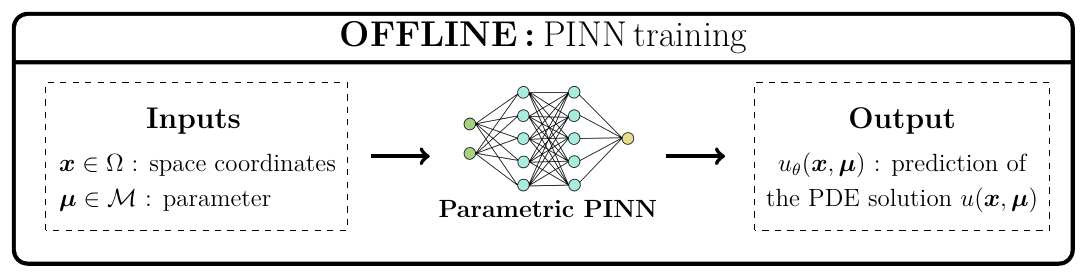}

    \includegraphics[scale=0.6]{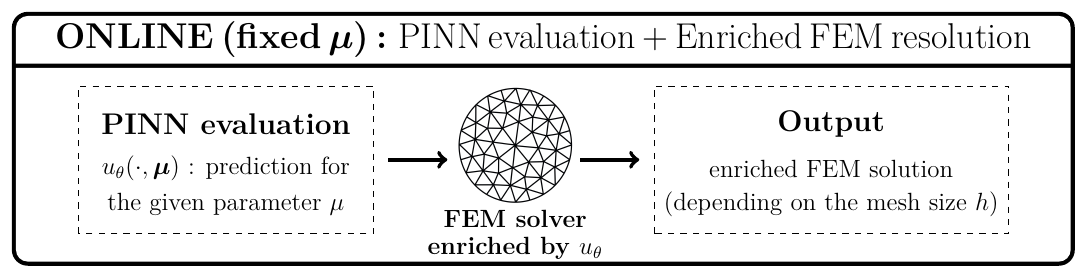}
    \caption{Pipeline of the enriched method. Top: offline phase (PINN training). Bottom: online phase (PINN evaluation + Enriched FEM resolution).}
    \label{fig:pipeline}
\end{figure}

\cref{fig:pipeline} presents the pipeline of our enriched method, described above.

The manuscript is organized as follows:
in \cref{sec:FEM}, we recall the classical finite element method applied to our problem, introducing the notations needed in the following sections.
In \cref{sec:additive_prior,sec:multiplicative_prior}, we present the two proposed enrichments and prove error estimates.
Both approaches rely on modifying the functions of the FEM approximation space, using information from prior knowledge of the solution.
This prior is first introduced in an additive way, and then in a multiplicative way.
Both approaches are compared in \cref{sec:comparison_add_mul}.
\cref{sec:prior_construction} is devoted to the construction of the prior, justifying the use of PINNs and recalling methods for improving their efficiency.
In \cref{sec:implementation_details}, we give details on the implementation.
Numerical simulations conclude this manuscript in \cref{sec:numerical_results} and show that the proposed methods can significantly reduce the computational cost of solving parametric problems.
A short conclusion, as well as plans for future work, are given in \cref{sec:conclusion}.
\cref{app:notations} reviews the notations used throughout the manuscript.

    \section{Continuous finite element method}\label{sec:FEM}
    The goal of this section is to recall the classical FEM,
and to introduce the notation that will be used throughout the paper.
\ROne{Recall that, in the online step, the goal is to perform a coarse finite element resolution of the PDE for a single parameter $\bm{\mu}$. Therefore, in this section and the next two concerning error estimates, we work a fixed $\bm{\mu}$, and we omit the dependence on $\bm{\mu}$ for conciseness.}
To solve the problem~\eqref{eq:ob_pde} under consideration for a fixed parameter $\bm{\mu}$
with homogeneous Dirichlet boundary conditions using
the continuous FEM, we rewrite it as the following variational problem:
\begin{equation}\label{eq:weakform}
	\text{Find } u \in V^0 \text{ such that, } \; \forall v\in V^0, \; a(u,v)=l(v),
\end{equation}
where $V^0=H^1_0(\Omega)$,
and where the bilinear form $a$ is given by
\begin{equation*}
	a(u,v)=
	\frac{1}{\text{Pe}} \int_{\Omega} D\RTwo{(\bm{x})} \nabla u\RTwo{(\bm{x})} \cdot  \nabla v\RTwo{(\bm{x})} \, \RTwo{\mathrm{d}\bm{x}} +
	\int_{\Omega} R\RTwo{(\bm{x})} \, u\RTwo{(\bm{x})} \, v\RTwo{(\bm{x})} \, \RTwo{\mathrm{d}\bm{x}}  +
	\int_{\Omega} v\RTwo{(\bm{x})} \, C\RTwo{(\bm{x})} \cdot \nabla u\RTwo{(\bm{x})} \, \RTwo{\mathrm{d}\bm{x}},
\end{equation*}
while the linear form $l$ reads
\begin{equation*}
	l(v)=\int_{\Omega} f\RTwo{(\bm{x})} \, v\RTwo{(\bm{x})} \, \RTwo{\mathrm{d}\bm{x}}.
\end{equation*}

\begin{rmrk}
	Note that since $a$ is continuous on $V^0\times V^0$ and coercive and $l$ is continuous on $V^0$,
	the existence and uniqueness of the solution $u$ are ensured by the Lax-Milgram \RTwo{theorem}.
\end{rmrk}

Let $\mathcal{T}_h$ be a mesh of the domain $\Omega$ composed of simplexes,
where $h$ denotes the characteristic size of the mesh, i.e.\ the biggest diameter of the simplexes. We suppose that $\mathcal{T}_h$ satisfies the Ciarlet condition (see e.g.\ \cite{Ern2004TheoryAP}) and that its boundary is exactly $\partial\Omega$.
Consider $V_h^0\subset V_h\subset V=H^1(\Omega)$  the two continuous Lagrange finite elements spaces of degree $k\geqslant 1$ defined by
\begin{equation}\label{eq:Vh}
	V_h = \left\{v_h\in C^0(\Omega),\; \forall K\in \mathcal{T}_h,\; v_h\vert_{K}\in\mathbb{P}_k\right\},
\end{equation}
and
\[
	V_h^0 = \left\{v_h\in C^0(\Omega),\; \forall K\in \mathcal{T}_h,\; v_h\vert_{K}\in\mathbb{P}_k,\; v_h\vert_{\partial\Omega}=0\right\},
\]
with $\mathbb{P}_k$ the space of polynomials with real coefficients of degree at most $k$.
The solution to~\eqref{eq:weakform} will be approximated by the solution $u_h$ to
\begin{equation}\label{eq:approachform}
	\text{Find } u_h \in V_h^0 \text{ such that, } \; \forall v_h\in V_h^0, \; a(u_h,v_h)=l(v_h).
\end{equation}

Let us now introduce some results used in the remainder of the paper.
We first define the Lagrange interpolation operator by
\begin{equation}\label{eq:Ih}
	\mathcal{I}_h  : C^0(\Omega) \ni v \mapsto \sum_{i=1}^{N_{\text{\rm dofs}}} v\big(\bm{x}^{(i)}\big) \psi_i\in V_h,
\end{equation}
with \smash{${(\bm{x}^{(i)})}_{i \in \{1,\ldots,N_{\text{\rm dofs}}\}}$}
the $N_{\text{\rm dofs}}$ degrees of freedom (dofs) associated to the mesh,
and \smash{${(\psi_i)}_{i \in \{1,\ldots,N_{\text{\rm dofs}}\}}$}
the associated Lagrange shape functions of degree $k$. \ROne{These Lagrange shape functions are the unique continuous piecewise polynomials of total degree at most $k$ satisfying the interpolation property}
\[
	\ROne{\psi_i(\bm{x}^{(j)}) = \delta_{ij}, \quad \forall i,j \in \{1,\ldots,N_{\text{\rm dofs}}\},}
\]
\ROne{where $\delta_{ij}$ is the Kronecker symbol.}

\begin{rmrk}
	In the whole manuscript, for a Sobolev space $H$, the notation $|\cdot|_{H}$ and $\|\cdot\|_{H}$ will represent respectively the semi-norm and the norm in $H$. \ROne{Namely, for $v \in H^1(\Omega)$, we set}
	\[
		\ROne{|v|_{H^1(\Omega)}^2 = \int_\Omega |\nabla v(\bm{x})|^2 \, \mathrm{d}\bm{x}
		\qquad \text{and} \qquad
		\|v\|_{H^1(\Omega)}^2 = \int_\Omega |v(\bm{x})|^2 \, \mathrm{d}\bm{x} + \int_\Omega |\nabla v(\bm{x})|^2 \, \mathrm{d}\bm{x}.}
	\]
\end{rmrk}

The following result gives a bound of the interpolation error:
\begin{thrm}[see e.g.\ \cite{Ern2004TheoryAP}]\label{th:interpol}
	There exists $C_q>0$ such that
	for all $v\in H^{q+1}(\Omega)$ and $1\leqslant q\leqslant k$,
	\begin{equation*}
		\|v-\mathcal{I}_h v\|_{H^1}\leqslant C_q h^q |v|_{H^{q+1}}.
	\end{equation*}
\end{thrm}

The next estimate is associated to the elliptic regularity:

\begin{thrm}[see e.g. {\cite[Theorem 4, p.\ 317]{evans2022partial}}]\label{th:ellip}
	\ROne{Suppose that the boundary $\partial \Omega$ is $\mathcal{C}^2$,
    $R\in L^{\infty}(\Omega)$, $C \in W^{1,\infty}(\Omega)^d$ and $D \in {(W^{1,\infty}(\Omega))}^{d\times d}$ in the definition of $\mathcal{L}$  \RTwo{given in~\eqref{eq:operatorL}}.
    There exists $C_e>0$ such that for all $\RTwo{\xi}\in L^2(\Omega)$
    and all associated weak solution $w\in H^1_0(\Omega)$ to 
	\begin{equation}\label{eq:dual}
		\mathcal{L}^* w=\xi
	\end{equation}
    with homogeneous Dirichlet boundary condition,
   we have $w\in H^2(\Omega)$
   and }
   \begin{equation*}
		\|w\|_{H^2}\leqslant C_{e} \|\xi\|_{L^2}.
	\end{equation*}
	Here $\mathcal{L}^*$ represents the adjoint of the operator $\mathcal{L}$.
\end{thrm}

\noindent These estimates, combined with Céa's Lemma, which uses the continuity and coercivity of $a$, give the following error estimate:
\begin{thrm}[see e.g.\ \cite{Ern2004TheoryAP}]\label{thm:classical_error_estimate}
	\ROne{Let $1\leqslant q\leqslant k$. Considering} $u\in H^{q+1}(\Omega)$ and $u_h\in V_h^0$ the solutions to~\eqref{eq:weakform} and~\eqref{eq:approachform}, one has
	\begin{equation*}
		|u-u_h|_{H^1}\leqslant C_q\dfrac{\gamma}{\alpha}h^{q} |u|_{H^{q+1}}
	\end{equation*}
	and
	\begin{equation*}
		\|u-u_h\|_{L^2}\leqslant C_eC_1C_q\dfrac{\gamma^2}{\alpha}h^{q+1} |u|_{H^{q+1}},
	\end{equation*}
	where $\gamma$ and $\alpha$ are respectively
	the constants of continuity and coercivity of $a$.
\end{thrm}

For the sake of simplicity, we consider an elliptic boundary
value problem with homogeneous Dirichlet conditions.
Obviously, we can use more general boundary conditions as Robin-like
conditions depending on the parameters. This will be investigated
in the numerical experiments to show that the proposed methodology
applies to a larger class of boundary value problems, see \cref{sec:Lap2DMixRing}.

    \section{Enriching the finite element method with additive priors}\label{sec:additive_prior}
    In this section, we assume that a prior knowledge of the solution to~\eqref{eq:ob_pde} is available. In what follows, we call this information a ``prior''.
This prior is denoted by $\bm{x} \mapsto u_{\theta}(\bm{x})$ with parameters $\theta$, and we assume that it can be constructed with the desired regularity $u_{\theta} \in H^{q+1}(\Omega)\cap H_0^1(\Omega)$ for $1\leqslant q\leqslant k$, where $k$ is the polynomial degree of the enriched FEM.
In this section and the following two, the prior is assumed to be a general function.
However, from \cref{sec:prior_construction} onwards, the prior will be the prediction of a parametric PINN.
In \cref{sec:modified_problem_add}, we first show how to use this prior to enriching classical finite element spaces.
Then, in \cref{sec:error_estimates_add}, we prove a convergence estimate for the resulting method.

\subsection{Construction of the modified problem}\label{sec:modified_problem_add}

In the general setting of FEM, we follow the Bubnov--Galerkin method~\cite{Ern2004TheoryAP}, where the basis functions and the numerical solutions are in the same space (see~\eqref{eq:approachform}, where both $u_h$ and $v_h$ are in $V_h^0$).
As we intend to enrich the classical approximation space, we exploit the idea formalized as the Petrov--Galerkin method (e.g.,~\cite{j2005introduction,brenner2008mathematical,demkowicz2023mathematical}), where the test and trial functions belong to different spaces.
This approach is often used for convection-dominated problems, see~\cite{ALMEIDA1997291}.
We propose to enrich the trial space using the prior $u_\theta$ by defining
\begin{equation}\label{eq:Vh_add}
	V_h^+ = \left\{
	u_h^+= u_{\theta} + p_h^+, \quad p_h^+ \in V_h^0
	\right\},
\end{equation}
and we use the space $V_h^0$ for the test functions.
Since we have assumed that $u_{\theta} \in H^{q+1}(\Omega)\cap H^1_0(\Omega)$,~$V_h^+$ is also a subset of $V^0$, like $V_h^0$.
Plugging this new trial space into the approximate problem~\eqref{eq:approachform}, we obtain the formulation
\begin{equation}\label{eq:add1}
	\text{Find } u_h^+ \in V_h^+ \text{ such that, } \; \forall v_h\in V_h^0, \; a(u_h^+,v_h)=l(v_h),
\end{equation}
which leads to the following approximation problem:
\begin{equation}\label{eq:approachform_add}
	\text{Find } p_h^+ \in V_h^0 \text{ such that, } \;
	\forall v_h \in V_h^0, \; a(p_h^+,v_h) = l(v_h) - a(u_{\theta},v_h).
\end{equation}
Therefore, we obtain a classical Galerkin approximation with a modified source term.

\subsection{Convergence analysis}\label{sec:error_estimates_add}

The objective is to prove that the FEM solution to problem~\eqref{eq:approachform_add} converges, with an error depending on the quality of the prior. \ROneN{Due to this dependence on the prior (neural network in the numerical section), we are not conducting an in-depth study of this gain. We are simply aiming to show that adding the prior does not degrade the convergence orders and that it is the high-order derivatives that control the gain in comparison with classical finite element method, which helps to guide choices regarding the prior.}

\begin{thrm}\label{lem:error_estimation_add}
	Let $u\in H^{q+1}(\Omega)$ be the solution to problem~\eqref{eq:weakform} and $u_{\theta}\in H^{q+1}(\Omega)\cap H_0^1(\Omega)$ be a prior on $u$.
	We consider $u_h^+\in V_h^+$ as the solution to the discrete problem~\ROne{\eqref{eq:add1}} with $V_h^+$ the modified trial space defined in~\eqref{eq:Vh_add}.
	The following estimates hold.
	For all $1\leqslant q\leqslant k$,
	\begin{equation}\label{eq:error_add}
		| u-u_h^+|_{H^1} \leqslant C_q\dfrac{\gamma}{\alpha} C_\text{\rm gain}^+ \, h^{q} |u|_{H^{q+1}}
	\end{equation}
	and, \ROne{supposing that dual problem \eqref{eq:dual} admits strong solutions (i.e. for all ${\xi}\in L^2(\Omega)$ there exists an associated solution $w\in H^1_0(\Omega)\cap H^2(\Omega)$ to \eqref{eq:dual}),}
	\begin{equation}\label{eq:error_addL2}
		\| u-u_h^+\|_{L^2} \leqslant C_e C_1 C_q\dfrac{\gamma^2}{\alpha} C_\text{\rm gain}^+ \, h^{q+1} |u|_{H^{q+1}},
	\end{equation}
	with $C_e$, $C_1$, $C_q$, $\gamma$, $\alpha$ defined in \cref{sec:FEM} and
	\begin{equation}\label{eq:gain_add}
		C_\text{\rm gain}^+= \frac{| u-u_{\theta} |_{H^{q+1}}}{| u |_{H^{q+1}}}.
	\end{equation}
\end{thrm}
\begin{rmrk}\label[rmrk]{rmk:gain_add}
	The constant \smash{$C_\text{\rm gain}^+$}
	represents the potential gain compared to the error of the classical FEM presented in \cref{thm:classical_error_estimate}. Note that this constant is the same in $L^2$ norm and $H^1$ semi-norm.
\end{rmrk}

\begin{proof}[Proof of \cref{lem:error_estimation_add}]
	\textbf{$H^1$-error:}   To prove~\eqref{eq:error_add}, we adapt the proof of Céa's lemma to the additive prior case.
	Considering the trial space defined in~\eqref{eq:Vh_add}, the numerical solution $u_h^+$ is given by
	\[
		u_h^+=u_{\theta}+p_h^+,
	\]
	with $p_h^+ \in V_h^0 \subset V$ solution to~\eqref{eq:approachform_add}.
	We have
	\begin{equation*}
		a(u-u_h^+, u-u_h^+)
        =
		a\big(u-u_h^+,(u-u_{\theta})-v_h\big)  +
		a\big(u-u_h^+, v_h-p_h^+\big),
        \quad \forall v_h \in V_h^0.
	\end{equation*}

	Let us first \ROne{treat} the second term on the right-hand side. By Galerkin orthogonality (difference of the continuous problem~\eqref{eq:weakform} and discrete problem~\eqref{eq:add1}),
    we obtain
	\[
		a\big(u-u_h^+, v_h-p_h^+\big)=0, \quad \forall v_h \in V_h^0.
	\]

	Denoting by $\alpha$ and $\gamma$ the
	coercivity and continuity constants of the bilinear form $a$, we have
	\begin{alignat*}{3}
		\alpha \big| u-u_h^+\big|_{H^1}^2 & \leqslant
		a\big(u-u_h^+,u-u_h^+\big) =
		a\big(u-u_h^+,(u-u_{\theta})-v_h\big),
        & & \quad \forall v_h \in V_h^0, \\
		& \leqslant \gamma \big| u-u_h^+\big|_{H^1} \big| (u-u_{\theta})-v_h \big|_{H^1},
		& & \quad \forall v_h \in V_h^0,
	\end{alignat*}
	which immediately leads to
	\[
		| u-u_h^+|_{H^1} \leqslant \frac{\gamma}{\alpha} \big| (u-u_{\theta})-v_h \big|_{H^1}, \quad \forall v_h \in V_h^0.
	\]
	We apply it to $v_h=\mathcal{I}_h(u-u_\theta) \in V_h^0$
	with $\mathcal{I}_h$ the Lagrange interpolation operator~\eqref{eq:Ih} in $V_h$, it holds
	using interpolation estimate given in \cref{th:interpol},
	\[
		|u-u_h^+|_{H^1} \leqslant C_q\frac{\gamma}{\alpha} h^{q} | u-u_{\theta} |_{H^{q+1}},
	\]
	with $C_q$ defined in \cref{sec:FEM}.

	The above expression can be rewritten as
	\begin{equation}\label{eq:trucbidule}
		| u-u_h^+|_{H^1} \leqslant C_q\frac{\gamma}{\alpha}  C_\text{\rm gain}^+ \, h^{q}|u|_{H^{q+1}} \,,
	\end{equation}
	with
	\[
		C_\text{\rm gain}^+ = \frac{| u-u_{\theta} |_{H^{q+1}}}{| u |_{H^{q+1}}},
	\]
	which completes the first part of the proof.

	\textbf{$L^2$-error:}
	We will follow the Aubin--Nitsche technique.
	\ROne{
    Consider $w\in H^1_0(\Omega)\cap H^2(\Omega)$ solution to
	\[
		\mathcal{L}^* w=u-u_h^+,
	\]
	with homogeneous Dirichlet boundary condition. Thanks to \cref{th:ellip}, one has
	\begin{equation}\label{eq:wH2}
		\|w\|_{H^2}\leqslant C_e\|u-u_h^+\|_{L^2}.
	\end{equation}}
	Using the Galerkin orthogonality and the continuity of the bilinear form $a$,
	\[
		\|u-u_h^+\|_{L^2}^2 = a(u-u_h^+,w-I_h w)
		\leqslant \gamma |u-u_h^+|_{H_1}|w-I_h w|_{H_1}.
	\]
	Thanks to \cref{th:interpol}
	and~\eqref{eq:wH2},
	\[
		|w-I_h w|_{H_1}\leqslant C_e C_1 h \|u-u_h^+\|_{L^2},
	\]
	which leads to the conclusion by using~\eqref{eq:trucbidule}.
\end{proof}

\begin{rmrk}\label[rmrk]{rmk:C_gain_additif}
	The gain constant $C_\text{\rm gain}^+$ defined in~\eqref{eq:gain_add} shows that the closer the prior is to the solution,
	the smaller is the error constant associated with the FEM while keeping the same order of accuracy.
	Therefore, as soon as \smash{$C_\text{\rm gain}^+ < 1$}, the FEM with additive prior will be more accurate than the classical one.
	While this gives us a particularly flexible constraint, our objective is to balance this gain by relaxing the contribution $h^q$, using a coarser grid and low-order polynomial, to reduce the computational cost of the FEM while maintaining accuracy.
	Nonetheless, the gain is related to the~$L^2$ error associated with the derivatives of $(q+1)$\textsuperscript{th} order (with $1\leqslant q\leqslant k$) of the prior.
	This shows that the prior must accurately approximate the derivatives of the solution in addition to the solution itself.
	This highlights that we need to build our prior by ensuring a good approximation of the derivatives of the solution.
	It also shows that the higher the order of the finite elements $k$ is, the better our prior should approximate the higher-order derivatives.
	Therefore, it is more appropriate to use only low-order FEM so that $k$ remains small.
\end{rmrk}

    \section{Enriching the finite element method with multiplicative priors}\label{sec:multiplicative_prior}
    This section employs the same assumptions as in \cref{sec:additive_prior}, namely that we have a sufficiently smooth prior $u_\theta$ on the solution~$u$ of the PDE~\eqref{eq:ob_pde}.
However, this prior will now be multiplied to elements of $V_h$ rather than added to them.
We construct the underlying modified problem in \cref{sec:modified_problem_mul}.
Then, similarly to the additive approach of \cref{sec:additive_prior}, error estimates are obtained in \cref{sec:error_estimates_mul}. \ROneN{The objectives of this section on error estimates are the same as in the previous section.}

\subsection{Construction of the modified problem}\label{sec:modified_problem_mul}

To construct the modified problem in this case, we must ensure that the prior $u_\theta$ never vanishes.
\RTwo{To that end, we lift the initial problem~\eqref{eq:ob_pde} by a constant $M \in \mathbb{R}_+$, chosen large enough to ensure that $u_M=u+M>0$, to get}
\begin{equation}\label{eq:ob_pde_M}
	\begin{dcases}
		\mathcal{L}(u_M)=f, &
		\text{\quad in } \Omega,          \\
		u_M = M,                  &
		\text{\quad on } \partial \Omega. \\
	\end{dcases}
\end{equation}
We then introduce the associated variational problem, defined by
\begin{equation}\label{eq:weakform2}
	\text{Find } u_M = u + M,~\text{with }u \in V^0 \text{ such that, } \;
	\forall v\in V^0,\; a(u_M,v)=l(v).
\end{equation}
Therefore, solving~\eqref{eq:weakform2}, we recover the solution $u$ of the initial problem~\eqref{eq:ob_pde} by setting \[u = u_M - M.\]
The prior 
\[
	u_{\theta,M}=u_\theta+M>0
\]
is associated with problem~\eqref{eq:weakform2}.

Let us now introduce a new modified finite element space, defined by
\begin{equation}\label{eq:Vh_mul}
	V_h^\times = \left\{
	u_{h,M}^\times = u_{\theta,M} \; p_h^\times,
	\quad p_h^\times \in 1+V_h^0
	\right\},
\end{equation}
with, for all $\bm{x}\in \Omega$, $u_{\theta,M}(\bm{x})\neq 0$.
From~\eqref{eq:weakform2}, this leads to the following approximate formulation:
\begin{equation}\label{eq:approachform_mul}
	\text{Find } p_h^\times \in 1+ V_h^0 \text{ such that, } \; \forall v_h \in V_h^0, \; a \big(u_{\theta,M} \; p_h^\times,u_{\theta,M}  v_h \big) = l(u_{\theta,M} v_h).
\end{equation}
Therefore, solving~\eqref{eq:approachform_mul}, we recover the solution $u_h^\times\in V_h^\times-M$ of the original problem~\eqref{eq:ob_pde} by setting
\[
	u_h^\times = u_{h,M}^\times - M.
\]
Based on the $N_{\text{dofs}}$ dofs $\big(\bm{x}^{(i)}\big)_{i \in \{1,\ldots,N_{\text{\rm dofs}}\}}$ of the mesh, we consider the interpolation operator on $V_h^\times$ given by
\begin{equation*}
	\tilde{\mathcal{I}}_h:
	C^0(\Omega) \ni v \mapsto
	\sum_{i=1}^{N_\text{dofs}}\frac{v\big(\bm{x}^{(i)}\big)}{u_{\theta,M}\big(\bm{x}^{(i)}\big)} \tilde{\psi_i} \in V_h^\times,
\end{equation*}
where the shape functions $\tilde{\psi_i}$ associated to $V_h^\times$ are defined by
\[
	\tilde{\psi_i}={u_{\theta,M}} \; {\psi}_i,
\]
with ${\psi}_i$ the classical shape functions presented in \cref{sec:FEM}.
Note that the new interpolation operator $\tilde{\mathcal{I}}_h$ is related to the classical Lagrange interpolation operator defined in~\eqref{eq:Ih} as follows
\begin{equation}\label{eq:relation_Ih_Ih_tilde}
	\forall v \in C^0(\Omega), \qquad
	\tilde{\mathcal{I}}_h(v)
	=
	u_{\theta,M} \;
	\mathcal{I}_h \left( \frac v {u_{\theta,M}} \right).
\end{equation}

\subsection{Convergence analysis}\label{sec:error_estimates_mul}

In this section, we finally prove that the modified FEM~\eqref{eq:approachform_mul} converges to the solution to~\eqref{eq:ob_pde_M},
and that it satisfies the same type of estimate as the classical one.
Equipped with the lifting trick from \cref{sec:modified_problem_mul},
we can state the following convergence \RTwo{theorem}.

\begin{thrm}\label{lem:error_estimate_multiplicative}
	Let $u_M\in H^{q+1}(\Omega)$ be the solution of the enhanced problem~\eqref{eq:weakform2}
    and $u_{\theta,M}\in M+H^{q+1}(\Omega)\cap H^1_0(\Omega)$ be a prior on $u_M$.
    We consider \smash{$u_{h,M}^\times\in V_h^\times$} the solution to the finite element problem~\eqref{eq:approachform_mul} with
    \smash{$V_h^\times$} the modified trial space defined in~\eqref{eq:Vh_mul},
    considering $\mathbb{P}_k$ polynomials.
    We define $u=u_M-M$ and $u_h^\times=u_{h,M}^\times-M$.
    Then, for all $1\leqslant q\leqslant k$
	\begin{equation}\label{eq:error_mul}
		| u-u_h^\times|_{H^1} \leqslant C_q\dfrac{\gamma}{\alpha} C_{\text{\rm gain},H^1}^{\times,M} h^{q}| u |_{H^{q+1}}
	\end{equation}
	and, \ROne{supposing that dual problem \eqref{eq:dual} admits strong solutions (i.e.  for all ${\xi}\in L^2(\Omega)$, there exists an associated solution $w\in H^1_0(\Omega)\cap H^2(\Omega)$ to \eqref{eq:dual}),}
	\begin{equation}\label{eq:error_mulL2}
        \| u-u_h^\times\|_{L^2} \leqslant C_e C_1 C_q\dfrac{\gamma^2}{\alpha} C_{\text{\rm gain},L^2}^{\times,M} \, h^{q+1} |u|_{H^{q+1}},
    \end{equation}
    with $C_e$, $C_1$, $C_q$, $\gamma$, $\alpha$ defined in \cref{sec:FEM}, and where
    \begin{equation}\label{eq:gain_mul}
		C_{\text{\rm gain},H^1}^{\times,M} = \left| \frac{u_M}{u_{\theta,M}} \right|_{H^{q+1}} \frac{\| u_{\theta,M}\|_{W^{1,\infty}}}{| u |_{H^{q+1}}},
	\end{equation}
	and
    \begin{equation}\label{eq:gain_mulL2}
		C_{\text{\rm gain},L^2}^{\times,M} =
		C_{\theta,M}\left| \frac{u_M}{u_{\theta,M}} \right|_{H^{q+1}} \frac{\| u_{\theta,M}\|_{W^{1,\infty}}^2}{| u |_{H^{q+1}}},
	\end{equation}
	with
    \begin{equation}\label{eq:CthetaM}C_{\theta,M}=\|u_{\theta,M}^{-1}\|_{L^{\infty}}
		+2|u_{\theta,M}^{-1}|_{W^{1,\infty}}
		+|u_{\theta,M}^{-1}|_{W^{2,\infty}}.
	\end{equation}
\end{thrm}

\begin{rmrk}\label[rmrk]{rmk:gain_mul}
The constants \smash{$C_{\text{\rm gain},H^1}^{\times,M}$}
and \smash{$C_{\text{\rm gain},L^2}^{\times,M}$}
represent the potential gains in both $H^1$ semi-norm and $L^2$ norm
when using the multiplicative approach, compared
to the error of the classical FEM presented in \cref{thm:classical_error_estimate} with $\mathbb{P}_k$ polynomials.
\end{rmrk}

\begin{proof}[Proof of \cref{lem:error_estimate_multiplicative}]
	\textbf{$H^1$-error:}
	Considering the trial space defined in~\eqref{eq:Vh_mul}, the numerical solution $u_{h,M}^\times$ is given by
	\[
		u_{h,M}^\times=u_{\theta,M} \; p_h^\times,
	\]
	with $p_h^\times\in 1+ V_h^0\subset 1 + V^0$  solution to~\eqref{eq:approachform_mul}.
	By coercivity of $a$,
		\[
		\alpha| u_M-u_{h,M}^\times|_{H^1}^2
		\leqslant a(u_M-u_{h,M}^\times,u_M-u_{h,M}^\times) .
	\]
	Thanks to~\RBoth{\eqref{eq:weakform2}} and~\eqref{eq:approachform_mul}, we have the following Galerkin orthogonality: for all $v_h\in V_h^0$,
	\begin{equation}\label{eq:orthGalmult}
		a(u_M-u_{h,M}^\times, u_{\theta,M}v_h)=0.
	\end{equation}
    \RTwo{Now, choosing
    \begin{equation*}
        v_h=p_h^{\times} - \mathcal{I}_h\left(\frac{u_M}{u_{\theta,M}}\right),
    \end{equation*}
    we note that $v_h$ is well-defined since the lifting trick ensures that $u_{\theta,M} > 0$, and that it belongs to $V_h^0$.
	With this choice, we deduce by the definition~\eqref{eq:relation_Ih_Ih_tilde} of $\tilde{\mathcal{I}}_h$ that}
	\[
		a(u_M-u_{h,M}^\times,u_M-u_{h,M}^\times)
		= a\left(u_M-u_{h,M}^\times,u_M-u_{\theta,M}\mathcal{I}_h\left(\frac{u_M}{u_{\theta,M}}\right)\right)
		= a(u_M-u_{h,M}^\times,u_M-\tilde{\mathcal{I}}_h(u_M)).
	\]
	By continuity of $a$,
	\begin{equation}\label{eq:tructruc}
		| u_M-u_{h,M}^\times|_{H^1} \leqslant \frac{\gamma}{\alpha}| u_M-\tilde{\mathcal{I}}_h(u_M)|_{H^1} .
	\end{equation}
	Again, using the definition~\eqref{eq:relation_Ih_Ih_tilde} of $\tilde{\mathcal{I}}_h$,
	\[
		|u_M-\tilde{\mathcal{I}}_h(u_M)|_{H^{1}} \leqslant
		\|u_{\theta,M}\|_{W^{1,\infty}} \left\|\frac{u_M}{u_{\theta,M}}-\mathcal{I}_h\left(\frac{u_M}{u_{\theta,M}}\right)\right\|_{H^{1}}.
	\]

	Finally, applying interpolation estimate given in \cref{th:interpol}, it holds
	\begin{equation}\label{eq:interpol tilde Ih}
		|u_M-\tilde{\mathcal{I}}_h(u_M)|_{H^1} \leqslant C_q
		\|u_{\theta,M}\|_{W^{1,\infty}} h^q \left|\frac{u_M}{u_{\theta,M}}\right|_{H^{q+1}} \,,
	\end{equation}
	with $C_q$ defined in \cref{sec:FEM}.
	Combining the last inequality with~\eqref{eq:tructruc}, we obtain
	\begin{equation}\label{eq:H1 mult}
		| u-u_h^\times|_{H^1} = | u_M-u_{h,M}^\times|_{H^1} \leqslant C_q \dfrac{\gamma}{\alpha} C_{\text{\rm gain},H^1}^{\times,M} h^{q}| u |_{H^{q+1}},
	\end{equation}
	with $C_{\text{\rm gain},H^1}^{\times,M}$ given in~\eqref{eq:gain_mul}, which conclude the first part of the proof.

	\textbf{$L^2$-error:}
	Again, we follow the Aubin-Nitsche strategy here. 
    \ROne{Consider $w-M\in H^1_0(\Omega)\cap H^2(\Omega)$ solution to
	\[
		\mathcal{L}^*w=u-u_h^{\times}=u_M-p_h^{\times}u_{\theta,M},
	\]
	with $w=M$ on $\partial\Omega$. Thanks to Theorem \ref{th:ellip}, one has
	\begin{equation*}
		\|w\|_{H^2}\leqslant C_e\|u-u_h^{\times}\|_{L^2}.
	\end{equation*}}    
    Then,  using the Galerkin orthogonality~\eqref{eq:orthGalmult} for $v_h=\mathcal{I}_h\left(\frac{u_M}{u_{\theta,M}}\right)$,
	\[
		\|u-u_h^{\times}\|_{L^2}^2
		=\|u_M-p_h^{\times}u_{\theta,M}\|_{L^2}^2
		=a(u_M-p_h^{\times}u_{\theta,M},w)
		=a(u_M-p_h^{\times}u_{\theta,M},w-\tilde I_h(w)).
	\]
    Hence, by continuity of $a$,
	\[
		\|u-u_h^{\times}\|_{L^2}^2 \leqslant \gamma |u_M-p_h^{\times}u_{\theta,M}|_{H^1}|w-\tilde I_h(w)|_{H^1}.
	\]
	Using~\eqref{eq:H1 mult} and~\eqref{eq:interpol tilde Ih} for $q=1$ to the term in the right-hand side,
	\begin{equation*}
		\|u-u_h^{\times}\|_{L^2}^2
		\leqslant C_1 C_q \dfrac{\gamma^2}{\alpha} \|u_{\theta,M}\|_{W^{1,\infty}} \left|\frac{w}{u_{\theta,M}}\right|_{H^{2}} C_{\text{\rm gain},H^1}^{\times,M} h^{q+1}| u |_{H^{q+1}}.
	\end{equation*}
	Moreover,
    \[
		\left|\frac{w}{u_{\theta,M}}\right|_{H^{2}}\leqslant C_{\theta,M}\|w\|_{H^{2}},
	\]
    with $C_{\theta,M}$ given in~\eqref{eq:CthetaM}.
	Thanks to the elliptic regularity, we obtain
    \begin{equation*}
        \| u-u_h^\times\|_{L^2} \leqslant C_e C_1 C_q\dfrac{\gamma^2}{\alpha} C_{\text{\rm gain},L^2}^{\times,M} \, h^{q+1} |u|_{H^{q+1}},
    \end{equation*}
	with $C_{\text{\rm gain},L^2}^{\times,M}$ defined in~\eqref{eq:gain_mulL2}.
\end{proof}

\begin{rmrk}\label[rmrk]{rmk:C_gain_multiplicatif}
	We note that the gain constants $C_{\text{\rm gain},H^1}^{\times,M}$ and $C_{\text{\rm gain},L^2}^{\times,M}$
	are similar to the constant $C_\text{\rm gain}^+$
	introduced in \cref{sec:additive_prior},
	in that, it depends on high-order derivatives of the prior.
	Hence, a high-quality prior will necessarily involve
	a good approximation of the derivatives of the exact solution,
	and \cref{rmk:C_gain_additif} also applies in the present context.
	The major difference with the additive approach lies
	in the choice of the lifting constant $M$.
	To better understand this dependency in $M$,
	the following section provides a
	study of the behavior of our two gain constants
	when $M$ goes to infinity.
	Moreover, the actual choice of $M$ will be
	numerically investigated in \cref{sec:numerical_results}.
\end{rmrk}

    \section{Comparison of the two enriched methods}\label{sec:comparison_add_mul}
    This section aims to compare the two proposed methods, namely the additive approach presented in \cref{sec:additive_prior} and the multiplicative one proposed in \cref{sec:multiplicative_prior}. Recall that the constant $M$ is chosen in the multiplicative approach so that $u_M>0$.
Let $u$ be the solution of problem~\eqref{eq:weakform} and $u_{\theta}\in H^{q+1}(\Omega)\cap H^1_0(\Omega)$ be a prior on $u$, with $1\leqslant q\leqslant k$ ($k$ the polynomial degree of the finite element method).
For clarity, we first recall the additive and multiplicative enrichments,
and their main error estimates.

\boldparagraph{Additive approach} We consider $u_h^+\in V_h^+$ as the solution to the finite element method associated to problem~\eqref{eq:approachform_add} with $V_h^+$ the modified trial space defined in~\eqref{eq:Vh_add}, considering $\mathbb{P}_k$ polynomials. Using \cref{lem:error_estimation_add}, \RBoth{we have for $1\leqslant q\leqslant k$, the following additive theoretical gain constant}
\begin{equation}\label{eq:C_gain_add_in_comparison}
	C_\text{\rm gain}^+= \frac{| u-u_{\theta} |_{H^{q+1}}}{| u |_{H^{q+1}}}.
\end{equation}

\boldparagraph{Multiplicative approach} Let $u_M=u+M$ be the solution of the enhanced problem~\eqref{eq:weakform2} and $u_{\theta,M}=u_\theta+M$ be a prior on $u_M$. We consider $u_{h,M}^\times\in V_h^\times$ the solution to the finite element  problem~\eqref{eq:approachform_mul} with $\smash{V_h^\times}$ the modified trial space defined in~\eqref{eq:Vh_mul}, considering $\mathbb{P}_k$ polynomials. Using \cref{lem:error_estimate_multiplicative}, \RBoth{we have for $1\leqslant q\leqslant k$, the following multiplicative theoretical gain constants}
\begin{equation}\label{eq:C_gain_mul_in_comparison}
    C_{\text{\rm gain},H^1}^{\times,M} = \left| \frac{u_M}{u_{\theta,M}} \right|_{H^{q+1}} \frac{\| u_{\theta,M}\|_{W^{1,\infty}}}{| u |_{H^{q+1}}},
\end{equation}
and
\begin{equation}\label{eq:gain_mulL2_in_comparison}
    C_{\text{\rm gain},L^2}^{\times,M} =
    C_{\theta,M}\left| \frac{u_M}{u_{\theta,M}} \right|_{H^{q+1}} \frac{\| u_{\theta,M}\|_{W^{1,\infty}}^2}{| u |_{H^{q+1}}},
\end{equation}
with $C_{\theta,M}$ given in~\eqref{eq:CthetaM}.

\boldparagraph{Comparison of the two approaches} The following result proves that the upper bound in~\eqref{eq:error_mul} and in~\eqref{eq:error_mulL2} converges respectively to the one in~\eqref{eq:error_add} and in~\eqref{eq:error_addL2} when $M$ goes to infinity. In other words, the multiplicative gain constants defined in~\eqref{eq:C_gain_mul_in_comparison} and~\eqref{eq:gain_mulL2_in_comparison} converge to the additive gain constant defined in~\eqref{eq:C_gain_add_in_comparison} when $M$ goes to infinity.

\begin{thrm}\label{thm:comparison_add_mul}
	We have
	\begin{equation}\label{eq:convCH1}
        C_{\text{\rm gain},H^1}^{\times,M}
		\underset{M\rightarrow\infty}{\longrightarrow}
		C^{+}_{\text{\rm gain}},
	\end{equation}
	and
	\begin{equation}\label{eq:convCL2}
		C_{\text{\rm gain},L^2}^{\times,M}
		\underset{M\rightarrow\infty}{\longrightarrow}
		C^{+}_{\text{\rm gain}}.
	\end{equation}
\end{thrm}

\begin{proof}
\textbf{Convergence in $H^1$:}
	According to the expressions~\eqref{eq:C_gain_add_in_comparison} and~\eqref{eq:C_gain_mul_in_comparison}
	of the gain constants, the objective of the proof is to show that
	\begin{equation*}
		\|u_{\theta,M}\|_{W^{1,\infty}} \left|\frac{u_M}{u_{\theta,M}}\right|_{H^{q+1}} \underset{M\rightarrow\infty}{\longrightarrow} |u-u_\theta|_{H^{q+1}}.
	\end{equation*}
	Denoting by
	\begin{equation*}
		E_\theta=u - u_{\theta},
	\end{equation*}
	the error made by the prior $u_{\theta}$
	when approximating the solution $u$, we have
	\begin{equation*}
		|u-u_\theta|_{H^{q+1}} = |E_\theta|_{H^{q+1}}.
	\end{equation*}
	On the one hand, we have,
	\begin{equation*}
		\left\|u_{\theta,M}\right\|_{W^{1,\infty}}
		=
		\left\|u_\theta + M \right\|_{W^{1,\infty}}
		=
		M
		\left\|1 + \frac {u_\theta} M \right\|_{W^{1,\infty}}.
	\end{equation*}
	On the other hand, we have,
	\begin{equation*}
		\left|\frac{u_M}{u_{\theta,M}}\right|_{H^{q+1}}
		=
		\left|\frac{u + M}{u_\theta + M}\right|_{H^{q+1}}
		=
		\left|\frac{u - u_\theta + u_\theta + M}{u_\theta + M}\right|_{H^{q+1}}
		=
		\left|1 + \frac{u - u_\theta}{u_\theta + M}\right|_{H^{q+1}}
		=
		\dfrac{1}{M}\left| \frac{E_\theta}{1 + \frac {u_\theta} M}\right|_{H^{q+1}}.
	\end{equation*}
	Multiplying these expressions, we obtain
	\begin{equation}\label{eq:error_mul_q2}
		\|u_{\theta,M}\|_{W^{1,\infty}} \left|\frac{u_M}{u_{\theta,M}}\right|_{H^{q+1}}= \underbrace{\left\|1+\frac{u_\theta}{M}\right\|_{W^{1,\infty}}}_\text{(\MakeUppercase{\romannumeral1})} \;
		\underbrace{\left|\frac{E_\theta}{1+\frac{u_\theta}{M}}\right|_{H^{q+1}}}_\text{(\MakeUppercase{\romannumeral2})}.
	\end{equation}
	We now estimate term by term the right-hand side of the
	above equality~\eqref{eq:error_mul_q2}, looking at their
	limits when $M$ goes to infinity.

    \smallskip

	\textbf{Term (\MakeUppercase{\romannumeral1}):} By decomposing the first term, we obtain
	\[
		\text{(\MakeUppercase{\romannumeral1})}=\left\|1+\frac{u_\theta}{M}\right\|_{W^{1,\infty}}=\left\|1+\frac{u_\theta}{M}\right\|_{L^\infty} + \frac{1}{M}\|\nabla u_\theta\|_{L^\infty}\underset{M\to\infty}{\longrightarrow} 1,
	\]
	since
	\begin{equation}\label{eq:majoration}
		1-\frac{\|u_\theta\|_{\RTwoN{L^\infty}}}{M}
		\leqslant
		\left\|1+\frac{u_\theta}{M}\right\|_{L^\infty}
		\leqslant
		1+\frac{\|u_\theta\|_{\RTwoN{L^\infty}}}{M}.
	\end{equation}

    \smallskip

	\textbf{Term (\MakeUppercase{\romannumeral2}):}
	Let us prove that
	\begin{equation*}
	\text{(\MakeUppercase{\romannumeral2})}=	\left|\frac{E_\theta}{1+\frac{u_\theta}{M}}\right|_{H^{q+1}}\underset{M\rightarrow\infty}{\longrightarrow} |E_\theta|_{H^{q+1}}.
	\end{equation*}
\RTwo{By the triangular inequality,
    	\begin{align*}
\left|\frac{E_\theta}{1+\frac{u_\theta}{M}}\right|_{H^{q+1}}^2
\leqslant\left|\frac{E_\theta}{1+\frac{u_\theta}{M}}-E_{\theta}\right|_{H^{q+1}}^2
+\left|E_{\theta}\right|_{H^{q+1}}^2
	\end{align*}
and
        	\begin{align*}
\left|E_{\theta}\right|_{H^{q+1}}^2
\leqslant\left|\frac{E_\theta}{1+\frac{u_\theta}{M}}-E_{\theta}\right|_{H^{q+1}}^2
+\left|\frac{E_\theta}{1+\frac{u_\theta}{M}}\right|_{H^{q+1}}^2.
	\end{align*}
    Hence, by definition of the semi-norm $\left|\,\cdot\,\right|_{H^{q+1}}$,}
    \begin{align*}
		\left|\left|\frac{E_\theta}{1+\frac{u_\theta}{M}}\right|_{H^{q+1}}^2-|E_\theta|_{H^{q+1}}^2\right| & \leqslant \left|\frac{E_\theta}{1+\frac{u_\theta}{M}}-E_\theta\right|_{H^{q+1}}^2                                    \\
		& = \left\|\nabla^{q+1}\left(\frac{E_\theta}{1+\frac{u_\theta}{M}}-E_\theta\right)\right\|_{L^2}^2
		\leqslant |\Omega| \left\|\nabla^{q+1}\left(\frac{E_\theta}{1+\frac{u_\theta}{M}}-E_\theta\right)\right\|_{L^\infty}^2.
	\end{align*}
	Then, using the general Leibniz rule, we have,

    \begin{align*}
		\left\|\nabla^{q+1}\left(\frac{E_\theta}{1+\frac{u_\theta}{M}}-E_\theta\right)\right\|_{L^\infty} & \leqslant
		\left\|\frac{\nabla^{q+1}E_\theta}{1+\frac{u_\theta}{M}}-\nabla^{q+1}E_\theta\right\|_{L^\infty}
		+\sum_{s=1}^{q+1} \begin{pmatrix} q+1 \\ s \end{pmatrix}
		\left\|\nabla^{q+1-s}E_\theta\right\|_{L^{\infty}}\left\|\nabla^{s}\left(\frac{1}{1+\frac{u_\theta}{M}}\right)\right\|_{L^\infty} \\
		& \leqslant
		\underbrace{ \left\|\frac{\nabla^{q+1}E_\theta}{1+\frac{u_\theta}{M}}-\nabla^{q+1}E_\theta\right\|_{L^\infty} }_{(1)}
		+\|E_\theta\|_{W^{q+1,{\infty}}} \sum_{s=1}^{q+1} \begin{pmatrix} q+1 \\ s \end{pmatrix}
		\underbrace{ \left\|\nabla^{s}\left(\frac{1}{1+\frac{u_\theta}{M}}\right)\right\|_{L^\infty} }_{(2)} .
	\end{align*}
	We now estimate terms (1) and (2) in the right-hand side of the above inequality.

	\textbf{Term (1):}
	Taking $M> \|u_{\theta}\|_{L^\infty}$, we obtain
	\begin{align*}
		(1)\leqslant	\left\|\frac{1}{1+\frac{u_\theta}{M}}-1\right\|_{L^\infty}\|\nabla^{q+1}E_\theta\|_{L^\infty}
		& \leqslant \frac{1}{M} \frac{\|u_\theta\|_{L^\infty}}{1-\frac{\|u_\theta\|_{L^\infty}}{M}} \|\nabla^{q+1}E_\theta\|_{L^\infty} \underset{M\to\infty}{\longrightarrow}0.
	\end{align*}

	\textbf{Term (2):} \RTwoN{Let us define $g=1+\frac{u_\theta}{M}$. The Fa\`{a} di Bruno formula (see~\cite{faadibruno2002}) gives the following expression for the $s$-th derivative of the inverse of $g$:}

	\begin{equation*}
		\RTwoN{\nabla^s \left(\frac{1}{g}\right) = \frac{s!}{g^{s+1}} \sum \frac{(-1)^\alpha \alpha!}{\prod_{i=1}^s (i!)^{m_i} m_i!} g^{s-\alpha} \prod_{i=1}^s \left(\nabla^i g\right)^{m_i},}
	\end{equation*}
	\RTwoN{where the sum covers all $s$-tuples $(m_1,\dots,m_s)$ satisfying the constraint $s = \sum_{i=1}^{s}im_i$ and $\alpha = \sum_{i=1}^{s}m_i$.}

	\RTwoN{Thus, using \eqref{eq:majoration},} there exists a constant $C>0$ such that, for any $M> \|u_{\theta}\|_{L^\infty}$ and $1\leqslant s\leqslant q+1$, the following estimate holds 
	\begin{align*}
	(2)	&\leqslant \RTwoN{C \, \frac{1}{\displaystyle{\min_{x\in\Omega}}\left|1+\frac{u_\theta(x)}{M}\right|^{s+1}} \sum \left\|1+\frac{u_\theta}{M}\right\|_{L^{\infty}}^{s-\alpha} \prod_{i=1}^s \left\|\nabla^i \left(1 + \frac{u_\theta}{M}\right)\right\|_{L^{\infty}}^{m_i}} \\
	&\leqslant \RTwoN{C \, \frac{1}{\left(1-\frac{\|u_\theta\|_{L^\infty}}{M}\right)^{s+1}} \sum \left(1+\frac{\|u_\theta\|_{L^\infty}}{M}\right)^{s-\alpha} \frac{1}{M^{\alpha}} \prod_{i=1}^s \left\|\nabla^i u_\theta\right\|_{L^{\infty}}^{m_i}}
	\underset{M\to\infty}{\longrightarrow}0,
	\end{align*}
	which leads to~\eqref{eq:convCH1}.

    \smallskip

\textbf{Convergence in $L^2$:} Using the convergence of the gain for the $H^1$ semi-norm, we only need to prove that
\[
	C_{\theta,M}\|u_{\theta,M}\|_{W^{1,\infty}}
	\underset{M\rightarrow\infty}{\longrightarrow}1,
\]
with $C_{\theta,M}$ given by
\[
	C_{\theta,M}=\|u_{\theta,M}^{-1}\|_{L^{\infty}}
		+2|u_{\theta,M}^{-1}|_{W^{1,\infty}}
		+|u_{\theta,M}^{-1}|_{W^{2,\infty}}.
\]
Since
\[
	\dfrac{1}{M}\|u_{\theta,M}\|_{W^{1,\infty}}\underset{M\rightarrow\infty}{\longrightarrow}1,
\]
we are only required to prove that
\[
	MC_{\theta,M}
\underset{M\rightarrow\infty}{\longrightarrow}1.
\]
Considering $M> \|u_{\theta}\|_{L^\infty}$, we have
\[
	M\|u_{\theta,M}^{-1}\|_{L^{\infty}}
	=\left\|\dfrac{1}{1+\dfrac{u_\theta}{M}}\right\|_{L^{\infty}}
	\leqslant \dfrac{1}{1-\dfrac{\|u_\theta\|_{L^\infty}}{M}}
	\underset{M\rightarrow\infty}{\longrightarrow}1.
\]
Moreover,
\[
	2M|u_{\theta,M}^{-1}|_{W^{1,\infty}}
	=2M\left\|\dfrac{\nabla u_{\theta}}{(u_{\theta}+M)^2}\right\|_{L^{\infty}}
	\underset{M\rightarrow\infty}{\longrightarrow}0.
\]
Similarly,
\[
	M|u_{\theta,M}^{-1}|_{W^{2,\infty}}
	\underset{M\rightarrow\infty}{\longrightarrow}0,
\]
which leads to the conclusion~\eqref{eq:convCL2}.

\end{proof}

    \section{Prior construction using parametric PINNs}\label{sec:prior_construction}
    We have introduced new finite element approximation spaces in \cref{sec:additive_prior,sec:multiplicative_prior} depending on the construction of priors.
Physics-Informed Neural Networks (PINNs) are a good choice to build such priors.
Indeed, since PINNs minimize the PDE residual, they inherently give a good approximation of the derivative of the solution, in addition to the solution itself (see e.g.~\cite{RAISSI2019686}).
This section is therefore dedicated to introducing PINNs in \cref{sec:PINNs_parametric_PDE}, and then to show how to improve them in \cref{sec:improve_PINNs}.

\subsection{Physics-Informed Neural Networks for parametric PDEs}\label{sec:PINNs_parametric_PDE}

Physics-Informed Neural Networks, or PINNs, were introduced by~\cite{RAISSI2019686} for solving a PDE with \RTwo{neural networks}.
The main idea is to recast a PDE as an optimization problem.
We illustrate the method on our problem~\eqref{eq:ob_pde}, which we now extend to non-homogeneous Dirichlet boundary conditions.
Moreover, recall that the problem of interest is a parametric PDE. Unlike classical PINNs, which are trained for specific physical parameters, parametric PINNs seek to learn a generalized solution covering a range of parameters.
They incorporate these parameters as additional inputs to the network, allowing greater flexibility in solving problems where physical conditions or properties vary.
\RTwo{Moreover, since they are based on neural networks, PINNs are ideally suited to solving such higher-dimensional problems. This advantage is compounded when one uses Monte-Carlo integration to estimate the loss functions, as its convergence rate is also insensitive to the dimension.}

Considering $p$ parameters \smash{$\bm{\mu} = (\mu_1, \dots, \mu_{p}) \in \mathcal{M} \subset \mathbb{R}^{p}$}, with some parameter space $\mathcal{M}$, the parametric PDE~\eqref{eq:ob_pde} reads, with non-homogeneous Dirichlet boundary conditions:
\begin{equation}\label{eq:parametric_PDE}
    \begin{dcases}
        \mathcal{L}\big(u(\bm{x},\bm{\mu});\bm{x},\bm{\mu}\big) = f(\bm{x},\bm{\mu}), & \bm{x}\in\Omega,          \\
        u(\bm{x},\bm{\mu}) = g(\bm{x},\bm{\mu}),                                      & \bm{x}\in\partial \Omega, \\
    \end{dcases}
\end{equation}
with $g$ the trace of a $H^2$ function on $\partial \Omega$.
Note that the solution of the equation depends on the parameters~$\bm{\mu}$, as do the operator $\mathcal{L}$ and the boundary conditions.
We then denote $u_{\theta}(\cdot, \bm{\mu})$ the approximate PINN prediction for given parameters~$\bm{\mu}$.

The first idea of PINNs comes from the observation that, by construction, neural networks with smooth activation functions
are nothing but smooth functions of their weights and inputs.
Therefore, neural networks form natural candidates for approximating solutions to PDEs, especially with the advent of automatic differentiation tools.
In our case, a PINN is a neural network that takes $d+p$ inputs, where $d$ is the dimension of the space variable $\bm{x} \in \Omega$ and $p$ is the number of parameters $\bm{\mu} \in \mathcal{M}$. We denote by $u_{\theta}(\bm{x},\bm{\mu})$ the output, where $\theta$ are the learnable weights of the network.
Classically, this neural network is a coordinate-based neural network, such as a multi-layer perceptron (MLP).
It depends on several used-defined hyperparameters,
to be specified in the numerical experiments.

Once this network is defined,
solving the PDE can be rewritten as a minimization problem on $\theta$,
namely finding the optimal weights $\theta^\star$
that satisfy the following minimization problem:
\begin{equation}\label{eq:minimization_problem}
    \theta^\star = \argmin_{\theta}
    \big( \omega_r J_r(\theta) + \omega_b J_b(\theta) + \omega_\text{data} J_\text{data}(\theta) \big),
\end{equation}
with $\omega_r$, $\omega_b$ and $\omega_\text{data}$ some weights to balance the different terms of the loss function.
In~\eqref{eq:minimization_problem}, the loss function has three terms: the residual loss function
\begin{equation}\label{eq:residual_loss_parametric}
    J_r(\theta) =
    \int_{\mathcal{M}}\int_{\Omega}
    \big| \mathcal{L}\big(u_\theta(\bm{x},\bm{\mu});\bm{x},\bm{\mu}\big)-f(\bm{x},\bm{\mu}) \big|^2 \, \mathrm{d}\bm{x} \, \mathrm{d}\bm{\mu},
\end{equation}
the boundary loss function
\begin{equation}\label{eq:boundary_loss_parametric}
    J_b(\theta) =
    \int_{\mathcal{M}}\int_{\partial \Omega} \big| u_\theta(\bm{x},\bm{\mu}) - g(\bm{x},\bm{\mu}) \big|^2 \, \mathrm{d}\bm{x} \, \mathrm{d}\bm{\mu},
\end{equation}
and the data loss function
\begin{equation}\label{eq:data_loss_parametric}
    J_\text{data}(\theta) =
    \frac 1 {N_\text{data}} \sum_{i=1}^{N_\text{data}} \big| u_\theta\big(\bm{x}_\text{data}^{(i)},\bm{\mu}_\text{data}^{(i)}\big) - u_\text{data}^{(i)} \big|^2,
\end{equation}
where \smash{${\big(\bm{x}_\text{data}^{(i)}, \bm{\mu}_\text{data}^{(i)}, u_\text{data}^{(i)}\big)}_{i=1,\dots,N_\text{data}}$} are $N_\text{data}$ known data points, with $u_\text{data}^{(i)}$ a reference solution at points \smash{$\bm{x}_\text{data}^{(i)}$} and for given parameters {$\bm{\mu}_\text{data}^{(i)}$}. These reference solutions can be the exact solutions of the parametric PDE, defined in this case by \smash{$u_\text{data}^{(i)}=u\big(\bm{x}_\text{data}^{(i)};\bm{\mu}_\text{data}^{(i)}\big)$}. They can also be an approximation produced by a numerical method, such as finite elements on a fine mesh.

\begin{rmrk}
    In \cref{sec:numerical_results}, the focus is on PINNs trained only with a residual loss function (with boundary conditions imposed exactly as presented in \cref{sec:exact_imposition_of_BC}). We will only consider the BC loss function in \cref{sec:Lap2Dlowbc}. Furthermore, we will not use the data loss function in PINNs except in \cref{sec:Lap1D} where we will seek to compare a true physics-informed network to a merely data-driven one.
\end{rmrk}

Solving the minimization problem~\eqref{eq:minimization_problem}
requires computing the gradient of the loss function with respect to~$\theta$, which involves calculating the integrals
in~\eqref{eq:residual_loss_parametric} and~\eqref{eq:boundary_loss_parametric}.
The most natural idea is to estimate them with a Monte-Carlo method, see e.g.~\cite{Caf1998}.
One could also use Gauss-type quadrature rules to evaluate integrals, as is done in Variational Physics-Informed
Neural Networks~\cite{KhaZhaKar2021}, but the limitation is the impossibility of selecting an adequate
quadrature order due to the unknown properties of the \RTwo{neural network} approximation.
For that purpose, we define so-called ``collocation points'' on $\Omega\times\mathcal{M}$ and its boundary $\partial\Omega\times \mathcal{M}$,
denoted respectively by \smash{${\big(\bm{x}_\text{col}^{(i)}, \bm{\mu}_\text{col}^{(i)}\big)}_{i=1,\dots,N_\text{col}}$} and \smash{${\big(\bm{x}_\text{bc}^{(i)}, \bm{\mu}_\text{bc}^{(i)}\big)}_{i=1,\dots,N_\text{bc}}$}.
Then, we approximate the residuals and boundary losses by
\begin{equation}\label{eq:residual_loss_parametric_MC}
    J_r(\theta) \simeq
    \frac{1}{N_\text{col}} \sum_{i=1}^{N_\text{col}} \big| \mathcal{L}\big(u_\theta(\bm{x}_\text{col}^{(i)},\bm{\mu}_\text{col}^{(i)});\bm{x}_\text{col}^{(i)},\bm{\mu}_\text{col}^{(i)}\big)-f(\bm{x}_\text{col}^{(i)},\bm{\mu}_\text{col}^{(i)})  \big|^2
\end{equation}
and
\begin{equation}\label{eq:boundary_loss_parametric_MC}
    J_b(\theta) \simeq
    \frac{1}{N_\text{bc}} \sum_{i=1}^{N_\text{bc}} \big| u_\theta\big(\bm{x}_\text{bc}^{(i)},\bm{\mu}_\text{bc}^{(i)}\big) - g\big(\bm{x}_\text{bc}^{(i)},\bm{\mu}_\text{bc}^{(i)}\big) \big|^2,
\end{equation}
\RTwo{where} $N_{\mathrm{col}}$ and $N_\mathrm{bc}$
are heuristically determined and should be large enough to ensure that the Monte-Carlo integration is accurate enough.
The precise values of these parameters will be given in the numerical experiments. In the case of complex geometries, one solution for obtaining a sample of \RTwo{these collocation} points in the $\Omega$ domain is to use a level-set function, denoted $\varphi$. This function, which vanishes on the boundary of $\Omega$, can be obtained differently. The authors of~\cite{Sukumar_2022} propose different approaches to obtain a level-set function analytically in the case of polygonal or curved geometries. Learning-based approaches have also been proposed, notably in e.g.~\cite{park2019deepsdflearningcontinuoussigned,sitzmann2020implicitneuralrepresentationsperiodic}.

Because of the minimization problem~\eqref{eq:minimization_problem}, the PINN $u_\theta$ does not exactly satisfy the boundary conditions.
Moreover, loss functions compete, which may require fine-tuning the coefficients between $J_r$ and $J_b$.
In addition, classical PINNs do not include information on higher-order derivatives.
As highlighted in \cref{rmk:C_gain_additif,rmk:C_gain_multiplicatif}, for our purposes, a good prior should yield a good approximation of the derivatives of the solution.
For these reasons, the following section recalls several improvements of classical PINNs in the literature.

\subsection{Improving PINN training and prediction}\label{sec:improve_PINNs}

This section focuses on several ways of improving PINNs:
exactly imposing the boundary conditions in \cref{sec:exact_imposition_of_BC},
adding a higher-order derivative term in the loss function in \cref{sec:sobolev_training},
and countering the spectral bias in \cref{sec:spectral_bias}.
Although these approaches are presented separately,
they can easily be combined with one another.

\subsubsection{Exact imposition of boundary conditions}\label{sec:exact_imposition_of_BC}

To avoid the issues of classical PINNs discussed in \cref{sec:PINNs_parametric_PDE},
the authors of~\cite{LagLikFot1998,FraMicNav2024}
propose a method to enforce inhomogeneous Dirichlet boundary conditions exactly.
To that end, they search the approximation $u_{\theta}$ of solution to~\eqref{eq:parametric_PDE} with the form:
for all $\bm{x} \in \Omega$ and $\bm{\mu} \in \mathcal{M}$
\begin{equation*}\label{eq:prior_with_levelset}
    u_{\theta}(\bm{x},\bm{\mu}) = \varphi(\bm{x}) w_{\theta}(\bm{x},\bm{\mu}) + g(\bm{x},\bm{\mu}),
\end{equation*}
where $\varphi$ and $w_\theta$ are, respectively, the level-set function and a neural network as defined in \cref{sec:PINNs_parametric_PDE}.
Thus $u_{\theta}$ will automatically satisfy the boundary conditions, since
$u_{\theta}(\bm{x},\bm{\mu}) = g(\bm{x},\bm{\mu})$ for all $\bm{x} \in \partial \Omega$ and $\bm{\mu} \in \mathcal{M}$.
Note that this level-set function can be used in a few different ways, firstly to generate a sample of points in $\Omega$, as shown in \cref{sec:PINNs_parametric_PDE}, and secondly to impose boundary conditions. However, to use it directly in the formulation of the prior, it will require a certain regularity. For example the signed distance function is not a good candidate. \RTwo{Similar methods to impose boundary conditions exist for Robin and Neumann conditions; see~\cite{Sukumar_2022}.}

In this case, only the residual and data loss functions are minimized,
and the minimization problem~\eqref{eq:minimization_problem} becomes
\begin{equation*}
    \theta^\star = \argmin_{\theta}
    \big( \omega_r J_r(\theta) + \omega_\text{data} J_\text{data}(\theta) \big),
\end{equation*}
with $\omega_r$ and $\omega_\text{data}$ some weights to balance the terms of the loss function.

\subsubsection{Sobolev training for PINNs}\label{sec:sobolev_training}

As presented in \cref{sec:PINNs_parametric_PDE}, PINNs approximate the PDE solution by directly incorporating the equations into their training. Despite their effectiveness, these models can sometimes struggle to learn correctly, especially when the solution or its derivatives are complicated.
The authors of~\cite{son2021sobolevtrainingphysicsinformed} have proposed an approach called Sobolev training to try and overcome these difficulties\ROne{, in the case of solutions with high regularity. Note that, in this case, the source term needs to be differentiable to ensure sufficient regularity and to compute the additional term.} This method simply imposes constraints not only on the solutions themselves, but also on their derivatives.
In the context of solving the problem~\eqref{eq:parametric_PDE} under consideration, Sobolev training is applied by adding a cost term $J_\text{\rm sob}$ to the initial minimization problem~\eqref{eq:minimization_problem}:
\begin{equation}\label{eq:minimization_problem_sobolev}
    \theta^\star = \argmin_{\theta}
    \big( \omega_r J_r(\theta) + \omega_\text{\rm sob} J_\text{\rm sob}(\theta) + \omega_b J_b(\theta) + \omega_\text{data} J_\text{data}(\theta) \big),
\end{equation}
with $J_r$, $J_b$ and $J_\text{data}$ defined as in~\eqref{eq:residual_loss_parametric},~\eqref{eq:boundary_loss_parametric} and~\eqref{eq:data_loss_parametric} respectively and $\omega_r$, $\omega_\text{\rm sob}$, $\omega_b$ and $\omega_\text{data}$ the weights to balance the different terms of the loss function.
The Sobolev loss function $J_\text{\rm sob}$ in~\eqref{eq:minimization_problem_sobolev} is defined by
\begin{equation}\label{eq:sobolev_loss}
    J_\text{\rm sob}(\theta) = \int_{\mathcal{M}}\int_{\Omega} \big|\nabla_{\bm{x}}\big( \mathcal{L}\big(u_\theta(\bm{x},\bm{\mu});\bm{x},\bm{\mu}\big)-f(\bm{x},\bm{\mu})\big)\big|^2
    \, \mathrm{d}\bm{x} \, \mathrm{d}\bm{\mu},
\end{equation}
where the integral is estimated by the Monte-Carlo method,
similarly to the other loss functions.

\begin{rmrk}\label[rmrk]{rmk:sob_training}
    \RTwo{Please note that this training has an additional cost, due to the calculation of higher-order derivatives. For example, considering exactly the network from the first 2D test case \cref{sec:Lap2Dlow} over 1000 epochs with Adam, Sobolev training takes two times longer than standard training (79 seconds versus 38 seconds).}
\end{rmrk}

\subsubsection{Overcoming the spectral bias}\label{sec:spectral_bias}

Multiple ways of overcoming the spectral bias of MLPs are available.
For instance, in~\cite{TanSri2020}, the authors introduce Fourier features to improve the network, while the authors of~\cite{DolHeiMisMos2024} rely on a domain decomposition-based approach.

In this work, when dealing with high-frequency solutions (i.e., solutions with more than three wavelengths propagating), we use the Fourier features from~\cite{TanSri2020}.
It relies on modifying the input of the neural network.
Indeed, the prior $u_\theta$ is now defined, for all $\bm{x} \in \Omega$ and $\bm{\mu} \in \mathcal{M}$, as
\begin{equation*}
    u_\theta(\bm{x},\bm{\mu}) =
    w_\theta\big(
    \bm{x},\bm{\mu};
    \sin(\pi a_1 \bm{x}),
    \cos(\pi b_1 \bm{x}),
    \dots,
    \sin(\pi a_{n_f} \bm{x}),
    \cos(\pi b_{n_f} \bm{x})
    \big),
\end{equation*}
\RTwo{where $(a_i)_i \in \mathbb{R}^{n_f}$ and $(b_i)_i \in \mathbb{R}^{n_f}$ are additional trainable parameters.}
This makes it possible to learn higher-frequency solutions, by also learning the frequency itself. This MLP with Fourier features (MLP w/ FF) needs the same parameters as the classical MLP (defined in \cref{rmk:PINN_notations}), but also the number $n_f$ of Fourier features.

    \section{Implementation details}\label{sec:implementation_details}
    Before moving on to numerical experiments,
we discuss some practical details
regarding the implementation of
the methods introduced in \cref{sec:additive_prior,sec:multiplicative_prior}.
In \cref{sec:using_PINN}, we first look at how to effectively plug the PINN prediction in the FEM solver.
Then, in \cref{sec:boundary_conditions}, we discuss the imposition of boundary conditions in the two proposed methods. 

The tools used to implement the methods and obtain the numerical results in \cref{sec:numerical_results} are, on the one hand, \texttt{PyTorch}~\cite{paszke2019pytorchimperativestylehighperformance} and \texttt{ScimBa}\footnote{\url{https://gitlab.inria.fr/scimba/scimba}} for the prior construction, in particular the implementation of the PINN, and, on the other hand, \texttt{FEniCS}~\cite{AlnBle2015} \RTwo{(version 2019.1.0) or \texttt{FEniCSx}~\cite{baratta_dolfinx_2023,scroggs_construction_2022,scroggs_basix_2022,alnaes_unified_2014} (version 0.8)} for the implementation of the finite element methods. For mesh generation, use either \texttt{FEniCS}\RTwo{/\texttt{FEniCSx}} or \texttt{mshr} mesh generators. \RTwo{Note that, we define the characteristic mesh size $h$ as the length of the longest edge, for the space dimension $d\in\{1,2,3\}$.}

\RTwo{Moreover,} in \cref{sec:numerical_results}, we do not specify the training times of the networks for each test case. To give the reader a rough idea, PINN training takes less than ten minutes on a \RTwo{NVIDIA RTX 2000 (8GB VRAM)}. For cases using the LBFGS optimizer, training takes a little longer, but remains under an hour. \RTwo{Regarding the EF part, the CPU considered depends on the numerical cost required for each test case; we ran the solvers using either Intel I7-13800H or AMD EPYC 7713.}

\begin{rmrk}\label[rmrk]{rmk:PINN_notations}
	To construct the PINN, we require some hyperparameters for the MLP and the training phase.
    In particular, the MLP activation function is denoted by $\sigma$, and ``\emph{layers}'' represents a sequence of integers describing the number of neurons associated with each layer of the MLP.
    For training, we will denote by ``\emph{lr}'' the learning rate and $n_\text{epochs}$ the number of epochs considered, as well as ``\emph{decay}'' the multiplicative factor of the learning rate decay considered every $20$ epochs thanks to \texttt{PyTorch}'s \texttt{StepLR} scheduler.
    Unless otherwise specified, the batch size will correspond to the number of collocation points chosen.
    The Adam optimizer~\cite{KinBa2015} will be used to train the network, but in some cases, we will switch to the LBFGS optimizer~\cite{nocedal_quasi_newton_2006} at the $n_\text{switch}$-th epoch.
\end{rmrk}

\subsection{Using PINN prediction effectively}\label{sec:using_PINN}

To be effective, our methods will in practice depend on the quality of the approximation of the prior's derivatives, computed from automatic differentiation, and its precise integration on the domain.

\boldparagraph{Automatic differentiation}

It is important to use the automatic differentiation offered by neural networks, enabling exact (in the sense of machine precision) derivative computation without having to manipulate complex symbolic expressions. In particular, in the context of PINNs, automatic differentiation will play a fundamental role in integrating the PDE under consideration. This automatic differentiation will enable the two improved finite element methods to use the exact derivatives of the prior $u_\theta$ and thus avoid introducing an additional error in the computation of the derivative.

\boldparagraph{Numerical integration}

\RTwo{Our approaches require the numerical integration, in the weak problem, of several functions with a closed-form expression (most notably the prior; also the right-hand side of the PDE for instance). This integration has to be done with sufficient precision for our methods to be effective.}
Thus, in e.g.\ the additive approach, a quadrature rule with a higher degree than the traditional FEM has to be applied to discretize the term $l(v_h) - a(u_\theta,v_h)$ in~\eqref{eq:approachform_add}.
This point, and the required degree of the quadrature rule,
will be studied in more detail in the first 2D test case considered in \cref{sec:Lap2Dlow}.

\begin{rmrk}
In practice, the source term in the additive approach will be computed in the strong way. For instance, in the case of the Laplacian equation with $u_{\theta}=0$ on $\partial\Omega$, the term
\[
    l(v_h)-a(u_{\theta},v_h)
    =
    \int_{\Omega}f\RTwo{(\bm{x})} v_h\RTwo{(\bm{x})} \, \RTwo{\mathrm{d}\bm{x}}
    -
    \int_{\Omega}\nabla u_\theta\RTwo{(\bm{x})} \cdot \nabla v_h\RTwo{(\bm{x})} \, \RTwo{\mathrm{d}\bm{x}}
\]
will be replaced by
\[
    \int_{\Omega}\big(f\RTwo{(\bm{x})}+\Delta u_{\theta}\RTwo{(\bm{x})}\big)v_h\RTwo{(\bm{x})} \, \RTwo{\mathrm{d}\bm{x}}.
\]
If $u_{\theta}$ is not equal to zero on $\partial\Omega$, one needs to include a boundary term.
\end{rmrk}

\subsection{Imposing boundary conditions}\label{sec:boundary_conditions}

In this section, we focus on the crucial question of imposing boundary conditions. We first look at this problem in the context of the additive approach presented in \cref{sec:additive_prior} and then in the context of the multiplicative approach presented in \cref{sec:multiplicative_prior}.

For simplicity, this section focuses on (non-homogeneous) Dirichlet conditions.
For our enriched FEM, just like in classical FEM,
these boundary conditions are imposed by manually eliminating essential dofs, see~\cite{Ern2004TheoryAP}, and more precisely by modifying the matrix and the right-hand side of the linear system. This approach is not needed for Neumann and Robin conditions.

\subsubsection{Additive approach}

In this first approach, if our Dirichlet problem satisfies
$u=g$ on $\partial \Omega$,
then $p_h^+$ has to satisfy
\[
    p_h^+ = g - u_{\theta} \text{\quad on } \partial \Omega,
\]
with $u_\theta$ the PINN prior.
This non-homogeneous boundary condition becomes homogeneous as soon as $u_\theta$ is exact at the boundary, or, in other words, as soon as the boundary conditions are imposed exactly in the PINN, as presented in \cref{sec:exact_imposition_of_BC}.

\begin{rmrk}\label[rmrk]{rem:bconcurved}
    However, in the case of curved geometries (e.g.\ disks) where the meshes do not coincide with the boundary of the geometry, problems occur when $k>2$, and especially on coarse meshes. In the numerical results, to avoid this problem and check the error estimates, we assume that $g=u$ on $\partial\Omega_h$ where $\Omega_h$ is the domain covered by the mesh. Furthermore, we need to be careful because even if, in PINN, the conditions are imposed exactly (as shown in \cref{sec:exact_imposition_of_BC}), the prediction $u_\theta$ will not be exact on \RTwo{the approximate boundary} $\partial\Omega_h$. We made this choice here to simplify the problem, but in practice, solutions exist to improve the quality of the results.
\end{rmrk}

\subsubsection{Multiplicative approach}\label{sec:multiplicative_BC}

In this second method, the boundary conditions are a bit more complex to handle. In the \cref{sec:multiplicative_prior}, we have denoted by
\[
    u_{h,M}^\times = u_{\theta,M} \; p_h^\times
\]
the solution obtained by the multiplicative approach to the modified problem~\eqref{eq:ob_pde_M} with the prior $u_{\theta,M}=u_\theta+M$. Therefore, we can recover the solution $u_h^\times$ of the original problem~\eqref{eq:ob_pde} by setting $u_h^\times = u_{h,M}^\times - M$.

\boldparagraph{Standard PINN}

In the case where our prior $u_\theta$ is the prediction of a standard PINN, as presented in \cref{sec:PINNs_parametric_PDE}, the boundary conditions are not imposed exactly. Thus, if our problem satisfies $u=g$ on $\partial \Omega$, then $p_h^\times$ has to satisfy
\[
    p_h^\times = \frac{g+M}{u_{\theta,M}} \text{\quad on } \partial \Omega.
\]

\boldparagraph{PINN with exact BC}

We now tackle the case where we exactly impose the boundary conditions in the PINN, as presented in \cref{sec:exact_imposition_of_BC}.
Then, supposing that $M>0$, we have $u_{\theta,M}=g+M$ on $\partial \Omega$ and therefore~$p_h^\times$ has to satisfy
\[
    p_h^\times = 1 \text{\quad on } \partial \Omega.
\]
However, there is a specific case when this condition is not necessarily true. Indeed, if the boundary conditions are homogeneous, then $g=0$ and $u_{\theta,M}=M$ on $\partial \Omega$.
Considering that $u_\theta>0$ in $\Omega$,
$M=0$ is a possible choice.
In this case, \smash{$u_{\theta,M}=u_{\theta,0}=0$} on $\partial \Omega$,
and \smash{$u_{h,M}^\times=u_{h,0}^\times=u_h^\times$} automatically satisfies the boundary conditions.
Hence, imposing a boundary condition on $p_h^\times$ becomes unnecessary.

\begin{rmrk}
    In the numerical results of \cref{sec:Ell1D}, one of these specific
    cases is considered in 1D.
    We will see that leaving $p_h^\times$ free will give better results here than imposing $p_h^\times=1$ on $\partial \Omega$. Indeed, this approach leaves more freedom to capture the correct derivatives at the boundary.
\end{rmrk}

    \section{Numerical results}\label{sec:numerical_results}
    This section is dedicated to validating the proposed method on several numerical experiments, which are \RBothN{mostly} instances of the \RBoth{problem~\eqref{eq:ob_pde}} with spatial dimension $d \in \{1,2\RTwo{,3}\}$, and with increasing complexity.
The idea is to compare, in different ways, the additive and multiplicative approaches presented in \cref{sec:additive_prior} and \cref{sec:multiplicative_prior} to the standard finite element method presented in \cref{sec:FEM}.
The multiplicative approach will only be considered in dimension $d=1$, showing that only in rare cases, and with a good choice of the lifting constant $M$, does it provide better results than the additive one.
We will also show that the approach proposed in \cref{sec:prior_construction} for choosing the prior is more efficient than more classical ones, still in the 1D case.

In \cref{sec:setup}, we present the two tests that will be performed for each test case. We are interested in \RBothN{three} 1D test cases ($d=1$): the Poisson problem with homogeneous Dirichlet conditions in \cref{sec:Lap1D}, a general elliptic system in a convection-dominated regime in \cref{sec:Ell1D} \RBothN{and a non-smooth transmission problem in \cref{sec:SingLap1D}}. We then consider three 2D cases ($d=2$). We start with a Poisson problem with homogeneous Dirichlet conditions on a square domain in \cref{sec:Lap2D}. We then continue with a more complicated elliptic problem, still with Dirichlet conditions and on a square domain, but with parameter-dependent anisotropy, in \cref{sec:Ell2D}. In \cref{sec:Lap2DMixRing}, we return to the 2D Poisson problem, but this time considering mixed boundary conditions on a ring-shaped domain. \RTwo{Finally, we consider a 3D Poisson problem with homogeneous Dirichlet conditions on a cube domain in \cref{sec:Lap3D}.}
To reproduce these results, an open-source code is available on GitHub\footnote{{\url{https://github.com/flecourtier/EnrichedFEMUsingPINNs}}}.

\RTwo{More specifically, in each test case, we will consider the additive approach and use a prior constructed by a PINN (where we \RBothN{mostly} impose boundary conditions using a level-set function). We will perform an analysis of error estimates and an evaluation of the gains of the enriched approach over the standard approach on a set of parameters (in $L^2$ norm). In \RBothN{the first two} 1D test cases (\cref{sec:Lap1D,sec:Ell1D}), we will also choose to consider the multiplicative approach. A study in the $H^1$ semi-norm will be performed in the first 1D and 2D test case (\cref{sec:Lap3D,sec:Lap2Dlow}). Additionally, this first 1D test case includes a comparative study between theoretical and numerical constants. \RBothN{In the third 1D test case, we will focus on a problem where the gradient is discontinuous.} Subsequently, we will test the approach with other networks, notably: a data-driven network (\cref{sec:Lap1D}) that will justify the use of PINNs, \RBothN{networks} with weakly imposed conditions (\RBothN{\cref{sec:SingLap1D,sec:Lap2Dlowbc}}), and an improved PINN where we used Sobolev training (\cref{sec:Lap2Dlowaug}). Two studies on the numerical costs of the additive approach will be carried out in \cref{sec:Lap2Dlow,sec:Lap3D_costs}. We will also perform a brief analysis of the influence of PINN quality, depending on the initialization of the network weights, the number of collocation points considered, and the number of epochs (\cref{sec:Lap1D_hyperparam,sec:Lap3D_badPINN}). In some test cases, we will also look at the visualization of the solutions obtained by the different FE approaches.}

\RTwoN{The first two 1D test cases show that the multiplicative approach only gives better results than the additive approach in certain very specific cases. The first test case also validates the theoretical results, in particular the importance of the derivatives of the prior by validating the choice of PINNs. The third 1D test case extends the theoretical framework to a non-smooth problem, showing that, once again, high-order derivatives have the greatest impact on the results. The first 2D test case (Poisson in a square) introduces training with Sobolev loss and discusses the numerical costs of FE approaches (in terms of number of degrees of freedom). It also discusses the imposition of weak or strong boundary conditions in the network and compares the gains between the “high frequency” and “low frequency” cases. The second 2D test case highlights the flexibility of the method for more complex problems, particularly with parameter-dependent anisotropy. The third 2D test case introduces the approach in a slightly more complex geometric setting (with an annulus) as well as Robin conditions. Finally, the 3D test case validates the approach in a larger dimension and analyzes the associated numerical costs. This time in terms of computation time, considering the complete parametric framework. It also introduces the influence of PINN quality on the final results. More generally, all test cases show that, with the right choice of priors, the additive approach yields significant gains over the standard approach. In particular, the average results obtained over a set of parameters show that for elements $\mathbb{P}_1$, this enriched approach provides the same accuracy as the standard approach with meshes $5$ (2D anisotropic problem) to $28$ (low-frequency Poisson with Sobolev training) times coarser, which leads to a gain in computation time (from a purely online perspective). The 3D test case takes a closer look at execution times, particularly in the parametric context in which the enriched approach is developed. More specifically, gains in terms of mesh size make it possible to achieve a speed-up from 19 parameter sets to reach a relative accuracy of $10^{-3}$ (in the $L^2$ norm).}

\subsection{Setup of the numerical experiments}\label{sec:setup}

For each of the proposed test cases,
we consider a parametric problem, on which we train a PINN as presented in \cref{sec:PINNs_parametric_PDE},
to resolve it on a set of parameters, denoted by $\mathcal{M}$.
Let $p=\dim(\mathcal{M})$ be the number of parameters,
and consider a \RTwo{set $\mathcal{S}$ of $n_p$ parameter instances}:
\begin{equation*}
    \mathcal{S}=\left\{\bm{\mu}^{(1)},\dots,\bm{\mu}^{(n_p)}\right\},
\end{equation*}
with, for $j=1,\dots,n_p$,
\begin{equation*}
    \bm{\mu}^{(j)}=\left(\mu_1^{(j)},\dots,\mu_{p}^{(j)}\right)\in \mathcal{M}.
\end{equation*}

In the following, we denote by $u^{(j)}$ a reference solution
to problem~\eqref{eq:parametric_PDE} for a given parameter $\bm{\mu}^{(j)}$
and by \smash{$u_h^{(j)}$} the solution obtained by
the standard finite element method~\eqref{eq:approachform},
where $V_h$ is the $\mathbb{P}_k$ Lagrange space defined in~\eqref{eq:Vh}
and $h$ is the characteristic mesh size.
We also denote by \smash{$u_\theta^{(j)}$} the solution obtained by the parametric PINN,
and by \smash{$u_{h,+}^{(j)}$} the solution obtained by the additive approach~\eqref{eq:approachform_add} with $V_h^+$ the $\mathbb{P}_k$ Lagrange space defined in~\eqref{eq:Vh_add}.
In some test cases, we also consider the solution \smash{$u_{h, M}^{(j)}$} of the multiplicative approach~\eqref{eq:approachform_mul} with $V_h^\times$ the $\mathbb{P}_k$ Lagrange space defined in~\eqref{eq:Vh_mul}, depending on the lifting constant $M$.

\begin{rmrk}\label[rmrk]{rk:ref}
    In the following, to estimate the error, we consider the reference solution to be either an analytical solution or a solution obtained with a very fine mesh and a high polynomial degree. More precisely, we need the characteristic size $h_\text{ref}$ associated with the reference mesh to be much smaller than the size associated with the current mesh $h$, i.e. $h_\text{ref}\ll h$, and we will consider $k_\text{ref}=3$ the polynomial degree associated with the reference solution.
\end{rmrk}

In each test case, we investigate two aspects; the first in~\cref{sec:setup_error_estimates} involves verifying the error estimates and the second in~\cref{sec:setup_gains} is the evaluation of the gains achieved by the proposed methods compared with the standard one.

\subsubsection{Error estimates}\label{sec:setup_error_estimates}

Consider a small \RTwo{set $\mathcal{S}$ with $n_p=2$ parameter instances}. Given a fixed parameter $\bm{\mu}^{(j)}$, $j=1,2$, we start by testing the error estimates obtained in \cref{lem:error_estimation_add} for the additive approach. In the case $d=1$, we will also be interested in the error estimates in \cref{lem:error_estimate_multiplicative} for the multiplicative approach. By varying the mesh size~$h$, we then estimate the errors obtained with the two methods.
To evaluate these errors, we compare the approximations to the reference solution $u^{(j)}$ (see \cref{rk:ref}).
We then define by
\begin{equation}\label{eq:error_rel_FEM}
    e_h^{(j)}=\frac{||u^{(j)}-u_h^{(j)}||_{L^2}}{||u^{(j)}||_{L^2}}
    \text{ and }
    e_\theta^{(j)}=\frac{||u^{(j)}-u_\theta^{(j)}||_{L^2}}{||u^{(j)}||_{L^2}},
\end{equation}
the $L^2$ relative error obtained for the standard FEM and the PINN respectively. We further define,
\begin{equation}\label{eq:error_rel_add}
    e_{h,+}^{(j)}=\frac{||u^{(j)}-u_{h,+}^{(j)}||_{L^2}}{||u^{(j)}||_{L^2}}
    \text{ and }
    e_{h,M}^{(j)}=\frac{||u^{(j)}-u_{h,M}^{(j)}||_{L^2}}{||u^{(j)}||_{L^2}},
\end{equation}
the $L^2$ relative errors obtained for the additive and multiplicative approach (depending on the lifting constant $M$), respectively. \ROne{Relative errors in the semi-norm $H^1$ can be defined in the same way.}

\subsubsection{Gains achieved with the enriched bases}\label{sec:setup_gains}

As we have trained the network to be parameter-dependent to predict a solution for a set of parameters, we are interested in the average gains we obtain with our enriched approaches compared to the PINN and the standard FEM. More precisely, for a fixed mesh size $h$ and a fixed polynomial degree $k$, for a \RTwo{set $\mathcal{S}$ of $n_p$ parameter instances}, the numerical gains obtained by the additive approach on PINN and standard FEM are respectively defined for $j=1,\dots,n_p$ by:
\begin{equation}\label{eq:gain_j}
    G_{+,\theta}^{(j)}=\frac{e_\theta^{(j)}}{e_{h,+}^{(j)}}
    \text{\quad and \quad}
    G_+^{(j)}=\frac{e_h^{(j)}}{e_{h,+}^{(j)}}, 
\end{equation}
with $e_\theta^{(j)}$, $e_h^{(j)}$ and $e_{h,+}^{(j)}$ respectively the $L^2$ relative errors obtained with the PINN, the standard FEM and the additive approach, defined in \cref{sec:setup_error_estimates}.
Similarly, the theoretical gains obtained by the multiplicative approach (depending on the lifting constant $M$) on PINN and standard FEM are respectively defined for $j=1,\dots,n_p$ by:
\begin{equation}\label{eq:gain_j_mul}
    G_{M,\theta}^{(j)}=\frac{e_\theta^{(j)}}{e_{h,M}^{(j)}} \quad \text{and} \quad G_M^{(j)}=\frac{e_h^{(j)}}{e_{h,M}^{(j)}}, 
\end{equation}
with $e_{h,M}^{(j)}$ the $L^2$ relative error obtained with the multiplicative approach (depending on the lifting constant $M$), defined in \cref{sec:setup_error_estimates}.
Therefore, we will be interested in the minimum, maximum, mean and standard deviation obtained on the following samples:
\begin{equation}\label{eq:gain_add_num}
    G_{+,\theta}=\left\{G_{+,\theta}^{(1)},\dots,G_{+,\theta}^{(n_p)}\right\} \quad \text{and} \quad G_+=\left\{G_+^{(1)},\dots,G_+^{(n_p)}\right\}, 
\end{equation}
which respectively represent the gains obtained with our additive approach over PINN and over standard FEM on the \RTwo{set $\mathcal{S}$ of parameter instances}. In the same way, we define $G_{M,\theta}$ and $G_M$, which respectively represent the gains obtained with our multiplicative approach over PINN and over standard FEM on \RTwo{$\mathcal{S}$} by:
\begin{equation}\label{eq:gain_mul_num}
    G_{M,\theta}=\left\{G_{M,\theta}^{(1)},\dots,G_{M,\theta}^{(n_p)}\right\} \quad \text{and} \quad G_M=\left\{G_M^{(1)},\dots,G_M^{(n_p)}\right\}. 
\end{equation}

\subsection{1D Poisson problem}\label{sec:Lap1D}

In this section, we will consider the problem~\eqref{eq:ob_pde} in its most
simplified Poisson form, with homogeneous Dirichlet boundary conditions.
In the 1D case ($d=1$), we have,
\begin{equation}\label{eq:Lap1D}
    \left\{
    \begin{aligned}
        -\partial_{xx} u & = f, \; &  & \text{in } \; \Omega \times \mathcal{M},         \\
        u                & = 0, \; &  & \text{on } \; \partial\Omega \times \mathcal{M},
    \end{aligned}
    \right.
\end{equation}
with $\Omega=\RBothN{(0,1)}$, $\partial\Omega$ its boundary and $\mathcal{M} \subset \mathbb{R}^p$ the parameter space (with $p$ the number of parameters).
{\ROne{We prescribe a family of exact solutions depending on the parameter vector $\bm{\mu}=(\mu_1,\mu_2,\mu_3)\in\mathcal{M}={[0,1]}^3$ ($p=3$ parameters)} defined by
}
\begin{equation*}
    u(x,\bm{\mu})=\mu_1\sin(2\pi x)+\mu_2\sin(4\pi x)+\mu_3\sin(6\pi x) \,.
\end{equation*}
{Note that the associated right-hand side $f$ in~\eqref{eq:Lap1D}
also depends on $\bm{\mu}$.
This problem is thus four-dimensional: one dimension in space and three dimensions
for the parameters $\bm{\mu}$.}

For this first test case, we construct two priors,
as detailed in \cref{sec:Lap1D_priors}.
The first one, denoted by $u_\theta$,
is built from a PINN as presented in \cref{sec:prior_construction}.
The second one, denoted by $u_\theta^\text{data}$
is constructed only from data (obtained from the analytical solution).
The aim is to show that using physics-informed training to construct the prior
leads to better results than data-driven training.
In this test case, we also compare the additive and multiplicative approaches (presented for several values of $M$), but only with $k=1$ polynomial order to remain concise.

We first present the error estimates obtained with the PINN prior in \cref{sec:Lap1D_error_estimations}. First, we check the orders of convergence of the two enriched approaches. Then, we verify the results expected in \cref{sec:comparison_add_mul} by comparing the theoretical gain constants of the additive and multiplicative approaches. Afterwards, in \cref{sec:Lap1D_derivatives}, we compare, for a given parameter, the derivatives obtained with the two priors and analyze the associated gains. \RTwo{Then}, we evaluate the gains obtained with the two priors in \cref{sec:Lap1D_gains} on a sample of parameters. \RTwo{Finally, in \cref{sec:Lap1D_hyperparam}, we analyze the influence of the weight initialization of the PINN prior on the gains obtained with our enriched approaches.}

\begin{rmrk}\label[rmrk]{rmk:N_nodes}
    In the following, the characteristic mesh size is $h=\frac{1}{N-1}$, where $N$ represents the number of considered nodes.
\end{rmrk}

\subsubsection{Construction of the two priors}\label{sec:Lap1D_priors}

The hyperparameters used to construct the two priors are presented in \cref{tab:paramtest1_1D}. We discuss below the specific differences in training
both priors.

\begin{table}[htbp]
    \centering
    \begin{tabular}{cc}
        \toprule
        \multicolumn{2}{c}{\textbf{Network - MLP}} \\
        \midrule
        \textit{layers} & $20,80,80,80,20,10$ \\
        \cmidrule(lr){1-2}
        $\sigma$ & sine \\
        \bottomrule
    \end{tabular}
    \hspace{0.3cm}
    \begin{tabular}{cccc}
        \toprule
        \multicolumn{4}{c}{\textbf{Training}} \\
        \midrule
        \textit{lr} & 9e-2 & $n_{epochs}$ & \num{10000} \\
        \cmidrule(lr){1-2} \cmidrule(lr){3-4}
        \textit{decay} & 0.99 \\
        \cmidrule(lr){1-2}
        $N_{\text{col}/\text{data}}$ & \num{5000} \\
        \bottomrule
    \end{tabular}
    \hspace{0.3cm}
    \begin{tabular}{cccccccc}
        \toprule
        \multicolumn{8}{c}{\textbf{Loss weights}} \\
        \midrule
        \multicolumn{4}{c}{PINN prior $u_\theta$} & \multicolumn{4}{c}{Data prior $u_\theta^\text{data}$} \\
        \cmidrule(lr){1-4} \cmidrule(lr){5-8}
        $\omega_r$ & 1 & $\omega_\text{data}$ & 0 & $\omega_r$ & 0 & $\omega_\text{data}$ & 1 \\
        \cmidrule(lr){1-2} \cmidrule(lr){3-4} \cmidrule(lr){5-6} \cmidrule(lr){7-8}
        $\omega_b$ & 0 & $\omega_\text{sob}$ & 0 & $\omega_b$ & 0 & $\omega_\text{sob}$ & 0 \\        
        \bottomrule
    \end{tabular}
    \caption{Network, training parameters (\cref{rmk:PINN_notations}) and loss weights for $u_\theta$ and $u_\theta^\text{data}$ in the \textit{1D Poisson problem}. Considering $N_\text{col}$ collocation points for the PINN prior and $N_\text{data}$ data for the data prior.}\label{tab:paramtest1_1D}
\end{table}

\boldparagraph{Physics-informed training}
For the first prior, we will consider a parametric PINN, depending on the problem parameters $\bm{\mu}$,
where we exactly impose the Dirichlet boundary conditions as presented in \cref{sec:exact_imposition_of_BC} and without using data in training.
\RBoth{Thus, we construct $u_\theta$ as in \eqref{eq:prior_with_levelset} with the level-set function $\varphi$ defined by}
\[
    \varphi(x)=x(x-1) \,,
\]
which vanishes exactly on $\partial\Omega$ \RBoth{and $g=0$}.
Since we impose the boundary conditions by using the level-set function, we will only consider the residual loss $J_r$ \RBoth{approached by a Monte-Carlo method as defined in~\eqref{eq:residual_loss_parametric_MC} with $N_\text{col}=\num{5000}$ collocation points uniformly chosen on $\Omega\times\mathcal{M}$.}
The hyperparameters are given in \cref{tab:paramtest1_1D};
we use the Adam optimizer~\cite{KinBa2015}.

\boldparagraph{Data-driven training}

For the second prior considered, noted $u_\theta^\text{data}$, a network is trained only on the data (constructed from the analytical solution). \RBoth{Therefore, we will only consider the data loss $J_\text{data}$ defined in~\eqref{eq:data_loss_parametric}, considering $N_\text{data}=\num{5000}$ points. As for the physics-informed training, we still consider the hyperparameters defined in \cref{tab:paramtest1_1D}  and construct the prior $u_\theta^\text{data}$ in the same way using the level-set function.}

\subsubsection{Error estimates --- with the PINN prior}\label{sec:Lap1D_error_estimations}

In this section, we look at the theoretical results of the additive and multiplicative approaches, considering the PINN prior $u_\theta$. First, we check the orders of convergence of \cref{lem:error_estimation_add,lem:error_estimate_multiplicative} (in the $L^2$ norm), associated with both methods. Next, we numerically verify \cref{thm:comparison_add_mul}, showing that the multiplicative correction converges, for sufficiently large $M$, towards the additive correction (in both the $L^2$ norm and $H^1$ semi-norm).

\boldparagraph{Convergence rate}

We test the error estimates of \cref{lem:error_estimation_add,lem:error_estimate_multiplicative} for the following two sets of parameters:
\begin{equation*}
    \bm{\mu}^{(1)}=(0.3,0.2,0.1) \quad \text{and} \quad \bm{\mu}^{(2)}=(0.8,0.5,0.8) \,,
\end{equation*}
with the PINN prior $u_\theta$.
For $j \in \{1, 2\}$, the aim is to compare, by varying the mesh size $h$, the $L^2$ relative errors \smash{$e_h^{(j)}$} obtained with the standard FEM, defined in~\eqref{eq:error_rel_FEM}, \smash{$e_{h,+}^{(j)}$} obtained with the additive approach, defined in~\eqref{eq:error_rel_add} and \smash{$e_{h,M}^{(j)}$} obtained with the multiplicative approach (taking $M=3$ and $M=100$), defined in~\eqref{eq:error_rel_add}. The results are presented in \cref{fig:case1_1D}
for polynomial orders $k=1$ and $k=2$, with $h$ depending
on the number of nodes $N \in \{16,32,64,128,256\}$ as
presented in \cref{rmk:N_nodes}. \ROne{We also present in \cref{fig:case1_1DH1} the results obtained for the semi-norm $H^1$.}

\begin{figure}[ht!]
    \centering
    \begin{subfigure}{0.48\linewidth}
        \centering
        \cvgFEMCorrMultOnedeg{fig_testcase1D_test1_cvg_FEM_case1_v1_param1_degree1.csv}{fig_testcase1D_test1_cvg_FEM_case1_v1_param1_degree2.csv}{fig_testcase1D_test1_cvg_Corr_case1_v1_param1_degree1.csv}{fig_testcase1D_test1_cvg_Mult_case1_v1_param1_degree1_M3.0.csv}{fig_testcase1D_test1_cvg_Mult_case1_v1_param1_degree1_M100.0.csv}{1e-5}
        \caption{Case of $\bm{\mu}^{(1)}$}\label{fig:case1param1_1D}
    \end{subfigure}
    \begin{subfigure}{0.48\linewidth}
        \centering
        \cvgFEMCorrMultOnedeg{fig_testcase1D_test1_cvg_FEM_case1_v1_param2_degree1.csv}{fig_testcase1D_test1_cvg_FEM_case1_v1_param2_degree2.csv}{fig_testcase1D_test1_cvg_Corr_case1_v1_param2_degree1.csv}{fig_testcase1D_test1_cvg_Mult_case1_v1_param2_degree1_M3.0.csv}{fig_testcase1D_test1_cvg_Mult_case1_v1_param2_degree1_M100.0.csv}{5e-6}
        \caption{Case of $\bm{\mu}^{(2)}$}\label{fig:case1param2_1D}
    \end{subfigure}
    \caption{Considering the \textit{1D Poisson problem} and the PINN prior $u_\theta$. Left -- Considering $\bm{\mu}^{(1)}$. $L^2$ \ROne{relative} error on $h$ obtained with standard FEM \smash{$e_h^{(1)}$} (solid lines) with $k=1$ and $k=2$, the additive approach \smash{$e_{h,+}^{(1)}$} (dashed lines) with $k=1$ and the multiplicative approach \smash{$e_{h,M}^{(1)}$} (dotted lines) with $k=1$ ($M=3$ and $M=100$). Right -- Same for $\bm{\mu}^{(2)}$.}\label{fig:case1_1D}
\end{figure}

\begin{figure}[ht!]
    \centering
    \begin{subfigure}{0.48\linewidth}
        \centering
        \cvgFEMCorrMultOnedegHun{fig_testcase1D_test1_cvgsemiH1_FEM_case1_v1_param1_degree1.csv}{fig_testcase1D_test1_cvgsemiH1_FEM_case1_v1_param1_degree2.csv}{fig_testcase1D_test1_cvgsemiH1_Corr_case1_v1_param1_degree1.csv}{fig_testcase1D_test1_cvgsemiH1_Mult_case1_v1_param1_degree1_M3.0.csv}{fig_testcase1D_test1_cvgsemiH1_Mult_case1_v1_param1_degree1_M100.0.csv}{2.5e-4}
        \caption{Case of $\bm{\mu}^{(1)}$}\label{fig:case1param1_1DH1}
    \end{subfigure}
    \begin{subfigure}{0.48\linewidth}
        \centering
        \cvgFEMCorrMultOnedegHun{fig_testcase1D_test1_cvgsemiH1_FEM_case1_v1_param2_degree1.csv}{fig_testcase1D_test1_cvgsemiH1_FEM_case1_v1_param2_degree2.csv}{fig_testcase1D_test1_cvgsemiH1_Corr_case1_v1_param2_degree1.csv}{fig_testcase1D_test1_cvgsemiH1_Mult_case1_v1_param2_degree1_M3.0.csv}{fig_testcase1D_test1_cvgsemiH1_Mult_case1_v1_param2_degree1_M100.0.csv}{2e-4}
        \caption{Case of $\bm{\mu}^{(2)}$}\label{fig:case1param2_1DH1}
    \end{subfigure}
    \caption{\ROne{Considering the \textit{1D Poisson problem} and the PINN prior $u_\theta$. Left -- Considering $\bm{\mu}^{(1)}$. Semi-$H^1$ relative error on $h$ obtained with standard FEM (solid lines) with $k=1$ and $k=2$, the additive approach (dashed lines) with $k=1$ and the multiplicative approach (dotted lines) with $k=1$ ($M=3$ and $M=100$). Right -- Same for $\bm{\mu}^{(2)}$.}}\label{fig:case1_1DH1}
\end{figure}

The results of \ROne{\cref{fig:case1_1D,fig:case1_1DH1}} show that all the enriched finite elements increase the accuracy of the method and that they also converge at the same rate as the classical approach (i.e., \RTwo{for polynomial approximation of order $k = 1$, the convergence order in the $L^2$ norm is $2$ and in the semi-norm $H^1$ is $1$}).
Furthermore, the theoretical analysis (which showed that the multiplicative correction has the same error as the additive one when $M \to \infty$) is confirmed for both sets of parameters.
In addition, \ROne{\cref{fig:case1_1D,fig:case1_1DH1}} also shows that this multiplicative enrichment can be less efficient for small $M$ when the $(k+1)$\textsuperscript{th} derivative of the solution is large.
Indeed, for the first parameter considered in \ROne{\cref{fig:case1param1_1D,fig:case1param1_1DH1}}, for which the second derivative takes lower values, we observe that the multiplicative approach with small $M$ is closer to the additive one than for the second set of parameters
considered in \ROne{\cref{fig:case1param2_1D,fig:case1param2_1DH1}}, for which the derivatives are larger.
Moreover, it seems that we gain almost one order of interpolation with the additive approach: the additive method with polynomial order $k=1$ gives \ROne{relative errors} close to the original FEM with $k=2$, although the rate of convergence is different.

\boldparagraph{Gain constants}

\begin{figure}[ht!]
    \centering
    \begin{minipage}{0.55\linewidth}
        \centering
        \includegraphics[scale=0.9]{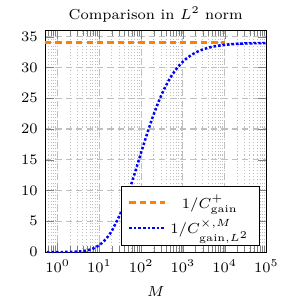}
        \includegraphics[scale=0.9]{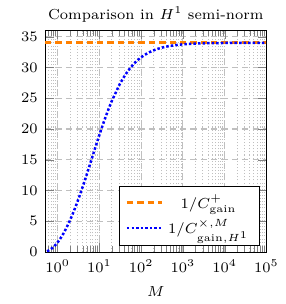}
    \end{minipage} \hfill
    \caption{Considering the \textit{1D Poisson problem} with $\bm{\mu}^{(1)}$, $k=1$ and the PINN prior $u_\theta$. Left -- Convergence of \cref{thm:comparison_add_mul} with the $L^2$ error. Right -- Convergence of \cref{thm:comparison_add_mul} with the \ROne{semi-}$H^1$ error.}\label{fig:Lap1D_gain_constants}
\end{figure}

We consider the first parameter $\bm{\mu}^{(1)}$ and the PINN prior $u_\theta$. We now evaluate the gain constants $C_\text{gain}^+$, $C_{\text{gain},H^1}^{\times,M}$and $C_{\text{gain},L^2}^{\times,M}$ (for different values of $M$), which are respectively defined in~\ROne{\cref{lem:error_estimation_add,lem:error_estimate_multiplicative}} for the additive and multiplicative approaches. The idea is to check the convergence of the two multiplicative gain constants towards the additive one, as proven in \cref{thm:comparison_add_mul}. The results are presented in \cref{fig:Lap1D_gain_constants} \RTwo{in $L^2$ norm and $H^1$ semi-norm, considering $M>\min|u_\theta|$.}

\cref{fig:Lap1D_gain_constants} shows that the multiplicative gain constant converges to the additive gain constant when $M$ increases, as expected in the theoretical results of \cref{thm:comparison_add_mul}.
\RTwo{Furthermore, we can compare the expected theoretical gains with the numerical gains obtained. For the additive approach, it would appear that the expected theoretical gain (approximately $34.01$) corresponds well with the numerical results of \cref{tab:case1_1D_both,tab:case1_1D_bothH1} (left subtable). As for the multiplicative approach, it appears that these constants align well when $M$ becomes large; however, for smaller $M$, the theoretical error estimates seem suboptimal. For example, for $M=3$, we find a theoretical gain of approximately $8.22$ in semi-norm $H^1$ compared to the $30$ obtained in practice in \cref{tab:case1_1D_bothH1} (left subtable). Moreover, in \cref{fig:Lap1D_gain_constants}, we see that the theoretical gain (in $L^2$ norm and $H^1$ semi-norm) converges to $0$ when $M$ tends to $\min|u_\theta|$, whereas in practice the numerical gains obtained have never been less than $1$. More precisely, for this test case in the $H^1$ semi-norm, we converge to $2.45$ numerically when $M$ is small.}

\subsubsection{Derivatives --- with both priors}\label{sec:Lap1D_derivatives}

\begin{figure}[ht!]
    \centering
    \includegraphics[scale=1]{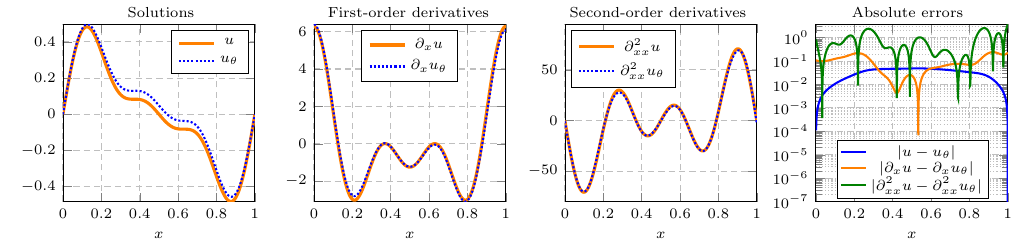}
    \caption{Considering the \textit{1D Poisson problem} with $\bm{\mu}^{(1)}$ and the PINN prior $u_\theta$, comparison between analytical solution and network prediction.
        From left to right: solution; first derivative; second derivative; errors.}\label{sinus1D_Pinns}
\end{figure}

\begin{figure}[ht!]
    \centering
    \includegraphics[scale=1]{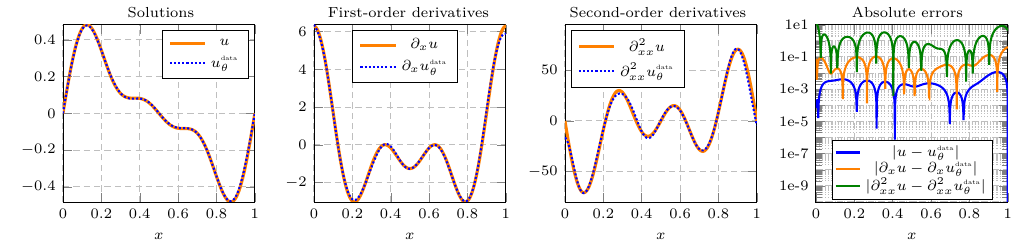}
    \caption{Considering the \textit{1D Poisson problem} with $\bm{\mu}^{(1)}$ and the data prior $u_\theta^\text{data}$,
        comparison between analytical solution and network prediction.
        From left to right: solution; first derivative; second derivative; errors.
    }\label{sinus1D_nn}
\end{figure}

To better explain the results of \cref{sec:Lap1D_error_estimations}, we compare the solution, the first- and second-order derivatives between the exact solution and the prediction of both priors, for selected parameter $\bm{\mu}^{(1)}$.
\cref{sinus1D_Pinns,sinus1D_nn} respectively present this comparison for the PINN prior $u_\theta$ and the data prior $u_\theta^\text{data}$.
We also compare the \ROne{relative} errors and gains obtained with these two priors for $N \in \{16,32\}$ in \ROne{\cref{tab:case1_1D_both,tab:case1_1D_bothH1}, respectively in $L^2$ norm and $H^1$ semi-norm.}
More precisely \ROne{for the $L^2$ norm}, we evaluate the additive error \smash{$e_{h,+}^{(1)}$} and the additive gain on FEM \smash{$G_+^{(1)}$}, respectively defined in~\eqref{eq:error_rel_add} and~\eqref{eq:gain_j}, for both PINN and data priors.
We also evaluate the multiplicative error \smash{$e_{h,M}^{(1)}$} and the multiplicative gain on FEM \smash{$G_M^{(1)}$} defined in~\eqref{eq:error_rel_add} and~\eqref{eq:gain_j_mul}, for both priors, with $M=3$ and $M=100$.

\begin{table}[ht!]
    \centering
    \GainsFixedMuTwoPriors{fig_testcase1D_test1_plots_FEM_param1.csv}{fig_testcase1D_test1_plots_compare_gains_param1.csv}
    \caption{Considering the \textit{1D Poisson problem} with $\bm{\mu}^{(1)}$, $k=1$ and $N \in \{16,32\}$. Left -- $L^2$ relative error obtained with FEM. Right -- Considering the PINN prior $u_\theta$ and the data prior $u_\theta^\text{data}$, $L^2$ relative errors and gains with respect to FEM, obtained with our methods. Our methods : additive approach, multiplicative approach with $M=3$ and $M=100$.}\label{tab:case1_1D_both}
\end{table}

\begin{table}[ht!]
    \centering
    \GainsFixedMuTwoPriors{fig_testcase1D_test1_plots_FEM_param1_H1.csv}{fig_testcase1D_test1_plots_compare_gains_param1_H1.csv}
    \caption{\ROne{Considering the \textit{1D Poisson problem} with $\bm{\mu}^{(1)}$, $k=1$ and $N \in \{16,32\}$. Left -- Semi-$H^1$ relative error obtained with FEM. Right -- Considering the PINN prior $u_\theta$ and the data prior $u_\theta^\text{data}$, Semi-$H^1$ relative errors and gains with respect to FEM, obtained with our methods. Our methods : additive approach, multiplicative approach with $M=3$ and $M=100$.}}\label{tab:case1_1D_bothH1}
\end{table}

The results reported in \cref{sinus1D_Pinns,sinus1D_nn} and \ROne{\cref{tab:case1_1D_both,tab:case1_1D_bothH1}} show that, even if the approach chosen to build the prior (physics-informed or data-driven training) gives a good approximation of the solution, the important point lies in the derivatives and mainly in the second-order derivatives, which are clearly better learned by PINNs.
Indeed, while the enriched FEM solution is more accurate using
the PINN prior (\ROne{\cref{tab:case1_1D_both,tab:case1_1D_bothH1}}), we see from
\cref{sinus1D_Pinns,sinus1D_nn} that the raw PINN
approximates the solution $u$ less accurately than the raw data-driven
the solution, but the PINN better approximates the derivatives.
As the error of the enriched FEM is mainly due to the $(k+1)$\textsuperscript{th} derivatives of the network being close to the $(k+1)$\textsuperscript{th} derivatives of the solution, this explains why the enriched FEM with data prior does not perform as well as with a PINN prior.
Therefore, PINNs have two advantages: they do not require training data and give better results.
\RTwo{Their main shortcoming is that training takes longer;
data could be used in addition to physics-informed training to speed up the process,
    but we do not expect it to improve the accuracy of the prior \RTwoN{derivatives}.}
However, we mention that if data are available for first- and second-order
derivatives, it could also be used to improve a purely data-driven prior.

\subsubsection{Gains achieved with the additive and multiplicative approaches -- with both priors}\label{sec:Lap1D_gains}

Considering a \RTwo{set $\mathcal{S}$ of $n_p=100$ parameter instances}, we now evaluate the gains $G_{+,\theta}$ and $G_+$ defined in~\eqref{eq:gain_add_num} with the PINN prior $u_\theta$ and with the data prior $u_\theta^\text{data}$.
We also compute $G_{M,\theta}$ and $G_M$, defined in~\eqref{eq:gain_mul_num}, similarly for both priors.
For fixed polynomial order $k=1$ and $N \in \{20,40\}$,
the results with the physics-informed prior $u_\theta$ and the data-driven prior $u_\theta^\text{data}$
are respectively presented in \cref{tab:case1_1D_PINNs,tab:case1_1D_data}.

\begin{table}[ht!]
    \centering
    \GainsTableOnedeg{fig_testcase1D_test1_gains_Tab_stats_case1_v1_degree1.csv}
    \caption{Considering the \textit{1D Poisson problem} and the PINN prior $u_\theta$. Left -- Gains in $L^2$ error of our methods with respect to PINN by taking $k=1$. Right -- Gains in $L^2$ error of our methods with respect to FEM by taking $k=1$. Our methods : additive approach, multiplicative approach with $M=3$ and $M=100$.}\label{tab:case1_1D_PINNs}
\end{table}

\begin{table}[ht!]
    \centering
    \GainsTableOnedegData{fig_testcase1D_test1_gains_Tab_stats_case1_v2_degree1.csv}
    \caption{Considering the \textit{1D Poisson problem} and the data prior $u_\theta^\text{data}$. Left -- Gains in $L^2$ error of our methods with respect to data Network by taking $k=1$. Right -- Gains in $L^2$ error of our methods with respect to FEM by taking $k=1$. Our methods : additive approach, multiplicative approach with $M=3$ and $M=100$.}\label{tab:case1_1D_data}
\end{table}

The previous results indicate that the average gain provided by the enriched FE with the PINN prior is significant, particularly when using the additive approach. These findings also confirm the behavior of the multiplicative prior method for varying values of $M$. In contrast, when applied with data-driven instead of physics-informed training, the same method does not yield similarly favourable results. Consequently, in the experiments we perform below, we only employ the PINN prior.

\RTwoN{Furthermore, it is important to note that we cannot expect enriched methods to have generalization or extrapolation properties for parameters $\mu \notin \mathcal{M}$. Indeed, these potential extrapolation properties would depend on those of the PINN. Since the PINN is a very smooth function, we might expect it to give reasonable results near the boundaries of $\mathcal{M}$. Intuitively, we might also expect that the further the parameters are from the training domain, the worse the results will be. However, there is no way to control this potential extrapolation property. For this test case, considering the additive approach with the PINN prior, $N=20$ and a sample of $100$ parameters in $[1,1.1]^3$ (close to $\mathcal{M}$), we find a mean gain of $9.18$, compared to the $140.74$ in \cref{tab:case1_1D_PINNs} for parameters in $\mathcal{M}$: the gain has been lowered by a factor of $15$. If we move away from the training area with parameters in $[1.1,1.5]^3$, we once again halve this factor with a mean gain of $4.48$, confirming the previous intuition.}

\subsubsection{\RTwo{Influence of the PINN initialization}}\label{sec:Lap1D_hyperparam}

\RTwo{In this section, we focus on the influence of network weight initialization on the quality of our results. To do this, we consider a set of 30 PINN initializations with the hyperparameters defined in \cref{tab:paramtest1_1D}. Thus, for each network, we are interested in the same way as in \cref{sec:Lap1D_gains}, in the gain of our enriched approaches on a set of $n_p=50$ parameter instances. For each prior, we evaluate the minimum, maximum, mean, and standard deviation on this sample. We then regroup these $30$ values in \cref{tab:case1_1D_30training}, in the form \textit{mean}$\pm$\textit{std}, where \textit{mean} and \textit{std} represent the mean and standard deviation, respectively, on the $30$ networks considered.}

\begin{table}[ht!]
    \centering
    \begin{filecontents*}{fig_testcase1D_test1_influence_weights_fem_summary.csv}
method,N,min_FEM,max_FEM,mean_FEM,std_FEM
Add,20,$23.04 \pm 5.66$,$306.85 \pm 50.27$,$148.56 \pm 24.05$,$61.57 \pm 10.46$
Add,40,$21.65 \pm 5.19$,$278.92 \pm 44.55$,$141.72 \pm 21.92$,$56.69 \pm 8.62$
Mult (M=3),20,$6.10 \pm 2.12$,$166.78 \pm 52.32$,$44.54 \pm 10.76$,$31.25 \pm 8.81$
Mult (M=3),40,$4.76 \pm 1.69$,$154.64 \pm 49.58$,$41.51 \pm 10.18$,$29.74 \pm 8.41$
Mult (M=100),20,$23.03 \pm 5.65$,$305.61 \pm 50.89$,$148.07 \pm 23.85$,$61.28 \pm 10.31$
Mult (M=100),40,$21.64 \pm 5.17$,$277.30 \pm 44.34$,$141.29 \pm 21.77$,$56.45 \pm 8.51$
\end{filecontents*}
\pgfplotstabletypeset[
    col sep=comma,
    string type,
    every head row/.style={
        before row={\toprule[1pt]
        & & \multicolumn{4}{c}{\textbf{Gains in $L^2$ rel error of our method w.r.t. FEM}} \\
        \cmidrule(lr){3-6}
        },
        after row=\cmidrule(lr){1-1} \cmidrule(lr){2-2} \cmidrule(lr){3-6}
    },
    every nth row={2}{before row=\cmidrule(lr){1-1} \cmidrule(lr){2-2} \cmidrule(lr){3-6}},
    every last row/.style={after row=\bottomrule[1pt]},
    columns/method/.style={column name=\textbf{method}},
    columns/N/.style={column name=\textbf{N}},
    columns/min_FEM/.style={column name=\textbf{min}},
    columns/max_FEM/.style={column name=\textbf{max}},
    columns/mean_FEM/.style={column name=\textbf{mean}},
    columns/std_FEM/.style={column name=\textbf{std}},
    columns={method,N,min_FEM,max_FEM,mean_FEM,std_FEM},
    ]{fig_testcase1D_test1_influence_weights_fem_summary.csv}
    
    \caption{\RTwo{Consider the \textit{1D Poisson problem} and a set of 30 PINN priors. $L^2$ error gains for $k=1$ of our methods compared to FEM : \textit{mean}$\pm$\textit{std} on the 30 networks considered. Our methods : additive approach, multiplicative approach with $M=3$ and $M=100$.}}\label{tab:case1_1D_30training}
\end{table}

\RTwo{In \cref{tab:case1_1D_30training}, we can see that the initialization of the PINN weights has an influence on the gains obtained with the different enriched approaches. More specifically, for the additive approach ($N=20$), for example, we can see that the minimum gain obtained on the sample averages $23$ across the 30 networks considered, with a standard deviation of $6$. It would therefore seem that there is always a gain in using the enriched approach, even if initialization influences its magnitude. On the other hand, even though the maximum is on average $307$, the standard deviation is close to $50$, which is not negligible. However, the results on the mean column seem close to the results obtained in \cref{tab:case1_1D_PINNs}. The results obtained with the multiplicative approach also appear to be consistent with those obtained previously.}

\subsection{1D general elliptic system and convection-dominated regime}\label{sec:Ell1D}

In this experiment, we consider the problem~\eqref{eq:ob_pde}
in a more complex form, still in a 1D ($d=1$) configuration:
\begin{equation*}
	\left\{
	\begin{aligned}
		\partial_x u-\frac{1}{\text{Pe}}\partial_{xx} u &= r, \; &  & \text{in } \; \Omega \times \mathcal{M}, \\
		u         & = 0, \;  &  & \text{on } \; \partial\Omega \times \mathcal{M},
	\end{aligned}
	\right.
\end{equation*}
with $\Omega=\RBothN{(0,1)}$ and $\partial\Omega$ its boundary, $r$ the reaction constant term, and $\text{Pe}$ the Péclet number, describing the ratio between the convection and the diffusion terms.
For all $x\in\Omega$, the analytical solution reads
\begin{equation}\label{eq:Ell1D_analytical}
	u(x,\bm{\mu})=r\left(x-\frac{e^{\text{Pe}\, x}-1}{e^{\text{Pe}}-1}\right) \,,
\end{equation}
with $p=2$ parameters $\bm{\mu}=(r,\text{Pe})\in\mathcal{M}=[1,2]\times[10,100]$.

\begin{rmrk}\label[rmrk]{rem:oscillations}
	In the large Péclet regime, i.e., for convection-dominated flows,
    the classical finite element method may generate oscillations
    when no specific treatment is applied,
    see e.g.~\cite{JohKnoNov2018}.
\end{rmrk}

In this test case, we construct only one prior, denoted $u_\theta$,
built from a PINN as presented in \cref{sec:prior_construction}.
We also compare the additive and multiplicative approaches by considering
polynomial order $k=1$, and since the solution is positive in $\Omega$,
we consider $M=0$ {for the multiplicative approach}.
We start by evaluating the error in \cref{sec:Ell1D_error_estimations},
then we compare the derivatives of the PINN prior and compare the different
approaches in \cref{sec:Ell1D_comparison}. Finally, we evaluate the gains obtained
in \cref{sec:Ell1D_gains} on a sample of parameters.
As we are dealing with a specific case, we will compare two methods
for imposing boundary conditions, as presented in \cref{sec:multiplicative_BC}:
the strong and the weak approaches.

\begin{rmrk}\label[rmrk]{rmk:Ell1D_N_nodes}
	As in \cref{sec:Lap1D}, the characteristic mesh size is $h=\frac{1}{N-1}$, where $N$ is the number of nodes considered.
\end{rmrk}

We consider a parametric PINN, depending on the problem parameters $\bm{\mu}$,
where we exactly impose the Dirichlet boundary conditions as presented in \cref{sec:exact_imposition_of_BC}. We define the prior $u_{\theta}$ and the level-set $\varphi$ as in~\cref{sec:Lap1D}.

\begin{table}[htbp]
    \centering
    \begin{tabular}{cc}
        \toprule
        \multicolumn{2}{c}{\textbf{Network - MLP}} \\
        \midrule
        \textit{layers} & $40,40,40,40,40$ \\
        \cmidrule(lr){1-2}
        $\sigma$ & tanh \\
        \bottomrule
    \end{tabular}
    \hspace{1cm}
    \begin{tabular}{cccc}
        \toprule
        \multicolumn{4}{c}{\textbf{Training}} \\
        \midrule
        \textit{lr} & 1e-3 & $n_{epochs}$ & \num{20000} \\
        \cmidrule(lr){1-2} \cmidrule(lr){3-4}
        \textit{decay} & 0.99 \\
        \cmidrule(lr){1-2}
        $N_\text{col}$ & \num{5000} \\
        \bottomrule
    \end{tabular}
    \hspace{1cm}
    \begin{tabular}{cccc}
        \toprule
        \multicolumn{4}{c}{\textbf{Loss weights}} \\
        \midrule
        $\omega_r$ & 1 & $\omega_\text{data}$ & 0 \\
        \cmidrule(lr){1-2} \cmidrule(lr){3-4}
        $\omega_b$ & 0 & $\omega_\text{sob}$ & 0 \\        
        \bottomrule
    \end{tabular}
    \caption{Network, training parameters (\cref{rmk:PINN_notations}) and loss weights for $u_\theta$ in the \textit{1D Elliptic case}.}\label{tab:paramtest2_1D}
\end{table}

\RBoth{Therefore, we will only consider the residual loss $J_r$ approached by a Monte-Carlo method as defined in~\eqref{eq:residual_loss_parametric_MC} with $N_\text{col}=\num{5000}$ collocation points (uniformly chosen on $\Omega\times\mathcal{M}$) and we seek to solve the minimisation problem \eqref{eq:minimization_problem}.}
The hyperparameters are given in \cref{tab:paramtest2_1D};
we use the Adam optimizer~\cite{KinBa2015}.

\subsubsection{Error estimates}\label{sec:Ell1D_error_estimations}

We start by testing the error estimates (\cref{lem:error_estimation_add,lem:error_estimate_multiplicative}) for the following two sets of parameters:
\begin{equation*}
	\bm{\mu}^{(1)}=(1.2,40) \quad \text{and} \quad \bm{\mu}^{(2)}=(1.5,90) \,,
\end{equation*}
by considering the PINN prior $u_\theta$.
For $j \in \{1, 2\}$, the aim is to compare for different mesh sizes $h$,
the $L^2$ relative errors \smash{$e_h^{(j)}$} obtained with the standard FEM, defined in~\eqref{eq:error_rel_FEM}, \smash{$e_{h,+}^{(j)}$} obtained with the additive approach and \smash{$e_{h,M}^{(j)}$} obtained with the multiplicative approach (taking $M=0$), defined in~\eqref{eq:error_rel_add}.  We will consider the two implementations of the boundary conditions for the multiplicative approach: the strong and the weak BC, as presented in \cref{sec:boundary_conditions}.
The results are presented in \cref{fig:case2_1D} by varying the mesh size $h$, \RBoth{considering $N\in\{16,32,64,128,256\}$ as presented in \cref{rmk:Ell1D_N_nodes}.}

\begin{figure}[ht!]
	\centering
	\begin{subfigure}{0.48\linewidth}
		\centering
		\cvgFEMCorrMultSWOnedeg{fig_testcase1D_test2_cvg_FEM_case2_v1_param1_degree1.csv}{fig_testcase1D_test2_cvg_FEM_case2_v1_param1_degree2.csv}{fig_testcase1D_test2_cvg_Corr_case2_v1_param1_degree1.csv}{fig_testcase1D_test2_cvg_Mult_case2_v1_param1_degree1_M0.0.csv}{fig_testcase1D_test2_cvg_Mult_case2_v1_param1_degree1_M0.0_weak.csv}{5e-6}
		\caption{Case of $\bm{\mu}^{(1)}$}
	\end{subfigure}
	\begin{subfigure}{0.48\linewidth}
		\centering
		\cvgFEMCorrMultSWOnedeg{fig_testcase1D_test2_cvg_FEM_case2_v1_param2_degree1.csv}{fig_testcase1D_test2_cvg_FEM_case2_v1_param2_degree2.csv}{fig_testcase1D_test2_cvg_Corr_case2_v1_param2_degree1.csv}{fig_testcase1D_test2_cvg_Mult_case2_v1_param2_degree1_M0.0.csv}{fig_testcase1D_test2_cvg_Mult_case2_v1_param2_degree1_M0.0_weak.csv}{6e-5}
		\caption{Case of $\bm{\mu}^{(2)}$}
	\end{subfigure}
	\caption{Considering the \textit{1D Elliptic case} and the PINN prior $u_\theta$. Left -- Considering $\bm{\mu}^{(1)}$. $L^2$ error on $h$ obtained with standard FEM \smash{$e_h^{(1)}$} (solid lines) with $k=1$ and $k=2$, the additive approach \smash{$e_{h,+}^{(1)}$} (dashed lines) with $k=1$ and the multiplicative approach \smash{$e_{h,M}^{(1)}$} (dotted lines) with $k=1$, considering strong and weak BC. Right -- Same for $\bm{\mu}^{(2)}$.}\label{fig:case2_1D}
\end{figure}

In \cref{fig:case2_1D}, we see that the enriched approaches seem to give better results than standard FEM except for the multiplicative approach with strong imposition of boundary conditions. Moreover, this approach, which imposes $p_h^\times=1$ on $\partial\Omega$, does not follow the expected convergence order \RTwo{(i.e., for $k=1$, the expected convergence order in the $L^2$ norm is $2$, compared to $1$ obtained numerically)}. 
The additive approach seems much less effective here than in the previous
experiment of \cref{sec:Lap1D}, whereas the multiplicative approach with weak BC seems to significantly improve the results obtained with standard FEM. 
\RTwo{Nevertheless, these two approaches appear to respect the expected convergence order in the $L^2$ norm, namely $2$ for polynomials of order $k=1$.}
A comparative study of the different methods is given in \cref{sec:Ell1D_comparison}. In addition, we see that the standard FEM with polynomial order $k=2$ is
clearly less accurate than the multiplicative approach using weak BC applied
with polynomial order $k=1$.

\subsubsection{Comparison of different approaches}\label{sec:Ell1D_comparison}

\begin{figure}[ht!]
	\centering
	\includegraphics[scale=1]{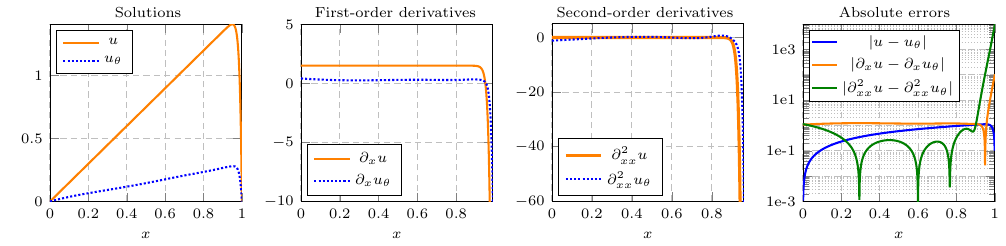}
	\caption{Considering the \textit{1D Ellipctic case} with $\bm{\mu}^{(2)}$ and the PINN prior $u_\theta$, comparison between analytical solution and network prediction.
    From left to right: solution; first derivative; second derivative; errors.}\label{fig:case2_1D_der}
\end{figure}

\begin{figure}[ht!]
	\centering
    \begin{subfigure}{0.48\linewidth} \centering
		\hspace{-12pt}\includegraphics[scale=0.84]{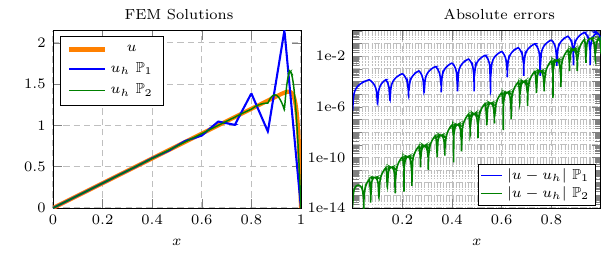}
		\caption{FEM solutions with polynomials order approximation $k=1$ and $k=2$, and absolute errors.}\label{fig:case2_1D_plots_fem}
    \end{subfigure} \hfill
    \begin{subfigure}{0.48\linewidth} \centering
		\hspace{-12pt}\includegraphics[scale=0.84]{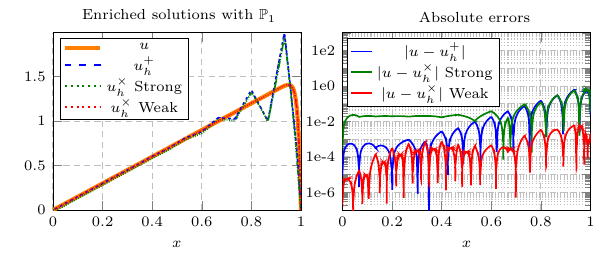}
		\caption{Enriched solutions with polynomial order approximation $k=1$, and absolute errors.}\label{fig:case2_1D_plots_add}
    \end{subfigure}
	\caption{Considering the \textit{1D Elliptic case} with $\bm{\mu}^{(2)}$, $N=16$ and the PINN prior $u_\theta$. Comparison of the solution obtained with the different methods with the analytical solution. For each enriched method, comparison of the correction term with the analytical one. Different methods : standard FEM, additive approach, multiplicative approach by taking $M=0$ (strong and weak BC).}\label{fig:case2_1D_plots}
\end{figure}

We now focus on the second parameter $\bm{\mu}^{(2)}$. We first look at the PINN prediction for this parameter and its derivatives in \cref{fig:case2_1D_der}.
As in the previous section, we consider the following approaches: standard FEM, the additive approach and the multiplicative approach (with $M=0$) with strong or weak imposition of boundary conditions.
We compare the different methods in \cref{tab:case2_1D_comparison}, where we can see the different errors obtained with the considered methods for $k=1$ and $N\in\{16,32\}$ as well as the gains obtained in comparison with standard FEM. Next, we take a closer look at the solutions obtained with the different approaches in \cref{fig:case2_1D_plots};
for each method, we compare the solution obtained ($u_h$ for standard FEM, $u_h^+$ for the additive approach and $u_h^\times$ for the multiplicative approach, with strong or weak BC imposition) with the analytical solution $u$. For the enriched methods, using the PINN prior, we will also compare the proposed correction; namely, for the additive approach, we will compare $p_h^+$ with $u-u_\theta$ and for the multiplicative one $p_h^\times$ with $u/u_\theta$ (with $u_\theta>0$ in $\Omega$).

In \cref{fig:case2_1D_der}, we can see that PINN has difficulties to capture the solution and that the prediction it provides is far from the analytical solution. As for its derivatives, they seem to be relatively inaccurate compared to the analytical.
Indeed, since the PINN is a smooth function, it has trouble approximating functions with very sharp gradients such as the one of~\eqref{eq:Ell1D_analytical}.

\begin{table}[ht!]
	\centering
	\GainsFixedMu{2}{fig_testcase1D_test2_plots_FEM.csv}{fig_testcase1D_test2_plots_compare_gains.csv}
	\caption{Considering the \textit{1D Elliptic case} with $\bm{\mu}^{(2)}$, $k=1$ and $N\in\{16,32\}$. Left -- $L^2$ relative error obtained with FEM. Right -- Considering the PINN prior $u_\theta$, $L^2$ relative errors and gains with respect to FEM, obtained with our methods. Our methods : additive approach, multiplicative approach by taking $M=0$ (strong and weak BC).}\label{tab:case2_1D_comparison}
\end{table}

In \cref{fig:case2_1D_plots,tab:case2_1D_comparison}, we present a comparison of the different approaches proposed. In \cref{fig:case2_1D_plots_fem}, we first notice the oscillations anticipated in \cref{rem:oscillations} for standard FEM at both polynomial
orders $k=1$ and $k=2$.
This behavior is also seen in the additive enrichment (blue dashed line in \cref{fig:case2_1D_plots_add}), which does not seem to give better results than standard FEM due to the derivatives presented in \cref{fig:case2_1D_der}, and for the multiplicative approach with strong boundary conditions (green dotted line in \cref{fig:case2_1D_plots_add}).
However, weakly imposing the BC gives the appropriate results (red dotted line in \cref{fig:case2_1D_plots_add}). \ROne{It seems that for this specific test case, the prediction obtained (close to a multiplicative constant of the solution) favors the multiplicative approach over the additive one.} 

\RTwoN{More specifically, since the accuracy of enriched approaches is completely dependent on the prior, it is not clear in advance which approach (additive or multiplicative) will yield the best results. However, this strong dependence means that we can intuitively choose the enriched approach based on our knowledge of PINNs. For example, as explained above, since PINNs are regular functions, they have difficulty approximating functions with very sharp gradients, as in this test case. Thus, knowing the nature of the test case in question, we can get an idea of the quality of the prior and therefore estimate which approach should be more compatible with it. To go further, if we want a more precise idea of the approach to focus on, we could, after training, calculate the theoretical gain constants using a coarse FEM as a reference solution, for example, which in a parametric approach could be seen as an offline cost. However, there is no theoretical guarantee that the constants will exhibit the same behavior for different parameters, particularly depending on the type of parameter considered.}

\subsubsection{Gains achieved with the additive and the multiplicative approaches}\label{sec:Ell1D_gains}

Considering a \RTwo{set $\mathcal{S}$ of $n_p=50$ parameter instances}, we will evaluate the gains $G_{+,\theta}$ and $G_+$ defined in~\eqref{eq:gain_add_num} with the PINN prior $u_\theta$. We will evaluate $G_{M,\theta}$ and $G_M$, defined in~\eqref{eq:gain_mul_num}, in the same way with $M=0$ for the multiplicative approach, by considering the two implementations of the boundary conditions; strong and weak BC. The results are presented in \cref{tab:case2_1D_PINNs} for fixed $k=1$ and $N \in \{20,40\}$ fixed.

\begin{table}[ht!]
	\centering
	\GainsTableOnedeg{fig_testcase1D_test2_gains_Tab_stats_case2_v1_degree1.csv}
	\caption{Considering the \textit{1D Ellipctic case}, $k=1$ and the PINN prior $u_\theta$. Left -- Gains in $L^2$ relative error of our methods with respect to PINN. Right -- Gains in $L^2$ relative error of our methods with respect to FEM. Our methods : additive approach, multiplicative approach with $M=0$ (strong and weak BC).}\label{tab:case2_1D_PINNs}
\end{table}

\cref{tab:case2_1D_PINNs} confirms the above results. The multiplicative approach with weak BC seems to give the best results on our \RTwo{set $\mathcal{S}$ of parameter instances}. The additive and multiplicative approaches with strong BC imposition do not appear to be very effective on this test case. In particular, even though the additive approach improves the standard FEM error by a factor of 3, we have seen in \cref{sec:Ell1D_comparison} that the solutions obtained do not correspond to the expected solution, whereas the multiplicative approach with low BCs does. In the following, we will only consider the additive approach, as it seems to be the most efficient one, except in special cases such as the one under consideration in this \cref{sec:Ell1D}. Indeed, the following test cases will not contain boundary layers or strong gradients.

\subsection{\RBothN{1D non-smooth transmission problem}}\label{sec:SingLap1D}

\RBothN{Let $\Omega=(0,\pi)$ and $\partial\Omega$ its boundary. In this test case, the physical properties of the medium are discontinuous with respect to space, and we define the interface point $I = \frac{\pi}{2}$ at the middle of the domain. We consider the folllowing one-dimensional non-smooth transmission problem:}
\begin{equation}\label{eq:SingLap1D}
    \RBothN{\left\{
    \begin{aligned}
        -\sigma \partial_{xx} u & = f, \; &  & \text{in } \;  \Omega\setminus\lbrace{I\rbrace} \times \mathcal{M}, \vphantom{\lim_{x\to I^+}}        \\
        u                & = 0, \; &  & \text{on } \; \partial\Omega \times \mathcal{M}, \vphantom{\lim_{x\to I^+}} \\
        \lim_{x\to I^+} \sigma(x,\bm{\mu}) \partial_x u(x,\bm{\mu}) & = \lim_{x\to I^-} \sigma(x,\bm{\mu}) \partial_x u(x,\bm{\mu}), \; &  & \bm{\mu} \in \mathcal{M}, \\
        \lim_{x\to I^+} u(x,\bm{\mu}) & = \lim_{x\to I^-} u(x,\bm{\mu}), \; &  & \bm{\mu} \in \mathcal{M},
    \end{aligned}
    \right.}
\end{equation}
\RBothN{where the diffusion coefficient $\sigma$ is defined by}
\begin{equation*}
    \RBothN{\sigma(x,\bm{\mu})=\left\{
    \begin{aligned}
        \sigma_1, & \; \text{if } x\in [0,I], \\
        \sigma_2, & \; \text{if } x\in (I,\pi],
    \end{aligned}
    \right.}
\end{equation*}
\RBothN{with $\bm{\mu}=(\sigma_1,\sigma_2)\in\mathcal{M}=[2.5,3.5]\times[0.5,1.5]$ (we have $p=2$ parameters).}
\RBothN{We prescribe a family of exact solutions (depending on the parameter vector) defined by}
\begin{equation*}
    \RBothN{u(x,\bm{\mu})=\left\{
    \begin{aligned}
        \frac{1}{\sigma_1}\sin(2x), & \text{ if } x\in (0,I), \\
        \frac{1}{\sigma_2}\sin(2x), & \text{ if } x\in (I,\pi),
    \end{aligned}
    \right.}
\end{equation*}
\RBothN{and we deduce the associated right-hand side $f(x,\bm{\mu})=4\sin(2x)$.}

\RBothN{At first, we will resonate in a non-parametric framework, focusing on the parameter $\bm{\mu}^{(1)}=(3,1)$. Therefore, in \cref{sec:SingLap1D_error_estimations}, we first present the error estimates obtained on this parameter, considering three different priors (whose construction is detailed in \cref{sec:SingLap1D_priors}). The first one (denoted by $u_\theta$) is a PINN with a generic MLP architecture. The second one (denoted by \smash{$u_\theta^\text{data}$}) is constructed from data on the solution and its derivatives (obtained from the analytical solution). For the third one (denoted by \smash{$u_\theta^\text{sing}$}), we use a PINN with the architecture proposed in \cite{reconns2024}, specifically designed to better capture the behavior of the solution around the interface by enriching the MLP with a singular function. The aim is to show that for this type of problem, using an architecture adapted to the problem leads to better results than a generic one. More specifically, we will show that a naive network is not suited to the enriched approach we propose. However, by enriching the network in a similar way that classical finite element methods are enriched to capture known discontinuities, we obtain results similar to the other test cases considered, for low values of the polynomial degree $k$. To do this, we will focus in \cref{sec:SingLap1D_derivatives} on the predictions and derivatives of the different non-parametric priors and analyze the results of the enriched approach with the three priors. Finally, in \cref{sec:SingLap1D_gains}, we will return to our parametric framework  by constructing the enriched prior \smash{$u_\theta^\text{sing}$} in a parametric way (considering $\bm{\mu}\in\mathcal{M}$) and evaluate the gains obtained with it on a sample of parameters.}

\begin{rmrk}\label[rmrk]{rmk:SingLap1D_N_nodes}
	\RBothN{As in the previous sections, the characteristic mesh size is $h=\frac{1}{N-1}$, where $N$ is the number of nodes considered. For the finite element methods, we will choose to consider that the interface $I$ lies on a mesh node (by taking odd values of $N$) to avoid integration issues.}
\end{rmrk}

\subsubsection{\RBothN{Construction of the different priors}}\label{sec:SingLap1D_priors}

\RBothN{Here, we will discuss the differences in training the three non-parametric networks (used in \cref{sec:SingLap1D_error_estimations,sec:SingLap1D_derivatives}) for the specific parameter $\bm{\mu}^{(1)}$. We will also add the construction of the parametric network \smash{$u_\theta^\text{sing}$} (considered in \cref{sec:SingLap1D_gains}). For PINNs, we will not strongly impose the boundary conditions (i.e. with a levelset function), but rather include them in the loss function.}

\boldparagraph{\RBothN{Non-parametric physics-informed training with classical architecture}}
\RBothN{The first non-parametric prior $u_\theta$ is a PINN with a generic MLP architecture. Since the problem is not defined on $I$, we can simply consider that the integral of the residual loss is divided between the two subdomains (to the left and right of $I$). By applying an independent Monte Carlo method, we obtain the same approximate residual loss $J_r$ defined in \eqref{eq:residual_loss_parametric_MC} (considering $N_\text{col}=\num{3000}$ collocation points uniformly chosen on $\Omega\setminus\lbrace{I\rbrace}$). Since the boundary conditions are not strongly imposed in this network, we also include the boundary loss $J_b$ defined in \eqref{eq:boundary_loss_parametric_MC}. It should be noted that, by definition of the network in question, the continuity condition is automatically validated, but the flow condition cannot be satisfied. In practice, the considered prior will aim to be accurate on the residual of transmission problem, while the solution itself will be very poorly represented. The hyperparameters used are given in \cref{tab:paramtest3_1D}.}

\begin{table}[htbp]
    \centering
    \begin{tabular}{cc}
        \toprule
        \multicolumn{2}{c}{\textbf{Network - MLP}} \\
        \midrule
        \textit{layers} & $20,20,20$ \\
        \cmidrule(lr){1-2}
        $\sigma$ & tanh \\
        \bottomrule
    \end{tabular}
    \hspace{0.3cm}
    \begin{tabular}{cccc}
        \toprule
        \multicolumn{4}{c}{\textbf{Training}} \\
        \midrule
        \textit{lr} & 1e-3 & $n_{epochs}$ & \num{2000} \\
        \cmidrule(lr){1-2} \cmidrule(lr){3-4}
        \textit{decay} & 0.99 \\
        \cmidrule(lr){1-2}
        $N_{\text{col}/\text{data}}$ & \num{3000} \\
        \bottomrule
    \end{tabular}
    \hspace{0.3cm}
    \begin{tabular}{cccccc}
        \toprule
        \multicolumn{6}{c}{\textbf{Loss weights}} \\
        \midrule
        \multicolumn{2}{c}{$u_\theta$} & \multicolumn{2}{c}{$u_\theta^\text{data}$} & \multicolumn{2}{c}{$u_\theta^\text{sing}$} \\
        \cmidrule(lr){1-2} \cmidrule(lr){3-4} \cmidrule(lr){5-6}
        $\omega_r$ & 1 & $\omega_\text{data}$ & 1 & $\omega_r$ & 1 \\
        \cmidrule(lr){1-2} \cmidrule(lr){3-4} \cmidrule(lr){5-6}
        $\omega_b$ & 20 & & & $\omega_b$ & 20 \\
        \cmidrule(lr){1-2} \cmidrule(lr){5-6}
        & & & & $\omega_\text{int}$ & 1 \\
        \bottomrule
    \end{tabular}
    \caption{\RBothN{Network, training parameters (\cref{rmk:PINN_notations}) and loss weights for non-parametric priors $u_\theta$, $u_\theta^\text{data}$ and $u_\theta^\text{sing}$ in the \textit{1D transmission problem}. Considering $N_\text{col}$ collocation points for PINNs and $N_\text{data}$ data for the data prior. Weights not specified in the table are defined to be zero.}}\label{tab:paramtest3_1D}
\end{table}

\boldparagraph{\RBothN{Non-parametric data-driven training with classical architecture}}
\RBothN{For the second non-parametric prior $u_\theta^\text{data}$, a  network is trained only on the data (constructed from the analytical solution), considering the same architecture as above. In \cref{sec:Lap1D}, we have already seen that a network based on data is not sufficient to correctly capture the derivatives of the solution. Hence, based on \cite{reconns2024}, we define the following modified data loss function, including information on the derivatives:}
\begin{equation*}
    \RBothN{J_\text{data}(\theta) =
    \frac 1 {N_\text{data}} \sum_{i=1}^{N_\text{data}} \big| u_\theta^\text{data}\big(\bm{x}_\text{data}^{(i)},\bm{\mu}^{(1)}\big) - u\big(\bm{x}_\text{data}^{(i)},\bm{\mu}^{(1)}\big) \big|^2 + \big| \partial_x u_\theta^\text{data}\big(\bm{x}_\text{data}^{(i)},\bm{\mu}^{(1)}\big) - \partial_x u\big(\bm{x}_\text{data}^{(i)},\bm{\mu}^{(1)}\big) \big|^2 \,,}
\end{equation*}
\RBothN{with $N_\text{data}=\num{3000}$ points uniformly chosen on $\Omega\setminus\lbrace{I\rbrace}$. The hyperparameters used are also given in \cref{tab:paramtest3_1D}.}

\boldparagraph{\RBothN{Non-parametric physics-informed training with enriched architecture}}

\RBothN{The third non-parametric prior $u_\theta^\text{sing}$ is a PINN with an enriched architecture specifically designed to capture the behavior of the solution around the interface. This architecture was proposed in \cite{reconns2024} under the name ReCoNNs (Regularity-Conforming Neural Networks). More precisely, we define the non-parametric network as follows:}
\begin{equation}\label{eq:SingLap1D_sing_architecture}
    \RBothN{u_\theta^\text{sing}(x,\bm{\mu}) = w_\theta^0(x,\bm{\mu}) + w_\theta^1(x,\bm{\mu})\frac{|x - I|}{2} \,,}
\end{equation}
\RBothN{with $w_\theta=(w_\theta^0,w_\theta^1)$ a classic MLP and $\bm{\mu}=\bm{\mu}^{(1)}$. This architecture allows us to capture the discontinuity of the first derivative of the solution at the interface $I$. The network is trained in a physics-informed way, minimizing the same residual loss function $J_r$ as for the first prior and the same boundary loss function $J_b$. The continuity condition is automatically satisfied by the construction of the network. For the flow condition, we then add the interface loss $J_\text{int}$ defined by:}
\begin{equation*}
    \RBothN{J_\text{int}(\theta) = \left|\lim_{x\to I^+} \sigma(x,\bm{\mu}^{(1)}) \partial_x u_\theta^\text{sing}(x,\bm{\mu}^{(1)}) - \lim_{x\to I^-} \sigma(x,\bm{\mu}^{(1)}) \partial_x u_\theta^\text{sing}(x,\bm{\mu}^{(1)})\right|^2,}
\end{equation*}
\RBothN{where the limits of network derivatives at the interface are calculated analytically, as detailled in \cite{reconns2024}. The minimization problem then becomes:}
\begin{equation*}
    \RBothN{\theta^\star = \argmin_{\theta}
    \big( \omega_r J_r(\theta) + \omega_b J_b(\theta) + \omega_\text{int} J_\text{int}(\theta) \big),}
\end{equation*}
\RBothN{with $\omega_r$, $\omega_b$ and $\omega_\text{int}$ some weights to balance the different terms of the loss function. The hyperparameters used are still given in \cref{tab:paramtest3_1D}.}

\boldparagraph{\RBothN{Parametric physics-informed training with enriched architecture}}

\RBothN{For \cref{sec:SingLap1D_gains}, we will consider the same enriched architecture as above, but in a parametric way. Thus, we define the parametric network as defined in \eqref{eq:SingLap1D_sing_architecture}, consdering $\bm{\mu}\in\mathcal{M}$ and the modified interface loss function:}
\begin{equation*}
    \RBothN{J_\text{int}(\theta) = \frac{1}{N_\text{col}} \sum_{i=1}^{N_\text{col}} \left|\lim_{x\to I^+} \sigma(x,\bm{\mu}_\text{col}^{(i)}) \partial_x u_\theta^\text{sing}(x,\bm{\mu}_\text{col}^{(i)}) - \lim_{x\to I^-} \sigma(x,\bm{\mu}_\text{col}^{(i)}) \partial_x u_\theta^\text{sing}(x,\bm{\mu}_\text{col}^{(i)})\right|^2.}
\end{equation*}
\RBothN{The hyperparameters used for the network and training are the same as for the non-parametric PINN prior with enriched architecture (defined in \cref{tab:paramtest3_1D}), except for the number of epochs set to $n_{epochs}=\num{10000}$ and the initial learning rate set to $\textit{lr}=1.8e$-$2$.}

\subsubsection{\RBothN{Error estimates --- with the three non-parametric priors}}\label{sec:SingLap1D_error_estimations}

\RBothN{We perform the same test as in the previous sections, considering only the additive approach with the three non-parametric priors presented in \cref{sec:SingLap1D_priors}.
The results in $L^2$ norm are presented for $\bm{\mu}^{(1)}$ in \cref{fig:case3_1D} for fixed $k\in \{1,2,3\}$ and by varying the number of nodes $N \in \{21,41,61,81,101\}$ (\cref{rmk:SingLap1D_N_nodes}).}

\RBothN{The results of the enriched approach are obtained by calculating the source term weakly, meaning only using the first derivatives of the prior. Also, like in the other test cases, the derivatives are calculated analytically using \texttt{PyTorch}. However, in the case of the particular architecture used in \smash{$u_\theta^\text{sing}$}, the derivatives produced are discontinuous, which must be taken into account during integration. More precisely, using an odd number of nodes $N$ (as stated in \cref{rmk:SingLap1D_N_nodes}), two cells of the mesh intersect at the interface $I$. Using a Galerkin discontinuous space to represent the first derivative of the prior, we will define the value at $I$ of the left (resp. right) cell by its value at $I-\epsilon$ (resp. $I+\epsilon$) with $\epsilon=10^{-6}$. This way, we ensure that the correct value of the derivative is used on each side of the interface during integration.}

\begin{figure}[ht!]
	\centering
	\hspace{-1cm}
	\begin{subfigure}{0.32\linewidth}
		\centering
		\resizebox{1.1\linewidth}{!}{
			\cvgFEMCorrThreePriors{1}{fig_testcase1D_test3_cvg_FEM_case3_v4_param1.csv}{fig_testcase1D_test3_cvg_Corr_case3_v3_param1.csv}{fig_testcase1D_test3_cvg_Corr_case3_v2_param1.csv}{fig_testcase1D_test3_cvg_Corr_case3_v4_param1.csv}{1e-6}{$u_\theta^\text{data}$}{$u_\theta^\text{sing}$}
		}
		\caption{$k=1$}
	\end{subfigure}
	\begin{subfigure}{0.32\linewidth}
		\centering
		\resizebox{1.1\linewidth}{!}{
			\cvgFEMCorrThreePriors{2}{fig_testcase1D_test3_cvg_FEM_case3_v4_param1.csv}{fig_testcase1D_test3_cvg_Corr_case3_v3_param1.csv}{fig_testcase1D_test3_cvg_Corr_case3_v2_param1.csv}{fig_testcase1D_test3_cvg_Corr_case3_v4_param1.csv}{5e-8}{$u_\theta^\text{data}$}{$u_\theta^\text{sing}$}
		}
		\caption{$k=2$}
	\end{subfigure}
	\begin{subfigure}{0.32\linewidth}
		\centering
		\resizebox{1.1\linewidth}{!}{
			\cvgFEMCorrThreePriors{3}{fig_testcase1D_test3_cvg_FEM_case3_v4_param1.csv}{fig_testcase1D_test3_cvg_Corr_case3_v3_param1.csv}{fig_testcase1D_test3_cvg_Corr_case3_v2_param1.csv}{fig_testcase1D_test3_cvg_Corr_case3_v4_param1.csv}{1e-9}{$u_\theta^\text{data}$}{$u_\theta^\text{sing}$}
		}
		\caption{$k=3$}
	\end{subfigure}
	\caption{\RBothN{Considering the \textit{1D transmission problem} with $\bm{\mu}^{(1)}$. Left -- $L^2$ relative error on $h$, obtained with the standard FEM \smash{$e_h^{(1)}$} (solid line) and the additive approach \smash{$e_{h,+}^{(1)}$} (dashed lines), with $k=1$, by considering the PINN prior with classic architecture $u_\theta$, the data-driven prior \smash{$u_\theta^\text{data}$} and the PINN with enriched architecture \smash{$u_\theta^\text{sing}$} (all in a non-parametric way). Middle -- Same with $k=2$. Right -- Same with $k=3$.}}\label{fig:case3_1D}
\end{figure}

\RBothN{The results presented in \cref{fig:case3_1D} show that the additive approach with the PINN prior and the enriched architecture \smash{$u_\theta^\text{sing}$} provides the best results. In particular, we observe that for polynomial order $k=1$, the errors obtained with PINNs are much lower than those obtained with the data network, and that the gain compared to standard FEM is significant. Furthermore, we also observe that the convergence rates are similar for all three priors and match those of the classical FEM. This confirms that using a network architecture adapted to the problem yields better results in our enriched finite element method. On the other hand, we can see that the data-driven network seems to deteriorate the results of the standard approach, particularly for $k\ne 1$, which is a behavior we did not observe elsewhere. As this test case does not enter into the theoretical framework of \cref{sec:additive_prior}, we cannot assert the causes of this; further study is warranted. However, from what we have seen so far, we have observed that the more we increase the polynomial order, the more the order of the derivatives that control the error increases. If the reasoning is similar here, a network trained solely by data could produce such a result. To better understand these results, we now analyze the predictions and derivatives of the three priors in \cref{sec:SingLap1D_derivatives}.}

\subsubsection{\RBothN{Derivatives --- with the three non-parametric priors}}\label{sec:SingLap1D_derivatives}

\begin{figure}[ht!]
    \centering
    \includegraphics[scale=1]{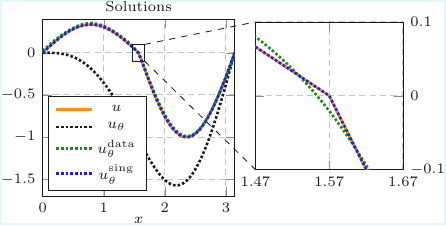}
    \hspace{-2pt}\includegraphics[scale=1]{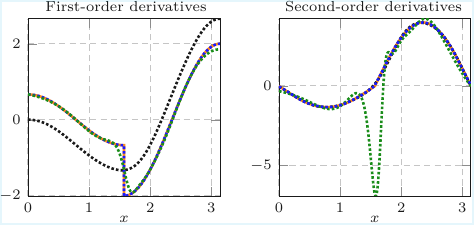}
    \caption{\RBothN{Considering the \textit{1D transmission problem} with $\bm{\mu}^{(1)}$ and the three different non-parametric priors, the PINN prior with classic architecture $u_\theta$ (black), the data-driven prior \smash{$u_\theta^\text{data}$} (green) and the PINN with enriched architecture \smash{$u_\theta^\text{sing}$} (blue). Comparison between analytical solution and network prediction. From left to right: solution; zoom on solution at interface; first derivative; second derivative.}}\label{plottest3}
\end{figure}

\RBothN{To better explain the results of \cref{sec:SingLap1D_error_estimations} for selected parameter $\bm{\mu}^{(1)}$, we compare the solution, the first- and second-order derivatives between the exact solution and the prediction of the three priors constructed in \cref{sec:SingLap1D_priors}. For analytical second derivatives, they are calculated independently on the left and right intervals. The results are presented in \cref{plottest3}.}

\RBothN{In \cref{plottest3}, we see that the prior that gives the best approximation of the solution is the PINN with the singular architecture, followed by the data-driven network, which seems to smooth the solution at the interface somewhat. The standard PINN, on the other hand, has great difficulty capturing the solution. In terms of first derivatives, we see that the enriched PINN is also very successful, as is the case for the solution. The data-driven network seems to give a smooth approximation of the first derivative, while the standard PINN still seems quite far off. However, the data-driven network has the most difficulty approximating second derivatives, followed by the standard PINN, while the singular network approximates them very well. In view of the results of \cref{sec:SingLap1D_error_estimations}, it would therefore seem that it is the higher-order derivatives that have the greatest impact on the results, even if the prior itself is not optimal. This conclusion also holds true in other test cases where our theoretical framework is applicable.}

\subsubsection{\RBothN{Gains achieved with the additive approach -- with the parametric enriched PINN prior}}\label{sec:SingLap1D_gains}

\RBothN{Considering a \RTwo{set $\mathcal{S}$ of $n_p=50$ parameter instances} and the parametric enriched PINN prior $u_\theta^\text{sing}$ defined in \cref{sec:SingLap1D_priors},
we now evaluate the gains $G_{+,\theta}$ and $G_+$ defined in~\eqref{eq:gain_add_num}.
The results are presented in \cref{tab:case3_1D}
for $k \in \{1,2,3\}$ and $N \in \{21, 41\}$.}

\begin{table}[ht!]
    \centering
    \GainsTableAlldeg{fig_testcase1D_test3_gains_Tab_stats_case3_v1.csv}
    \caption{\RBothN{Considering the \textit{1D transmission problem}, $k\in\{1,2,3\}$ and the parametric PINN prior with enriched architecture \smash{$u_\theta^\text{sing}$}. Left -- Gains in $L^2$ relative error of the additive method with respect to PINN. Right -- Gains in $L^2$ relative error of our approach with respect to FEM.}}\label{tab:case3_1D}
\end{table}

\RBothN{The results in \cref{tab:case3_1D} confirm the results obtained for $\mu^{(1)}$ in \cref{sec:SingLap1D_error_estimations}. It would appear that the additive approach with the enriched PINN prior $u_\theta^\text{sing}$ provides significant gains, particularly for $k=1$. For $k=3$, the minimum obtained is greater than $1$, so we see that, unlike when using the data-driven prior, the results of the standard method do not appear to be degraded. On the other hand, as the degree of the polynomial increases, the average gains appear to decrease, and quite considerably so.}

\subsection{2D Poisson problem in a square domain}\label{sec:Lap2D}

We now consider the problem of \cref{sec:Lap1D}
but in two dimensions ($d=2$), with,
\begin{equation}\label{eq:Lap2D}
	\left\{
	\begin{aligned}
		-\Delta u & = f, \; &  & \text{in } \; \Omega \times \mathcal{M}, \\
		u         & =0, \;  &  & \text{on } \; \partial\Omega \times \mathcal{M},
	\end{aligned}
	\right.
\end{equation}
with $\Delta$ the Laplace operator on the domain
${\Omega=\RBothN{(-0.5 \pi, 0.5 \pi)}}^2$ with boundary $\partial\Omega$,
and $\mathcal{M} \subset \mathbb{R}^p$ the parameter space (with $p$ the number of parameters).
We define the right-hand side $f$ such that the solution is given by
\begin{equation}\label{eq:analytical_solution_Lap2D}
	u(\bm{x},\bm{\mu})=\exp\left(-\frac{{(x-\mu_1)}^2+{(y-\mu_2)}^2}{2}\right)\sin(\kappa x)\sin(\kappa y),
\end{equation}
with $\bm{x}=(x,y)\in\Omega$ and some parameters $\bm{\mu}=(\mu_1,\mu_2) \in \mathcal{M}={[-0.5,0.5]}^p$, hence with $p=2$ parameters.
With an abuse of language as well,
we refer to the quantity $\kappa$ in~\eqref{eq:analytical_solution_Lap2D}
as the frequency of the solution, in the sense that it characterizes
the number of oscillations in the solution.

We start with a ``low frequency'' case in \cref{sec:Lap2Dlow}, taking $\kappa = 2$ and considering a PINN where we impose the Dirichlet boundary conditions as presented in \cref{sec:exact_imposition_of_BC}., i.e. using a level-set function.
To further improve the prior quality, we introduce an augmented loss function in \cref{sec:Lap2Dlowaug} by using the Sobolev training presented in \cref{sec:sobolev_training}.
Afterwards, we test another loss in \cref{sec:Lap2Dlowbc} that includes the Dirichlet condition, or in other words, that does not use a level-set function.
Finally, we consider a ``higher frequency''
case in \cref{sec:Lap2Dhigh}, with $\kappa = 8$.

\begin{rmrk}\label[rmrk]{rmk:Lap2D_N_nodes}
	In the following, the characteristic mesh size $h=\frac{\pi\sqrt{2}}{N-1}$ is defined as a function of $N$, considering a cartesian mesh of $N^2$ nodes for our squared 2D domain of length $\pi$.
\end{rmrk}

\subsubsection{Low-frequency case}\label{sec:Lap2Dlow}

We consider a ``low-frequency'' problem, taking $\kappa = 2$.
In this section, we consider the additive approach, as presented in \cref{sec:additive_prior}, by considering the PINN prior $u_\theta$. We start by testing the error estimates in (in $L^2$ norm \ROne{and $H^1$ semi-norm}) with polynomial order
$k\in\{1,2,3\}$, then we compare the different approaches. We evaluate the gains obtained on a sample of parameters. Then, we compare the numerical costs of the different methods. Finally we discuss the importance of integrating analytical functions, as presented in \cref{sec:using_PINN}.


Since the problem under consideration is parametric
we deploy a parametric PINN,
which depends on both the space variable $\bm{x}=(x,y) \in \Omega$
and the parameters $\bm{\mu}=(\mu_1,\mu_2) \in \mathcal{M}$.
Moreover, we strongly impose the Dirichlet boundary conditions,
as presented in \cref{sec:PINNs_parametric_PDE}.
To do this, we use the prior \RBoth{$u_\theta$ defined in~\eqref{eq:prior_with_levelset}, where we choose the level-set function $\varphi$ defined by}
\begin{equation*}
	\varphi(\bm{x})=(x+0.5\pi)(x-0.5\pi)(y+0.5\pi)(y-0.5\pi).
\end{equation*}
\RBoth{Thus, we will only consider the residual loss $J_r$ approached by a Monte-Carlo method as defined in~\eqref{eq:residual_loss_parametric_MC} with $N_\text{col}=\num{6000}$ collocation points uniformly chosen on $\Omega\times\mathcal{M}$. The parametric network} is defined as an MLP with the hyperparameters defined in \cref{tab:paramtest1_2D}; we use the Adam optimizer and then switch to the LBFGS optimizer after the $n_\text{switch}$-th epoch.

\begin{table}[htbp]
    \centering
    \begin{tabular}{cc}
        \toprule
        \multicolumn{2}{c}{\textbf{Network - MLP}} \\
        \midrule
        \textit{layers} & $40,60,60,60,40$ \\
        \cmidrule(lr){1-2}
        $\sigma$ & sine \\
        \bottomrule
    \end{tabular}
    \hspace{1cm}
    \begin{tabular}{cccc}
        \toprule
        \multicolumn{4}{c}{\textbf{Training - with LBFGS}} \\
        \midrule
        \textit{lr} & 1.7e-2 & $n_\text{epochs}$ & \num{5000} \\
        \cmidrule(lr){1-2} \cmidrule(lr){3-4}
        \textit{decay} & 0.99 & $n_\text{switch}$ & \num{1000} \\
        \cmidrule(lr){1-2} \cmidrule(lr){3-4}
        $N_\text{col}$ & \num{6000} \\
        \bottomrule
    \end{tabular}
    \hspace{1cm}
    \begin{tabular}{cccc}
        \toprule
        \multicolumn{4}{c}{\textbf{Loss weights}} \\
        \midrule
        $\omega_r$ & 1 & $\omega_\text{data}$ & 0 \\
        \cmidrule(lr){1-2} \cmidrule(lr){3-4}
        $\omega_b$ & 0 & $\omega_\text{sob}$ & 0 \\        
        \bottomrule
    \end{tabular}
    \caption{Network, training parameters (\cref{rmk:PINN_notations}) and loss weights for $u_\theta$ in the \textit{2D low-frequency case}.}\label{tab:paramtest1_2D}
\end{table}

\boldparagraph{Error estimates}

We start by testing the error estimates of \cref{lem:error_estimation_add} for the following two sets of parameters, randomly selected in $\mathcal{M}$:
\begin{equation*}
	\bm{\mu}^{(1)}=(0.05,0.22) \quad \text{and} \quad \bm{\mu}^{(2)}=(0.1,0.04)\;,
\end{equation*}
by considering the PINN prior $u_\theta$.
\RBoth{We perform the same tests as in the previous sections, but this time considering only the additive approach.} The results are presented in \ROne{\cref{fig:case1,fig:case1H1} (respectively in norm $L^2$ and semi-norm $H^1$)} for fixed $k \in \{1,2,3\}$ with $h$ depending on $N\in\{16,32,64,128,256\} $, as presented in \cref{rmk:Lap2D_N_nodes}.

\begin{figure}[ht!]
	\centering
	\begin{subfigure}{0.48\linewidth}
		\centering
		\cvgFEMCorrAlldeg{fig_testcase2D_test1_cvg_FEM_case1_v1_param1.csv}{fig_testcase2D_test1_cvg_Corr_case1_v1_param1.csv}{1e-10}
		\caption{Case of $\bm{\mu}^{(1)}$}
	\end{subfigure}
	\begin{subfigure}{0.48\linewidth}
		\centering
		\cvgFEMCorrAlldeg{fig_testcase2D_test1_cvg_FEM_case1_v1_param2.csv}{fig_testcase2D_test1_cvg_Corr_case1_v1_param2.csv}{1e-10}
		\caption{Case of $\bm{\mu}^{(2)}$}
	\end{subfigure}
	\caption{Considering the \textit{2D low-frequency case} and the PINN prior $u_\theta$. Left -- $L^2$ relative error on $h$, obtained with the standard FEM $e_h^{(1)}$ (solid lines) and the additive approach $e_{h,+}^{(1)}$ (dashed lines) for $\bm{\mu}^{(1)}$, with $k \in \{1,2,3\}$. Right -- Same for $\bm{\mu}^{(2)}$.}\label{fig:case1}
\end{figure}

\begin{figure}[ht!]
	\centering
	\begin{subfigure}{0.48\linewidth}
		\centering
		\cvgFEMCorrAlldegHun{fig_testcase2D_test1_cvgsemiH1_FEM_case1_v1_param1.csv}{fig_testcase2D_test1_cvgsemiH1_Corr_case1_v1_param1.csv}{1e-7}
		\caption{Case of $\bm{\mu}^{(1)}$}
	\end{subfigure}
	\begin{subfigure}{0.48\linewidth}
		\centering
		\cvgFEMCorrAlldegHun{fig_testcase2D_test1_cvgsemiH1_FEM_case1_v1_param2.csv}{fig_testcase2D_test1_cvgsemiH1_Corr_case1_v1_param2.csv}{1e-7}
		\caption{Case of $\bm{\mu}^{(2)}$}
	\end{subfigure}
	\caption{\ROne{Considering the \textit{2D low-frequency case} and the PINN prior $u_\theta$. Left -- Semi-$H^1$ relative error on $h$, obtained with the standard FEM (solid lines) and the additive approach (dashed lines) for $\bm{\mu}^{(1)}$, with $k \in \{1,2,3\}$. Right -- Same for $\bm{\mu}^{(2)}$.}}\label{fig:case1H1}
\end{figure}

In \ROne{\cref{fig:case1,fig:case1H1}}, we observe the expected behavior.
Indeed, the error decreases with the correct order of accuracy as the mesh size $h$ decreases \RTwo{(i.e., with a slope of $k+1$ in the $L^2$ norm and $k$ in the $H^1$ semi-norm)}.
This observation is valid for both the classical and enriched FEM.
Moreover, we observe that the error constant of the additive approach is significantly lower than that of the classical FEM. In \RTwo{the following paragraph}, we will compare these different approaches in more detail. As noted in the 1D case, we can see that the additive enriched approach for $k=1$ (resp. $k=2$) seems to give \ROne{relative errors close to the original FEM} for $k=2$ (resp. $k=3$)\ROne{, although the rate of convergence is different}. 

\boldparagraph{Comparison of different approaches}

\begin{figure}[ht!]
	\centering
    \includegraphics[scale=1]{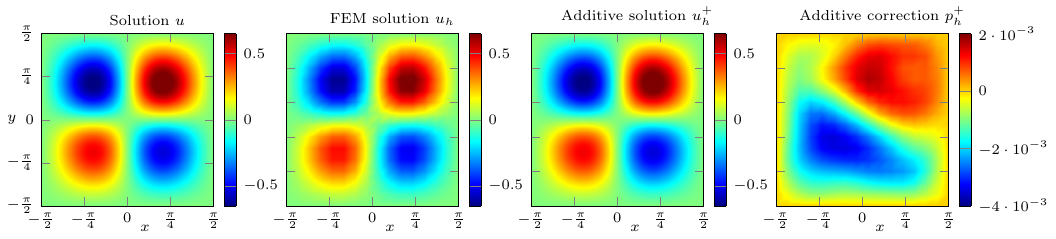}

    \includegraphics[scale=1]{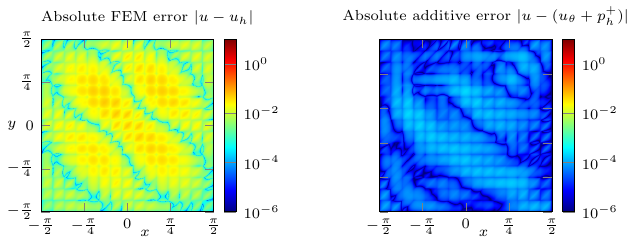}

	\caption{Considering the \textit{2D low-frequency case} with $\bm{\mu}^{(1)}$, $k=1$, $N=16$ and the PINN prior $u_\theta$. Comparison of the solution obtained with the standard FEM and the additive approach with the analytical solution. For the additive method, comparison of the correction term with the analytical one.}\label{fig:case1_2D_plots}
\end{figure}

We now focus on the first parameter $\bm{\mu}^{(1)}$. We compare the standard FEM with the additive approach, first in the \cref{tab:case1_2D_comparison} where we can see the different errors obtained with the different methods for $k=1$ fixed and $N\in\{16,32\}$ as well as the gains obtained in comparison with standard FEM. Next, we take a closer look at the solution obtained with the different approaches in \cref{fig:case1_2D_plots}; for each method, we compare the solution obtained ($u_h$ for standard FEM and $u_h^+$ for the additive approach) with the analytical solution $u$. For the enriched method, using the PINN prior $u_\theta$, we will also compare the proposed correction; namely, for the additive approach, we will compare $p_h^+$ with $u-u_\theta$.

\begin{table}[ht!]
	\centering
	\GainsFixedMu{1}{fig_testcase2D_test1_plots_FEM.csv}{fig_testcase2D_test1_plots_compare_gains.csv}
	\caption{Considering the \textit{2D low-frequency case} with $\bm{\mu}^{(1)}$, $k=1$ and $N\in\{16,32\}$. Left -- $L^2$ relative error obtained with FEM. Right -- Considering the PINN prior $u_\theta$, $L^2$ relative errors and gains with respect to FEM, obtained with the additive approach.}\label{tab:case1_2D_comparison}
\end{table}

In \cref{tab:case1_2D_comparison}, we observe that the additive approach significantly improves the error of the standard FEM, with gains of around $260$ for $N=16$ and $k=1$, which is equivalent to refining the mesh by a factor of $16$ for $\mathbb{P}_1$ elements. Indeed, in this case, our enriched approach gives much better results than standard FEM. In \cref{fig:case1_2D_plots}, we observe that the solution obtained with the additive approach is very close to the analytical solution, with a correction term that is also very close to the analytical one. This shows the effectiveness of the additive approach in this case.

\boldparagraph{Gains achieved with the additive approach}

Considering a \RTwo{set $\mathcal{S}$ of $n_p=50$ parameter instances},
we now evaluate the gains $G_{+,\theta}$ and $G_+$ defined in~\eqref{eq:gain_add_num}.
The results are presented in \cref{tab:case1_2D}
for $k \in \{1,2,3\}$ and $N \in \{20, 40\}$.

\begin{table}[ht!]
	\centering
	\GainsTableAlldeg{fig_testcase2D_test1_gains_Tab_stats_case1_v1.csv}
	\caption{Considering the \textit{2D low-frequency case}, $k\in\{1,2,3\}$ and the PINN prior $u_\theta$. Left -- Gains in $L^2$ relative error of the additive method with respect to PINN. Right -- Gains in $L^2$ relative error of our approach with respect to FEM.}\label{tab:case1_2D}
\end{table}

In \cref{tab:case1_2D}, we observe (left subtable) that our method
significantly improves the error of the PINN,
especially for large values of $k$,
where the enrichment is performed in a richer approximation space.
Moreover, we also observe (right subtable) significant
gains with respect to classical FEM.
For instance, as expected from the results of \RTwo{the previous paragraph}, the mean gains for $k=1$ are around $270$, which corresponds to refining the mesh approximately $16$ times for $\mathbb{P}_1$ elements. This means that our $\mathbb{P}_1$ enhanced bases capture the solution as accurately as classical $\mathbb{P}_1$ bases with a mesh sixteen times finer.
For $k=2$ and $k=3$, the mean gains are around $134$ and $61$, respectively, which corresponds to refining the mesh approximately $5$ times for $\mathbb{P}_2$ elements and $2.8$ times for $\mathbb{P}_3$ elements. A natural follow-up question consists in assessing the impact of the PINN quality on our results. This will be the subject of the \cref{sec:Lap2Dlowaug}.

\boldparagraph{Costs of the different methods}

To more accurately assess the benefits of using the enriched methods, we look in this section at the costs of the different methods proposed, considering the parameter $\bm{\mu}^{(1)}$. Thus, we will consider that the cost of using the PINN prior $u_\theta$, corresponds to the total number of weights of the network considered. In this case, it is given as an MLP with the hyperparameters defined in \cref{tab:paramtest1_2D}, for a total of $N_\text{weights}=\num{12461}$ weights.

For the different finite element methods, we endeavor to determine,
for a fixed polynomial degree $k$, the characteristic mesh size $h$
(depending on $N$ as described in \cref{rmk:Lap2D_N_nodes})
required to reach a given fixed error $e$.
In \cref{tab:case1_2D_costs}, we study, for $k\in\{1,2,3\}$, considering standard FEM and the additive approach, the $N$ required to achieve the same error $e$. More precisely, the characteristic mesh size required by standard FEM so that \smash{$e_h^{(1)}\approx e$} and the one required by the additive approach so that \smash{$e_{h,+}^{(1)}\approx e$}. Depending on the polynomial degree $k$, we can also determine the number of degrees of freedom $N_\text{dofs}$ associated with each case. \RTwo{These results are obtained by interpolating the convergence curves of \cref{fig:case1} for the different methods for a given $e$.}

\begin{table}[ht!]
	\centering
	\coststableallq{fig_testcase2D_test1_costs_TabDoFs_case1_v1_param1.csv}
	\caption{Considering the \textit{2D low-frequency case} with $\bm{\mu}^{(1)}$, $k\in\{1,2,3\}$ and the PINN prior $u_\theta$. Left -- Characteristic $N$ (associated to the characteristic mesh size $h$) required to reach a fixed error $e$ for standard FEM and the additive approach. Right -- Number of degrees of freedom $N_\text{dofs}$ associated with each case.}\label{tab:case1_2D_costs}
\end{table}

In \cref{tab:case1_2D_costs}, we see that the additive approach proposed in \cref{sec:additive_prior} requires a much coarser mesh than standard FEM to achieve the same $e$ error. This is due to the error estimations of \cref{lem:error_estimation_add} which show that the error of the enhanced FEM is significantly lower than that of the classical FEM (depending on the quality of the prior). This is also reflected in the number of degrees of freedom required to achieve the same error $e$.

However, the enriched approaches proposed require using the prior PINN $u_\theta$, which also includes its inference cost. For this reason, it is also interesting to study these same costs on a set of parameters, say of size $n_p=100$. Since we are in the context of parametric PINN, we can estimate that the computational cost of solving~\eqref{eq:Lap2D} on this sample of $n_p$ parameters corresponds, for the additive approach, to $n_p$ times its number of dofs plus the cost of using PINN (i.e. its total number of weights), thus $n_p\times N_\text{dofs}+N_\text{weights}$ (with $N_\text{dofs}$ the number of dofs associated to the additive approach). For standard FEM, this cost is equivalent to $n_p$ times its estimated number of dofs $n_p\times N_\text{dofs}$ (with $N_\text{dofs}$ the number of dofs associated with standard FEM). We will then compare these costs for a set of $n_p=100$ parameters, considering the same error $e$ to be achieved for both methods. The results are presented in \cref{tab:case1_2D_costs100}. \RTwo{Note that these results (right subtable) are, in fact, only an estimate of the real cost of solving $n_p$ problems. In practice, the number of degrees of freedom $N_\text{dofs}$ associated with each method depends on the parameter itself. The error to be achieved $e$ will require more or less fine meshes for each parameter.}

\begin{table}[ht!]
	\centering
	\coststableallqhundred{fig_testcase2D_test1_costs_TabDoFsParam_case1_v1_param1_nparams100.csv}
	\caption{Considering the \textit{2D low-frequency case}, $k\in\{1,2,3\}$ and the PINN prior $u_\theta$. Left -- Total costs of standard FEM and the additive approach to reach an error $e$ for a set of $n_p=1$ parameter. Right -- Same for a set of $n_p=100$ parameters.}\label{tab:case1_2D_costs100}
\end{table}

In \cref{tab:case1_2D_costs100}, for $n_p=1$ (left subtable), we can see that the cost of the additive method is generally lower than that of standard FEM, even though they are of the same range. However, it is important to note that these are not entirely comparable: in fact, a large part of the cost of the additive method lies in PINN prediction, which will be more or less well estimated depending on a number of hyper-parameters (number of epochs, learning rate, etc.). If we then take $n_p=100$ (right subtable), we can see that the cost of standard FEM becomes radically higher than with the additive approach. This is why the improved approach is particularly interesting for solving the~\eqref{eq:Lap2D} problem on a set of parameters.

\boldparagraph{Integration of analytical functions}

This section aims to discuss one of the important points, presented in \cref{sec:using_PINN}, that enables PINN to be used effectively and can make our enriched methods more or less effective.
Indeed, according to~\eqref{eq:approachform_add}, we have to integrate $f+\Delta u_\theta$ multiplied by the test function.
To perform this integration, we first interpolate this term on a polynomial space and then integrate it exactly.
The degree of this polynomial approximation is an important parameter to make our technique effective.

\begin{figure}[ht!]
	\centering
	\includegraphics[scale=0.8]{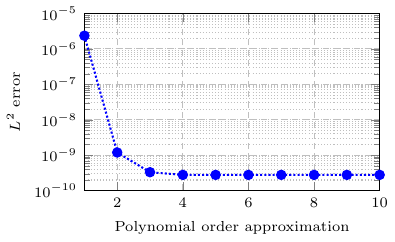}
	\caption{Considering the \textit{2D low-frequency case} with $\bm{\mu}^{(1)}$, $k=3$, $N=128$ and the PINN prior $u_\theta$. Considering the additive approach, $L^2$ error $e_{h,+}^{(1)}$ with respect to the degree of polynomial approximation of $f+\Delta u_\theta$.}\label{fig:case1_2D_highdeg}
\end{figure}

The goal here is simply to show that for enriched approaches to be effective, particularly the additive method, this polynomial approximation must be of a sufficiently high degree.
Consider the parameter $\bm{\mu}^{(1)}$, a polynomial degree $k=3$ and a number of nodes $N^2$ with $N=128$.
In \cref{fig:case1_2D_highdeg}, we display the $L^2$ \RBoth{relative} error of the additive approach with respect to the degree of polynomial approximation of $f+\Delta u_\theta$.

In \cref{fig:case1_2D_highdeg}, we observe that the error decreases as the degree of polynomial approximation increases. This shows the importance of properly interpolating analytical functions in the context of enriched methods.

\subsubsection{Low-frequency case --- Sobolev training}\label{sec:Lap2Dlowaug}

This section focuses on the same problem as in \cref{sec:Lap2Dlow}. The aim here is to show that the network quality has a non-negligible impact on the results obtained with our method and that if the network is better, our results will be, too. To this end, we defined a new prior $u_\theta^\text{sob}$ by using the Sobolev training presented in \cref{sec:sobolev_training}, where the derivatives of the solution should be better approximated than by a standard training and compared it with the PINN prior $u_\theta$ defined in \cref{sec:Lap2Dlow}. We start by testing the error estimation with $k\in\{1,2,3\}$ polynomial order and evaluate the gains obtained on the same sample of parameters as in \cref{sec:Lap2Dlow}.


We deploy here \RBoth{a parametric PINN, denoted by $u_\theta^\text{sob}$, where we consider the residual loss $J_r$ and the Sobolev loss $J_\text{\rm sob}$, respectively defined in~\eqref{eq:residual_loss_parametric} and ~\eqref{eq:sobolev_loss}, whose integrals are both approximated by a Monte-Carlo method with $N_\text{col}=\num{6000}$ collocation points. We then solve the minimisation problem defined in~\eqref{eq:minimization_problem_sobolev} considering that the Dirichlet boundary conditions are strongly imposed as in \cref{sec:Lap2Dlow}.} The hyperparameters are defined in \cref{tab:paramtest1v7_2D}; we use the Adam optimizer.

\begin{rmrk}
	Adding the Sobolev loss can make training more difficult, so we only consider \num{3000} epochs but a batch size of \num{2000}. This means that for each epoch, the weights will be updated 3 times (because $N_\text{col}=\num{6000}$).
\end{rmrk}

\begin{table}[htbp]
    \centering
    \begin{tabular}{cc}
        \toprule
        \multicolumn{2}{c}{\textbf{Network - MLP}} \\
        \midrule
        \textit{layers} & $40,60,60,60,40$ \\
        \cmidrule(lr){1-2}
        $\sigma$ & sine \\
        \bottomrule
    \end{tabular}
    \hspace{1cm}
    \begin{tabular}{cccc}
        \toprule
        \multicolumn{4}{c}{\textbf{Training - with LBFGS}} \\
        \midrule
        \textit{lr} & 1.7e-2 & $n_\text{epochs}$ & \num{3000} \\
        \cmidrule(lr){1-2} \cmidrule(lr){3-4}
        \textit{decay} & 0.99 & batch size & \num{2000} \\
        \cmidrule(lr){1-2} \cmidrule(lr){3-4}
        $N_\text{col}$ & \num{6000} \\
        \bottomrule
    \end{tabular}
    \hspace{1cm}
    \begin{tabular}{cccc}
        \toprule
        \multicolumn{4}{c}{\textbf{Loss weights}} \\
        \midrule
        $\omega_r$ & 1 & $\omega_\text{data}$ & 0 \\
        \cmidrule(lr){1-2} \cmidrule(lr){3-4}
        $\omega_b$ & 0 & $\omega_\text{sob}$ & 0.1 \\        
        \bottomrule
    \end{tabular}
    \caption{Network, training parameters (\cref{rmk:PINN_notations}) and loss weights for $u_\theta^\text{sob}$ in the \textit{2D low-frequency case}.}\label{tab:paramtest1v7_2D}
\end{table}

\boldparagraph{Error estimates}

For simplicity, we consider the first parameter $\bm{\mu}^{(1)}=(0.05,0.22)$ presented in \cref{sec:Lap2Dlow}. 
\RBoth{We perform the same test as in the previous section, considering the two priors, $u_\theta$ and $u_\theta^\text{sob}$, in the enriched approach.}
The results are presented in \cref{fig:case1v7} for fixed $k\in \{1,2,3\}$ \RBoth{(in $L^2$ norm)}.

\begin{figure}[ht!]
	\centering
	\hspace{-1cm}
	\begin{subfigure}{0.32\linewidth}
		\centering
		\resizebox{1.1\linewidth}{!}{
			\cvgFEMCorrTwoPriors{1}{fig_testcase2D_test1_cvg_FEM_case1_v1_param1.csv}{fig_testcase2D_test1_cvg_Corr_case1_v1_param1.csv}{fig_testcase2D_test1_v7_cvg_Corr_case1_v7_param1.csv}{8e-7}{$u_\theta^\text{sob}$}
		}
		\caption{$k=1$}
	\end{subfigure}
	\begin{subfigure}{0.32\linewidth}
		\centering
		\resizebox{1.1\linewidth}{!}{
			\cvgFEMCorrTwoPriors{2}{fig_testcase2D_test1_cvg_FEM_case1_v1_param1.csv}{fig_testcase2D_test1_cvg_Corr_case1_v1_param1.csv}{fig_testcase2D_test1_v7_cvg_Corr_case1_v7_param1.csv}{4e-9}{$u_\theta^\text{sob}$}
		}
		\caption{$k=2$}
	\end{subfigure}
	\begin{subfigure}{0.32\linewidth}
		\centering
		\resizebox{1.1\linewidth}{!}{
			\cvgFEMCorrTwoPriors{3}{fig_testcase2D_test1_cvg_FEM_case1_v1_param1.csv}{fig_testcase2D_test1_cvg_Corr_case1_v1_param1.csv}{fig_testcase2D_test1_v7_cvg_Corr_case1_v7_param1.csv}{4e-11}{$u_\theta^\text{sob}$}{3}
		}
		\caption{$k=3$}
	\end{subfigure}
	\caption{Considering the \textit{2D low-frequency case} with $\bm{\mu}^{(1)}$. Left -- $L^2$ relative error on $h$, obtained with the standard FEM $e_h^{(1)}$ (solid line) and the additive approach $e_{h,+}^{(1)}$ (dashed lines), with $k=1$, by considering the PINN prior with standard training $u_\theta$ and Sobolev training $u_\theta^\text{sob}$. Middle -- Same with $k=2$. Right -- Same with $k=3$.}
	\label{fig:case1v7}
\end{figure}

We observe that Sobolev training improves the results obtained with the $L^2$ training,
for $k\in\{1,2,3\}$.
This shows the impact of the quality of the network prediction on our method.
To further investigate this, we evaluate the gains obtained with the Sobolev training on the same sample of parameters as in \cref{sec:Lap2Dlow}.

\boldparagraph{Gains achieved with the additive approach}

Considering the same \RTwo{set $\mathcal{S}$ of $n_p=50$ parameter instances} as in \cref{sec:Lap2Dlow}, we now evaluate the gains $G_{+,\theta}$ and $G_+$ defined in~\eqref{eq:gain_add_num} considering the PINN prior $u_\theta^\text{sob}$ using Sobolev training. The results are presented in \cref{tab:case1v7} for $k \in \{1,2,3\}$ and $N \in \{20,40\}$.

\begin{table}[ht!]
	\centering
	\GainsTableAlldeg{fig_testcase2D_test1_v7_gains_Tab_stats_case1_v7.csv}
	\caption{Considering the \textit{2D low-frequency case}, $k\in\{1,2,3\}$ and the PINN prior $u_\theta^\text{sob}$ (Sobolev training). Left -- Gains in $L^2$ relative error of the additive method with respect to PINN. Right -- Gains in $L^2$ relative error of our approach with respect to FEM.}\label{tab:case1v7}
\end{table}

The gains reported in \cref{tab:case1v7} show that, compared to $L^2$ training,
Sobolev training increases the mean gains by a factor of about $3$.
This corresponds to almost half an additional mesh refinement for the $\mathbb{P}_1$ elements. We also note that this Sobolev training is particularly interesting for higher polynomial degrees, with standard $L^2$ training having lower gains than for $k=1$.
\subsubsection{Low-frequency case --- Boundary loss training}\label{sec:Lap2Dlowbc}

\RBoth{In this section}, we focus on the same problem as in \cref{sec:Lap2Dlow} and \cref{sec:Lap2Dlowaug}. We now turn to a standard PINN, denoted by $u_\theta^\text{bc}$, where we impose the boundary conditions in the loss function (no longer with the level-set function). The aim here is to show that our enriched methods also work with priors that do not have exact boundary conditions. To this end, we start by testing the error estimation and evaluate the gains obtained on the same sample of parameters as in \cref{sec:Lap2Dlow}.

\begin{figure}[ht!]
    \centering
    \includegraphics[width=0.3\textwidth]{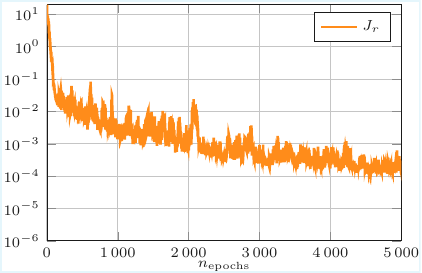} \hspace{0.5cm}
    \includegraphics[width=0.3\textwidth]{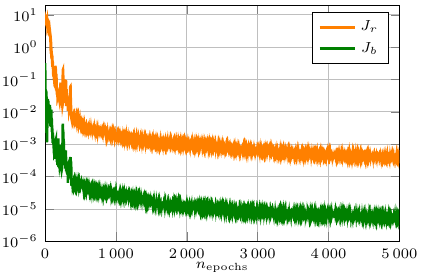}
    \caption{\RTwoN{Training loss curves obtained in the \textit{2D low-frequency case}. Left -- Considering the PINN prior $u_\theta$ with strong boundary conditions (\cref{sec:Lap2Dlow}). Right -- Considering the PINN prior $u_\theta^\text{bc}$ with boundary loss training (\cref{sec:Lap2Dlowbc})	.}}
    \label{fig:loss_curves}
\end{figure}

\RBoth{We deploy a new parametric PINN, denoted $u_{\theta}^\text{bc}$, where we consider the residual loss function $J_r$ and the boundary loss function $J_b$, respectively defined in~\eqref{eq:residual_loss_parametric_MC} and~\eqref{eq:boundary_loss_parametric_MC}. These integrals are both approximated by a Monte-Carlo method considering $N_\text{col}=\num{6000}$ collocation points and $N_\text{bc}=\num{2000}$ boundary collocation points. We insist that the collocation points are resampled at each epoch during the training process. Here, we no longer strongly impose the Dirichlet boundary conditions. The hyperparameters are defined in \cref{tab:paramtest1v3_2D}; we use the Adam optimizer and then switch to the LBFGS optimizer after the $n_\text{switch}$-th epoch.} \RTwo{The training loss curves are presented in \cref{fig:loss_curves} for both priors $u_\theta$ and $u_\theta^\text{bc}$.}

\begin{table}[htbp]
    \centering
    \begin{tabular}{cc}
        \toprule
        \multicolumn{2}{c}{\textbf{Network - MLP}} \\
        \midrule
        \textit{layers} & $40,60,60,60,40$ \\
        \cmidrule(lr){1-2}
        $\sigma$ & sine \\
        \bottomrule
    \end{tabular}
    \hspace{1cm}
    \begin{tabular}{cccc}
        \toprule
        \multicolumn{4}{c}{\textbf{Training - with LBFGS}} \\
        \midrule
        \textit{lr} & 1.7e-2 & $n_\text{epochs}$ & \num{5000} \\
        \cmidrule(lr){1-2} \cmidrule(lr){3-4}
        \textit{decay} & 0.99 & $n_\text{switch}$ & \num{1000} \\
        \cmidrule(lr){1-2} \cmidrule(lr){3-4}
        $N_\text{col}$ & \num{6000} & $N_\text{bc}$ & \num{2000} \\
        \bottomrule
    \end{tabular}
    \hspace{1cm}
    \begin{tabular}{cccc}
        \toprule
        \multicolumn{4}{c}{\textbf{Loss weights}} \\
        \midrule
        $\omega_r$ & 1 & $\omega_\text{data}$ & 0 \\
        \cmidrule(lr){1-2} \cmidrule(lr){3-4}
        $\omega_b$ & 30 & $\omega_\text{sob}$ & 0 \\        
        \bottomrule
    \end{tabular}
    \caption{Network, training parameters (\cref{rmk:PINN_notations}) and loss weights for $u_\theta^\text{bc}$ in the \textit{2D low-frequency case}.}\label{tab:paramtest1v3_2D}
\end{table}

\boldparagraph{Error estimates}

We consider the first parameter $\bm{\mu}^{(1)}=(0.05,0.22)$ presented in \cref{sec:Lap2Dlow} \RBoth{and perform the same test as before, considering the two priors $u_\theta$ and $u_\theta^\text{bc}$.} The results are presented in \cref{fig:case1v3} for fixed $k\in \{1,2,3\}$ \RBoth{(in $L^2$ norm)}.

\begin{figure}[H]
	\centering
	\hspace{-1cm}
	\begin{subfigure}{0.32\linewidth}
		\centering
		\resizebox{1.1\linewidth}{!}{
			\cvgFEMCorrTwoPriors{1}{fig_testcase2D_test1_cvg_FEM_case1_v1_param1.csv}{fig_testcase2D_test1_cvg_Corr_case1_v1_param1.csv}{fig_testcase2D_test1_v3_cvg_Corr_case1_v3_param1.csv}{3e-6}{$u_\theta^\text{bc}$}
		}
		\caption{$k=1$}
	\end{subfigure}
	\begin{subfigure}{0.32\linewidth}
		\centering
		\resizebox{1.1\linewidth}{!}{
			\cvgFEMCorrTwoPriors{2}{fig_testcase2D_test1_cvg_FEM_case1_v1_param1.csv}{fig_testcase2D_test1_cvg_Corr_case1_v1_param1.csv}{fig_testcase2D_test1_v3_cvg_Corr_case1_v3_param1.csv}{1e-8}{$u_\theta^\text{bc}$}
		}
		\caption{$k=2$}
	\end{subfigure}
	\begin{subfigure}{0.32\linewidth}
		\centering
		\resizebox{1.1\linewidth}{!}{
			\cvgFEMCorrTwoPriors{3}{fig_testcase2D_test1_cvg_FEM_case1_v1_param1.csv}{fig_testcase2D_test1_cvg_Corr_case1_v1_param1.csv}{fig_testcase2D_test1_v3_cvg_Corr_case1_v3_param1.csv}{8e-11}{$u_\theta^\text{bc}$}
		}
		\caption{$k=3$}
	\end{subfigure}
	\caption{Considering the \textit{2D low-frequency case} with $\bm{\mu}^{(1)}$. Left -- $L^2$ relative error on $h$, obtained with the standard FEM \smash{$e_h^{(1)}$} (solid line) and the additive approach \smash{$e_{h,+}^{(1)}$} (dashed lines), with $k=1$, by considering the PINN prior with standard training $u_\theta$ and the BC loss training \smash{$u_\theta^\text{bc}$}. Middle -- Same with $k=2$. Right -- Same with $k=3$.}\label{fig:case1v3}
\end{figure}

We can see in \cref{fig:case1v3} that the additive approach also works when the prior is not exact on the boundary, as here with $u_\theta^\text{bc}$. In particular, for $k\in\{1,2,3\}$ and the parameter $\bm{\mu}^{(1)}$, our enriched approach using the prior $u_\theta^\text{bc}$ seems to give very similar results to those obtained with $u_\theta$ even if the approach with level-set is, for this case, slightly better.

\boldparagraph{Gains achieved with the additive approach}

Considering the same \RTwo{set $\mathcal{S}$ of $n_p=50$ parameter instances} as in \cref{sec:Lap2Dlow}, we now evaluate the gains $G_{+,\theta}$ and $G_+$ defined in~\eqref{eq:gain_add_num} considering the PINN prior $u_\theta^\text{bc}$ using BC loss training. The results are presented in \cref{tab:case1v3} for $k \in \{1,2,3\}$ and $N \in \{20,40\}$.

\begin{table}[ht!]
	\centering
	\GainsTableAlldeg{fig_testcase2D_test1_v3_gains_Tab_stats_case1_v3.csv}
	\caption{Considering the \textit{2D low-frequency case}, $k\in\{1,2,3\}$ and the PINN prior $u_\theta^\text{bc}$ (BC loss training). Left -- Gains in $L^2$ relative error of the additive method with respect to PINN. Right -- Gains in $L^2$ relative error of our approach with respect to FEM.}\label{tab:case1v3}
\end{table}

The gains reported in \cref{tab:case1v3} show that for this test case, the use of the prior $u_\theta$ (using the level-set) in our enriched approach, seems to give better gains than those of \cref{tab:case1_2D}, considering the current prior $u_\theta^\text{bc}$. This may be due to the addition of the $\omega_\text{bc}$ hyperparameter for balancing losses in training, which may make training less efficient. However, the results obtained with the current prior $u_\theta^\text{bc}$ are still very good, and the gains are still significant compared to the standard FEM.

\subsubsection{High-frequency case}\label{sec:Lap2Dhigh}

To increase in complexity, we investigate a higher-frequency problem by taking $\kappa=8$ in~\eqref{eq:analytical_solution_Lap2D}.
In this section, we start by testing the error estimates. Then, we compare the different methods and evaluate the gains obtained on a sample of parameters.


This time, we use the Fourier features from \cite{TanSri2020} as presented in \cref{sec:spectral_bias} to construct the PINN prior $u_\theta$. The hyperparameters are defined in \cref{tab:paramtest2_2D};
we use the Adam optimizer and then switch to the LBFGS optimizer after the $n_\text{switch}$-th epoch. We consider $N_\text{col}=\num{6000}$ collocation points, uniformly chosen on $\Omega$. We impose the Dirichlet boundary conditions as in \cref{sec:Lap2Dlow} using the level-set function and the same residual loss.

\begin{table}[htbp]
    \centering
    \begin{tabular}{cc}
        \toprule
        \multicolumn{2}{c}{\textbf{Network - MLP w/ FF}} \\
        \midrule
        \textit{layers} & $40,60,60,60,40$ \\
        \cmidrule(lr){1-2}
        $\sigma$ & sine \\
        \cmidrule(lr){1-2}
        $n_f$ & 40 \\
        \bottomrule
    \end{tabular}
    \hspace{1cm}
    \begin{tabular}{cccc}
        \toprule
        \multicolumn{4}{c}{\textbf{Training - with LBFGS}} \\
        \midrule
        \textit{lr} & 1.7e-2 & $n_\text{epochs}$ & \num{20000} \\
        \cmidrule(lr){1-2} \cmidrule(lr){3-4}
        \textit{decay} & 0.99 & $n_\text{switch}$ & \num{1000} \\
        \cmidrule(lr){1-2} \cmidrule(lr){3-4}
        $N_\text{col}$ & \num{6000} \\
        \bottomrule
    \end{tabular}
    \hspace{1cm}
    \begin{tabular}{cccc}
        \toprule
        \multicolumn{4}{c}{\textbf{Loss weights}} \\
        \midrule
        $\omega_r$ & 1 & $\omega_\text{data}$ & 0 \\
        \cmidrule(lr){1-2} \cmidrule(lr){3-4}
        $\omega_b$ & 0 & $\omega_\text{sob}$ & 0 \\        
        \bottomrule
    \end{tabular}
    \caption{Network, training parameters (\cref{rmk:PINN_notations}) and loss weights for $u_\theta$ in the \textit{2D high-frequency case}.}\label{tab:paramtest2_2D}
\end{table}

\boldparagraph{Error estimates}

We perform the same test as in \cref{sec:Lap2Dlow}, with a standard training (Sobolev training is not considered).
The results are displayed in \cref{fig:case2},
where we observe the expected behavior.
Indeed, all schemes have the correct order of accuracy,
and the enhanced FEM has a significantly lower error constant than the classical FEM.

\begin{figure}[ht!]
	\centering
	\begin{subfigure}{0.48\linewidth}
		\centering
		\cvgFEMCorrAlldeg{fig_testcase2D_test2_cvg_FEM_case2_v1_param1.csv}{fig_testcase2D_test2_cvg_Corr_case2_v1_param1.csv}{3e-8}
		\caption{$\bm{\mu}^{(1)}$ parameter.}
	\end{subfigure}
	\begin{subfigure}{0.48\linewidth}
		\centering
		\cvgFEMCorrAlldeg{fig_testcase2D_test2_cvg_FEM_case2_v1_param2.csv}{fig_testcase2D_test2_cvg_Corr_case2_v1_param2.csv}{3e-8}
		\caption{$\bm{\mu}^{(2)}$ parameter.}
	\end{subfigure}
	\caption{Considering the \textit{2D high-frequency case} and the PINN prior $u_\theta$. Left -- $L^2$ relative error on $h$, obtained with the standard FEM $e_h^{(1)}$ (solid lines) and the additive approach $e_{h,+}^{(1)}$ (dashed lines) for $\bm{\mu}^{(1)}$, with $k \in \{1,2,3\}$. Right -- Same for $\bm{\mu}^{(2)}$.}\label{fig:case2}
\end{figure}

\boldparagraph{Comparison of different approaches}

\begin{figure}[ht!] \centering

    \includegraphics[scale=1]{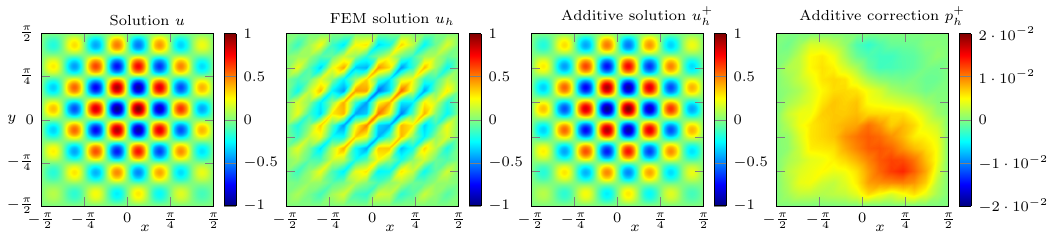}

    \includegraphics[scale=1]{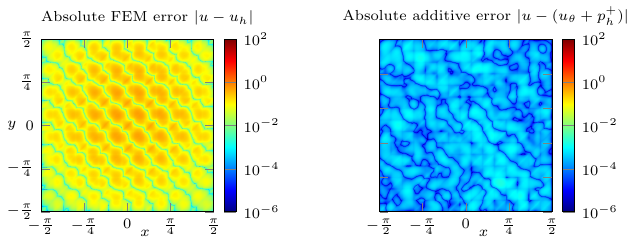}

	\caption{Considering the \textit{2D high-frequency case} with $\bm{\mu}^{(1)}$, $k=1$, $N=16$ and the PINN prior $u_\theta$. Comparison of the solution obtained with the standard FEM and the additive approach with the analytical solution. For the additive method, comparison of the correction term with the analytical one.}\label{fig:case2_2D_plots}
\end{figure}

\RBoth{Considering the first parameter $\bm{\mu}^{(1)}$, we} perform the same comparison as in \cref{sec:Lap2Dlow} for the high-frequency case. \RBoth{The results are presented in \cref{tab:case2_2D_comparison} and \cref{fig:case2_2D_plots}.}

\begin{table}[ht!]
	\centering
	\GainsFixedMu{1}{fig_testcase2D_test2_plots_FEM.csv}{fig_testcase2D_test2_plots_compare_gains.csv}
	\caption{Considering the \textit{2D high-frequency case} with $\bm{\mu}^{(1)}$, $k=1$ and $N\in\{16,32\}$. Left -- $L^2$ relative error obtained with FEM. Right -- Considering the PINN prior $u_\theta$, $L^2$ relative errors and gains with respect to FEM, obtained with the additive approach.}\label{tab:case2_2D_comparison}
\end{table}

We can see here that the gains obtained in \cref{tab:case2_2D_comparison} are much better than for the ``low frequency'' case presented in \cref{sec:Lap2Dlow}. This is, in fact, due to FEM's difficulty in approximating the solution for high frequencies, especially on coarse meshes. In fact, for the same choice of parameters, the FEM error on this high-frequency problem is 10 times worse than on the low-frequency one, which explains why our gains are so much greater. This also makes the use of the proposed enriched methods particularly interesting.
Moreover, we note that while FEM provides a reasonable approximation of the mean of the solution (as evidenced by the second figure on the top row of \cref{fig:case2_2D_plots}), it is unable to correctly resolve the small-scale oscillating behavior of the solution.
The additive correction restores this ability, and the new solution (third figure on the top row of \cref{fig:case2_2D_plots}) is much better able to capture the oscillations.

\boldparagraph{Gains achieved with the additive approach}

We now evaluate the gains $G_{+,\theta}$ and $G_+$, defined in~\eqref{eq:gain_add_num},
using the same \RTwo{set $\mathcal{S}$ of $n_p=50$ parameter instances}. The results are reported in \cref{tab:case2} for $k \in \{1,2,3\}$ and $N \in \{20, 40\}$.

\begin{table}[ht!]
	\centering
	\GainsTableAlldeg{fig_testcase2D_test2_gains_Tab_stats_case2_v1.csv}
	\caption{Considering the \textit{2D high-frequency case}, $k\in\{1,2,3\}$ and the PINN prior $u_\theta$. Left -- Gains in $L^2$ relative error of the additive method with respect to PINN. Right -- Gains in $L^2$ relative error of our approach with respect to FEM.}\label{tab:case2}
\end{table}

The same results can be observed in \cref{tab:case2} as for the $\bm{\mu}^{(1)}$ parameter. These could be improved by considering a Sobolev training, as in the ``low-frequency'' case presented in \cref{sec:Lap2Dlowaug}.

\subsection{2D anisotropic elliptic problem on a square}\label{sec:Ell2D}

In this section, we will consider the~\eqref{eq:ob_pde} problem in a more complex form than in \cref{sec:Lap2D}, by considering the following elliptic problem with homogeneous Dirichlet boundary conditions, in the 2D case ($d=2$),
\begin{equation*}
	\left\{
	\begin{aligned}
		-\text{div}(D\nabla u) & = f, \; &  & \text{in } \; \Omega, \\
		u         & =0, \;  &  & \text{on } \; \partial\Omega,
	\end{aligned}
	\right.
\end{equation*}
with $\Omega=\RBothN{(0,1)}^2$, $\partial\Omega$ its boundary and $\mathcal{M} \subset \mathbb{R}^p$ the parameter space (with $p$ the number of parameters).
Considering $\bm{x}=(x,y)\in\Omega$, we define $p=4$ parameters $\bm{\mu}=(\mu_1,\mu_2,\epsilon,\sigma)\in\mathcal{M}=[0.4, 0.6]\times [0.4, 0.6]\times [0.01,1]\times [0.1,0.8]$. We define the (symmetric and positive definite) diffusion matrix $D$ by
\begin{equation*}
	D(\bm{x},\bm{\mu})=\begin{pmatrix}
		\epsilon x^2+y^2 & (\epsilon-1)xy \\
		(\epsilon-1)xy & x^2+\epsilon y^2
	\end{pmatrix}
\end{equation*}
and the right-hand side $f$ by
\begin{equation*}
	f(\bm{x},\bm{\mu})=\exp\left(-\frac{{(x-\mu_1)}^2+{(y-\mu_2)}^2}{0.025\sigma^2}\right).
\end{equation*}
Note that the matrix $D$ has eigenvalues $x^2 + y^2$
and $\epsilon(x^2 + y^2)$,
leading to a diffusion process whose anisotropy
increases as $\epsilon$ decreases.

In this section, we consider the following three sets of parameters:
\begin{equation*}
	\bm{\mu}^{(1)}=(0.51,0.54,0.52,0.55), \quad \bm{\mu}^{(2)}=(0.48,0.53,0.41,0.89) \quad \text{and} \quad \bm{\mu}^{(3)}=(0.46,0.52,0.05,0.12).
\end{equation*}

In \cref{sec:Ell2D_error_estimates}, we start by testing the error estimates of \cref{lem:error_estimation_add}. Then, we compare the two different approaches in \cref{sec:Ell2D_comparison}. Finally, we evaluate the gains achieved with the additive approach as in \cref{sec:Ell2D_gains}.

\begin{rmrk}\label[rmrk]{rmk:Ell2D_N_nodes}
	In the following, the characteristic mesh size $h=\frac{\sqrt{2}}{N-1}$ is defined as a function of $N$, considering a cartesian mesh of $N^2$ nodes.
\end{rmrk}

We consider a parametric PINN where we exactly impose the Dirichlet boundary conditions as presented in \cref{sec:exact_imposition_of_BC}. \RBoth{Thus, we construct $u_\theta$ as in \eqref{eq:prior_with_levelset} with the level-set function $\varphi$ defined by}
\[
    \varphi(\bm{x})=x(x-1)y(y-1).
\]
Since we impose the boundary conditions by using the level-set function, we will only consider the residual loss $J_r$ \RBoth{approached by a Monte-Carlo method as defined in~\eqref{eq:residual_loss_parametric_MC} with $N_\text{col}=\num{8000}$ collocation points uniformly chosen on $\Omega\times\mathcal{M}$.}
The hyperparameters are given in \cref{tab:paramtest3_2D};
we use the Adam optimizer~\cite{KinBa2015}.

\begin{table}[htbp]
    \centering
    \begin{tabular}{cc}
        \toprule
        \multicolumn{2}{c}{\textbf{Network - MLP}} \\
        \midrule
        \textit{layers} & $40, 60, 60, 60, 40$ \\
        \cmidrule(lr){1-2}
        $\sigma$ & tanh \\
        \bottomrule
    \end{tabular}
    \hspace{1cm}
    \begin{tabular}{cccc}
        \toprule
        \multicolumn{4}{c}{\textbf{Training}} \\
        \midrule
        \textit{lr} & 1.6e-2 & $n_{epochs}$ & \num{15000} \\
        \cmidrule(lr){1-2} \cmidrule(lr){3-4}
        \textit{decay} & 0.99 \\
        \cmidrule(lr){1-2}
        $N_\text{col}$ & \num{8000} \\
        \bottomrule
    \end{tabular}
    \hspace{1cm}
    \begin{tabular}{cccc}
        \toprule
        \multicolumn{4}{c}{\textbf{Loss weights}} \\
        \midrule
        $\omega_r$ & 1 & $\omega_\text{data}$ & 0 \\
        \cmidrule(lr){1-2} \cmidrule(lr){3-4}
        $\omega_b$ & 0 & $\omega_\text{sob}$ & 0 \\        
        \bottomrule
    \end{tabular}
    \caption{Network, training parameters (\cref{rmk:PINN_notations}) and loss weights for $u_\theta$ in the \textit{2D Elliptic case}.}\label{tab:paramtest3_2D}
\end{table}

\begin{rmrk}\label[rmrk]{rmk:ref_sol_Ell2D}
	Here, we do not know the analytical solution associated with the problem under consideration. So, in order to analyze the results obtained, we will define $u$ as a reference solution $u_{\text{ref}}$ obtained from a FEM solver on an over-refined mesh of characteristic mesh size $h_\text{ref}$ and with $k_\text{ref}$ polynomial order. In this section, we set $N_{\text{ref}}=\num{1000}$ (and the associated characteristic mesh size $h_{\text{ref}}$, as defined in \cref{rmk:Ell2D_N_nodes}) and $k_{\text{ref}}=3$. \RTwo{The linear systems resulting from the finite element discretization is solved using the Conjugate Gradient iterative method (tolerance : $10^{-13}$, maximum number of iterations : $\num{10000}$) preconditioned by the algebraic multigrid (AMG).}
\end{rmrk}

\subsubsection{Error estimates}\label{sec:Ell2D_error_estimates}

We first test the error estimates (\cref{lem:error_estimation_add}) for the following two sets of parameters, randomly chosen from~$\mathcal{M}$: $\bm{\mu}^{(1)}$ and $\bm{\mu}^{(2)}$. \RBoth{We perform the same tests as in the \cref{sec:Lap2D}, considering only the additive approach.} The results are presented in \cref{fig:case3} for a fixed $k \in \{1,2,3\}$ with $N \in \{16,32,64,128,256\}$, as presented in \cref{rmk:Ell2D_N_nodes}.

\begin{figure}[ht!]
	\centering
	\begin{subfigure}{0.48\linewidth}
		\centering
		\cvgFEMCorrAlldeg{fig_testcase2D_test3_cvg_FEM_case3_v1_param1.csv}{fig_testcase2D_test3_cvg_Corr_case3_v1_param1.csv}{1e-9}
		\caption{Case of $\mu^{(1)}$}
	\end{subfigure}
	\begin{subfigure}{0.48\linewidth}
		\centering
		\cvgFEMCorrAlldeg{fig_testcase2D_test3_cvg_FEM_case3_v1_param2.csv}{fig_testcase2D_test3_cvg_Corr_case3_v1_param2.csv}{5e-8}
		\caption{Case of $\mu^{(2)}$}
	\end{subfigure}
	\caption{Considering the \textit{2D elliptic case} and the PINN prior $u_\theta$. Left -- $L^2$ relative error on $h$, obtained with the standard FEM $e_h^{(1)}$ (solid lines) and the additive approach $e_{h,+}^{(1)}$ (dashed lines) for $\bm{\mu}^{(1)}$, with $k \in \{1,2,3\}$. Right -- Same for $\bm{\mu}^{(2)}$.}\label{fig:case3}
\end{figure}

As in the other test cases, the two approaches tested appear to respect the correct slopes of \cref{thm:classical_error_estimate} and \cref{lem:error_estimation_add}. The additive approach seems to be more efficient than the standard FEM for polynomial orders $k \in \{1,2,3\}$ and for the two sets of parameters considered.

\subsubsection{Comparison of different approaches}\label{sec:Ell2D_comparison}

\begin{figure}[ht!] \centering
    \includegraphics[scale=1]{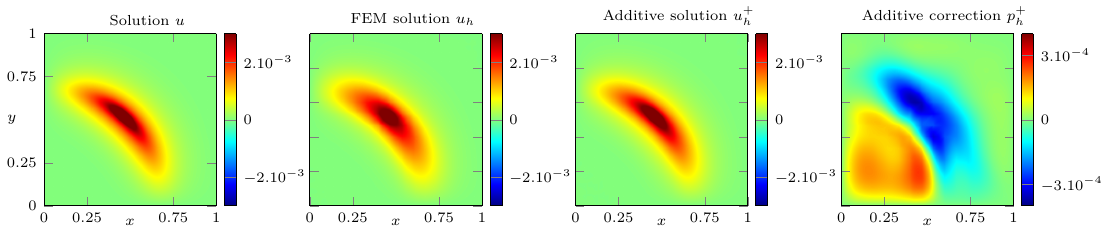}

    \includegraphics[scale=1]{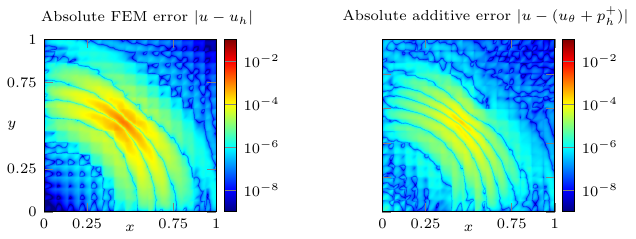}

	\caption{Considering the \textit{2D elliptic case} with $\bm{\mu}^{(3)}$, $k=2$, $N=16$ and the PINN prior $u_\theta$. Comparison of the solution obtained with the standard FEM and the additive approach with the analytical solution. For the additive method, comparison of the correction term with the analytical one.}\label{fig:case3_2D_plots}
\end{figure}

We perform the same comparison as in \cref{sec:Lap2Dlow} for this elliptic case. We focus on the third parameter $\bm{\mu}^{(3)}$ by taking a closer look at the solution obtained with the different approaches in \cref{fig:case3_2D_plots} considering $N=16$ and $k=2$.

We observe that the enriched FEM provides a more accurate solution compared to the standard FEM one. The results indicate that the additive approach is particularly effective in capturing the solution's finer details. This demonstrates its potential in solving anisotropic problems with higher accuracy than standard methods.

\subsubsection{Gains achieved with the additive approach}\label{sec:Ell2D_gains}

Considering a \RTwo{set $\mathcal{S}$ of $n_p=50$ parameter instances}, we will evaluate the gains $G_{+,\theta}$ and $G_+$ defined in~\eqref{eq:gain_add_num}. The results are presented in \cref{tab:case3} for $k \in \{1,2,3\}$ fixed and $N \in \{20,40\}$ fixed.

\begin{table}[ht!]
	\centering
	\GainsTableAlldeg{fig_testcase2D_test3_gains_gains_table_case3.csv}
	\caption{Considering the \textit{2D elliptic case}, $k\in\{1,2,3\}$ and the PINN prior $u_\theta$. Left -- Gains in $L^2$ relative error of the additive method with respect to PINN. Right -- Gains in $L^2$ relative error of our approach with respect to FEM.}\label{tab:case3}
\end{table}

As in the previous test cases, the additive approach seems to be more efficient than the standard FEM for the three polynomial orders $k \in \{1,2,3\}$ and for the two mesh sizes $N \in \{20,40\}$. However, the gains obtained are less significant than in the previous test cases. This is due to the fact that the problem under consideration is more complex, and thus, the prior $u_\theta$ is less accurate.
Namely, PINNs have trouble converging on such highly anisotropic problems,
see e.g.~\cite{ZhaAlK2022}.

\subsection{2D Poisson problem on an annulus, with mixed boundary conditions}\label{sec:Lap2DMixRing}

This section concerns the problem~\eqref{eq:ob_pde}, considering the Poisson problem with mixed (Dirichlet and Robin) boundary conditions defined in two space dimensions ($d=2$) by
\begin{equation*}
	\left\{
	\begin{aligned}
		-\Delta u & = f, \; &  & \text{in } \; \Omega \times \mathcal{M}, \\
		u         & = g, \;  &  & \text{on } \; \Gamma_E \times \mathcal{M}, \\
        \smash{\frac{\partial u}{\partial n}}+u  & = g_R, \;  &  & \text{on } \; \Gamma_I \times \mathcal{M},
	\end{aligned}
	\right.
\end{equation*}
with $\mathcal{M} \subset \mathbb{R}^p$ the parameter space (with $p$ the number of parameters). We consider~$\Omega$ to be an annulus centered at the origin, defined by the unit circle with a circular hole of radius $0.25$. Or, in other words, $\Omega$ is defined by
\[
	\Omega = \left\{ \bm{x} \in \mathbb{R}^2 \; \big| \; 0.25 \,\RBothN{<}\, \sqrt{x^2+y^2} \,\RBothN{<}\, 1 \right\}.
\]

We then define $\partial\Omega=\Gamma_I\cup\Gamma_E$ the boundary of $\Omega$, with $\Gamma_I$ the inner boundary (the hole) and $\Gamma_E$ the outer boundary (the unit circle). We consider the analytical solution defined for all $\bm{x}=(x,y)\in\Omega$ by
\begin{equation*}
	u(\bm{x},\bm{\mu})= 1 - \frac{\ln\big(\mu_1\sqrt{x^2+y^2}\big)}{\ln(4)},
\end{equation*}
with some parameters $\bm{\mu}=\mu_1\in[2.4, 2.6]$ ($p=1$ parameter), and the associated right-hand side $f=0$. The Dirichlet condition $g$ on $\Gamma_E$ and the Robin condition $g_R$ on $\Gamma_I$ are defined by
\begin{equation*}
	g(\bm{x},\bm{\mu})=1 - \frac{\ln(\mu_1)}{\ln(4)} \quad \text{and} \quad g_R(\bm{x},\bm{\mu})=2 + \frac{4-\ln(\mu_1)}{\ln(4)}.
\end{equation*}
These boundary conditions are thus parameter-dependent,
contrary to the previous cases,
which makes the problem more complex. \RTwo{To avoid geometric errors, we apply \cref{rem:bconcurved}, by considering that $g=u$ on $\Gamma_{E,h}$ and $g_R=\frac{\partial u}{\partial n}$ on~$\Gamma_{I,h}$, with $\Gamma_{E,h}$ and $\Gamma_{I,h}$ the respective outer and inner boundaries of $\Omega_h$, the domain covered by the mesh. Note also that $u_\theta$ is not exact on these approximate boundaries.}

In this section, we consider the additive approach, as presented in \cref{sec:additive_prior}, by considering the PINN prior $u_\theta$. We start by testing the error estimates in \cref{sec:Lap2DAnn_error_estimations}. Then, we compare the different approaches in \cref{sec:Lap2DAnn_comparison} and evaluate the gains obtained in \cref{sec:Lap2DAnn_gains} on a sample of parameters.

Since the problem under consideration is parametric,
we deploy a parametric PINN,
which depends on both the space variable $\bm{x}=(x,y) \in \Omega$
and the parameters $\bm{\mu}=\mu_1 \in \mathcal{M}$. To improve the derivatives' quality, we consider the Sobolev training presented in \cref{sec:sobolev_training}. Moreover, we strongly impose the Dirichlet boundary conditions,
as explained in \cref{sec:exact_imposition_of_BC}, by using the formulation proposed in~\cite{Sukumar_2022}.
To do this, we define the prior
\begin{equation}\label{eq:mixedformulation}
	u_{\theta} = \frac{\varphi_E}{\varphi_E+\varphi_I^2}\left[w_\theta+\varphi_I\big(w_\theta-\nabla\varphi_I\cdot\nabla w_\theta-g_R\big)\right] + \frac{\varphi_I^2}{\varphi_E+\varphi_I^2}g+\varphi_E\varphi_I^2w_\theta,
\end{equation}
where $w_\theta$ is the neural network under consideration and $\varphi_I$ and $\varphi_E$ are respectively the signed distance functions to $\Gamma_I$ and $\Gamma_E$ defined by
\begin{equation*}
	\varphi_I(\bm{x})=\sqrt{x^2+y^2}-0.25, \quad \varphi_E(\bm{x})=1-\sqrt{x^2+y^2},
\end{equation*}
which cancels out exactly on $\Gamma_I$ and $\Gamma_E$. \RTwo{Note that the level-sets considered are signed distance functions in this specific test case, which is not the case in the other test cases. In this test case, this is necessary because of the formulation proposed in~\eqref{eq:mixedformulation} by~\cite{Sukumar_2022}.}

In this case, we consider an MLP with $5$ layers and a tanh activation function with the hyperparameters defined in \cref{tab:paramtest5_2D};
we use the Adam optimizer~\cite{KinBa2015}. \RBoth{We consider the residual and Sobolev losses in the same way as in~\cref{sec:Lap2Dlowaug} with $N_\text{col}=\num{6000}$ collocation points.}

\begin{table}[ht!]
    \centering
    \begin{tabular}{cc}
        \toprule
        \multicolumn{2}{c}{\textbf{Network - MLP}} \\
        \midrule
        \textit{layers} & $40,40,40,40,40$ \\
        \cmidrule(lr){1-2}
        $\sigma$ & tanh \\
        \bottomrule
    \end{tabular}
    \hspace{1cm}
    \begin{tabular}{cccc}
        \toprule
        \multicolumn{4}{c}{\textbf{Training - with LBFGS}} \\
        \midrule
        \textit{lr} & 1e-2 & $n_\text{epochs}$ & \num{4000} \\
        \cmidrule(lr){1-2} \cmidrule(lr){3-4}
        \textit{decay} & 0.99 & $n_\text{switch}$ & \num{3000} \\
        \cmidrule(lr){1-2} \cmidrule(lr){3-4}
        $N_\text{col}$ & \num{6000} \\
        \bottomrule
    \end{tabular}
    \hspace{1cm}
    \begin{tabular}{cccc}
        \toprule
        \multicolumn{4}{c}{\textbf{Loss weights}} \\
        \midrule
        $\omega_r$ & 1 & $\omega_\text{data}$ & 0 \\
        \cmidrule(lr){1-2} \cmidrule(lr){3-4}
        $\omega_b$ & 0 & $\omega_\text{sob}$ & 0.1 \\        
        \bottomrule
    \end{tabular}
    \caption{Network, training parameters (\cref{rmk:PINN_notations}) and loss weights for $u_\theta$ in the \textit{2D Laplacian case on an Annulus}.}\label{tab:paramtest5_2D}
\end{table}

\subsubsection{Error estimates}\label{sec:Lap2DAnn_error_estimations}

We start by testing the error estimation of \cref{lem:error_estimation_add} for the following two sets of parameters,
uniformly selected from $\mathcal{M}$:
\begin{equation*}
	\bm{\mu}^{(1)}=(2.51) \quad \text{and} \quad \bm{\mu}^{(2)}=(2.54)\;,
\end{equation*}
by considering the PINN prior $u_\theta$. \RBoth{We perform the same tests as in the previous sections, considering only the additive approach.}
The results are presented in \cref{fig:case5} for fixed $k \in \{1,2,3\}$.

\begin{figure}[ht!]
	\centering
	\begin{subfigure}{0.48\linewidth}
		\centering
		\cvgFEMCorrAlldeg{fig_testcase2D_test5_cvg_FEM_case5_v2_param1.csv}{fig_testcase2D_test5_cvg_Corr_case5_v2_param1.csv}{1e-10}
		\caption{Case of $\bm{\mu}^{(1)}$}
	\end{subfigure}
	\begin{subfigure}{0.48\linewidth}
		\centering
		\cvgFEMCorrAlldeg{fig_testcase2D_test5_cvg_FEM_case5_v2_param2.csv}{fig_testcase2D_test5_cvg_Corr_case5_v2_param2.csv}{1e-10}
		\caption{Case of $\bm{\mu}^{(2)}$}
	\end{subfigure}
	\caption{Considering the \textit{2D Laplacian case on an Annulus} and the PINN prior $u_\theta$. Left -- $L^2$ relative error on $h$, obtained with the standard FEM $e_h^{(1)}$ (solid lines) and the additive approach $e_{h,+}^{(1)}$ (dashed lines) for $\bm{\mu}^{(1)}$, with $k \in \{1,2,3\}$. Right -- Same for $\bm{\mu}^{(2)}$.}\label{fig:case5}
\end{figure}

As expected, we see in \cref{fig:case5} that the error estimates are confirmed by the numerical results obtained with the standard FEM and the additive approach. The error decreases with the correct order of convergence for these two methods. Furthermore, the enriched approach provides a better accuracy than the standard FEM, as expected.
\subsubsection{Comparison of different approaches}\label{sec:Lap2DAnn_comparison}

\begin{figure}[ht!] 
	\centering
    \includegraphics[scale=1]{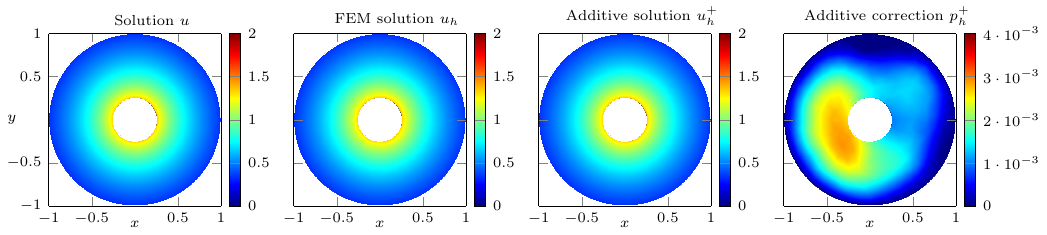}

	\includegraphics[scale=1]{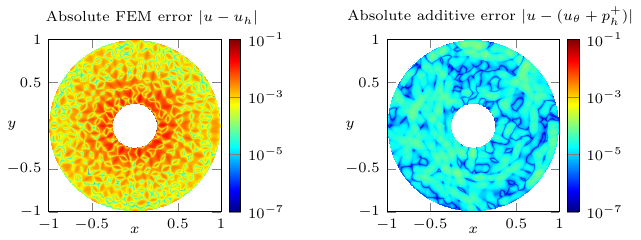}

	\caption{Considering the \textit{2D Laplacian case on an Annulus} with $\bm{\mu}^{(1)}$, $k=1$, $h\simeq 1.67\cdot 10^{-1}$ and the PINN prior $u_\theta$. Comparison of the solution obtained with the standard FEM and the additive approach with the analytical solution. For the additive method, comparison of the correction term with the analytical one.}\label{fig:case5_2D_plots}
\end{figure}

We perform the same comparison as in \cref{sec:Lap2Dlow} for this \ROne{L}aplacian case on an annulus. We focus on the first parameter $\bm{\mu}^{(1)}$ by taking a closer look at the solution obtained with the different approaches in \cref{fig:case5_2D_plots} considering $h\simeq 1.67\cdot 10^{-1}$ and $k=1$.

Once again, we observe that the enriched approach provides a significant improvement in accuracy compared to the standard FEM. This demonstrates the effectiveness of incorporating neural network priors in the case of mixed boundary conditions on more complex geometries than squares (here, on an annulus).

\subsubsection{Gains achieved with the additive approach} \label{sec:Lap2DAnn_gains}

Considering a \RTwo{set $\mathcal{S}$ of $n_p=50$ parameter instances},
we now evaluate the gains $G_{+,\theta}$ and $G_+$ defined in~\eqref{eq:gain_add_num}.
The results are presented in \cref{tab:case4}
for $k \in \{1,2,3\}$ and $h \in \{1.33\cdot 10^{-1},6.90\cdot 10^{-2}\}$.

\begin{table}[ht!]
	\centering
	\GainsTableAlldegh{fig_testcase2D_test5_gains_Tab_stats_case5_v2.csv}
	\caption{Considering the \textit{2D Laplacian case on an Annulus}, $k\in\{1,2,3\}$ and the PINN prior $u_\theta$. Left -- Gains in $L^2$ relative error of the additive method with respect to PINN. Right -- Gains in $L^2$ relative error of our approach with respect to FEM.}\label{tab:case4}
\end{table}

As in previous sections, the PINN-enriched approach seems to give better results than standard FEM. For $k=1$ we see a gain of $50$ on average for this test case, which is equivalent to refining the mesh by a factor of $7$ for $\mathbb{P}_1$ elements.

\subsection{\RTwo{3D Poisson problem in a cube domain}}\label{sec:Lap3D}

\RTwo{We now consider the problem \eqref{eq:Lap2D} of \cref{sec:Lap2Dlow}
but in three dimensions ($d=3$),
in the space domain \smash{$\Omega=\RBothN{(-0.5 \pi, 0.5 \pi)}^3$}
and the parameter domain $\mathcal{M} = [-0.5,0.5]^3$.
We define the right-hand side $f$ such that the solution is given
for all $\bm{x}=(x,y,z)\in\Omega$ and
$\bm{\mu}=(\mu_1,\mu_2,\mu_3)\in\mathcal{M}$ by
\begin{equation*}
    u(\bm{x},\bm{\mu})=\exp\left(-\frac{{(x-\mu_1)}^2+{(y-\mu_2)}^2+{(z-\mu_3)}^2}{2}\right)\sin(2x)\sin(2y)\sin(2z).
\end{equation*}}

\RTwo{In this section, we study the additive approach with $k=1$. In \cref{sec:Lap3D_error_estimations}, we start by testing the error estimates. Then, in \cref{sec:Lap3D_costs}, we conduct a study of the computation times of the different methods. In \cref{sec:Lap3D_gains}, we evaluate the gains achieved with the enriched approach. Finally, we discuss the importance of the prior quality in \cref{sec:Lap3D_badPINN}.}

\begin{rmrk}\label[rmrk]{rmk:Lap3D_N_nodes}
	\RTwo{In this section, the cube $\Omega$ of side length $\pi$ is discretized using a Cartesian mesh with $N^3$ nodes. Consequently, the characteristic mesh size is defined as a function of $N$ by $h=\frac{\pi\sqrt{3}}{N-1}$.}
\end{rmrk}

\RTwo{We deploy the same type of parametric PINN as in the previous sections, where we strongly impose the Dirichlet boundary conditions. This PINN depends on both the spatial variable $\bm{x}$ and the parameters $\bm{\mu}$.
We use the prior $u_\theta$ defined in~\eqref{eq:prior_with_levelset}, where we choose the level-set function $\varphi$ defined by
\begin{equation*}
	\varphi(\bm{x})=(x+0.5\pi)(x-0.5\pi)(y+0.5\pi)(y-0.5\pi)(z+0.5\pi)(z-0.5\pi).
\end{equation*}
The residual loss function~\eqref{eq:residual_loss_parametric_MC} is approximated with $N_\text{col}=\num{40000}$ collocation points uniformly chosen in~$\Omega\times\mathcal{M}$. The hyperparameters are defined in \cref{tab:paramtest1_3D}; we use the Adam optimizer and then switch to the LBFGS optimizer after the $n_\text{switch}$-th epoch.}

\begin{table}[htbp]
    \centering
    \begin{tabular}{cc}
        \toprule
        \multicolumn{2}{c}{\textbf{Network - MLP}} \\
        \midrule
        \textit{layers} & $40,60,60,60,40$ \\
        \cmidrule(lr){1-2}
        $\sigma$ & tanh \\
        \bottomrule
    \end{tabular}
    \hspace{1cm}
    \begin{tabular}{cccc}
        \toprule
        \multicolumn{4}{c}{\textbf{Training - with LBFGS}} \\
        \midrule
        \textit{lr} & 1.7e-2 & $n_\text{epochs}$ & \num{5000} \\
        \cmidrule(lr){1-2} \cmidrule(lr){3-4}
        \textit{decay} &  & $n_\text{switch}$ & \num{2000} \\
        \cmidrule(lr){1-2} \cmidrule(lr){3-4}
        $N_\text{col}$ & \num{40000} \\
        \bottomrule
    \end{tabular}
    \hspace{1cm}
    \begin{tabular}{cccc}
        \toprule
        \multicolumn{4}{c}{\textbf{Loss weights}} \\
        \midrule
        $\omega_r$ & 1 & $\omega_\text{data}$ & 0 \\
        \cmidrule(lr){1-2} \cmidrule(lr){3-4}
        $\omega_b$ & 0 & $\omega_\text{sob}$ & 0 \\        
        \bottomrule
    \end{tabular}
    \caption{\RTwo{Network, training parameters (\cref{rmk:PINN_notations}) and loss weights for $u_\theta$ in the \textit{3D Poisson problem}.}}\label{tab:paramtest1_3D}
\end{table}

\begin{rmrk}
	\RTwo{The results are obtained using the iterative conjugate gradient method (tolerance: $10^{-8}$, maximum number of iterations: \num{1000}) preconditioned by \ROneN{preconditioned by HYPRE's BoomerAMG (a classical algebraic multigrid method) with a strong connection threshold of $0.7$ and $4$ levels of non-local aggregation using $2$ paths (recommended options for the Poisson problem in \texttt{FEniCSx}\footnote{{\url{https://github.com/FEniCS/performance-test}}}). For both FE approaches, the right-hand side is interpolated on a $\mathbb{P}^2$ space. This means that} for the additive approach, the derivatives of the network are evaluated analytically (using \texttt{torch.autograd}) on the degrees of freedom associated to the space $\mathbb{P}^2$. To give an order of magnitude, the evaluation of the second derivatives of the prior, on the mesh with $N=100$ (i.e., \num{7880599} evaluation points for a $\mathbb{P}^2$ space) takes $1.01$ seconds on an NVIDIA H100 GPU.}
\end{rmrk}

\subsubsection{\RTwo{Error estimates}}\label{sec:Lap3D_error_estimations}

\RTwo{We start by testing the error estimates of \cref{lem:error_estimation_add} for the following two sets of parameters, randomly selected in $\mathcal{M}$:
\begin{equation*}
	\bm{\mu}^{(1)}=(0.05,0.22,0.1) \quad \text{and} \quad \bm{\mu}^{(2)}=(0.04,-0.08,0.15),
\end{equation*}
by considering the PINN prior $u_\theta$ and by performing the same tests as in the previous sections.
The results are presented in \cref{fig:case13D} with $h$ depending on $N\in\{20,40,60,80,100\}$, as presented in \cref{rmk:Lap3D_N_nodes}. We make the same observations as in the previous sections, both in terms of convergence rate and gains.}

\begin{figure}[ht!]
	\centering
	\begin{subfigure}{0.48\linewidth}
		\centering
		\cvgFEMCorrOnedeg{fig_testcase3D_cvg_FEM_case1_v2_param1_degree1.csv}{fig_testcase3D_cvg_Corr_case1_v2_param1_degree1.csv}{2e-5}
		\caption{Case of $\bm{\mu}^{(1)}$}
	\end{subfigure}
	\begin{subfigure}{0.48\linewidth}
		\centering
		\cvgFEMCorrOnedeg{fig_testcase3D_cvg_FEM_case1_v2_param2_degree1.csv}{fig_testcase3D_cvg_Corr_case1_v2_param2_degree1.csv}{2e-5}
		\caption{Case of $\bm{\mu}^{(2)}$}
	\end{subfigure}
	\caption{\RTwo{Considering the \textit{3D Poisson problem} and the PINN prior $u_\theta$. Left -- $L^2$ relative error on $h$, obtained with the standard FEM $e_h^{(1)}$ (solid lines) and the additive approach $e_{h,+}^{(1)}$ (dashed lines) for $\bm{\mu}^{(1)}$, with $k=1$. Right -- Same for $\bm{\mu}^{(2)}$.}}\label{fig:case13D}
\end{figure}

\subsubsection{\RTwo{Computation times of the different methods}}\label{sec:Lap3D_costs}


\RTwo{Just like in \cref{sec:Lap2Dlow}, we are now interested in the numerical costs of the different approaches, but this time in terms of computation time. \ROneN{First, we perform a study in a non-parametric framework, considering only the parameter $\bm{\mu}^{(1)}$. Thus,} in \cref{fig:case1_3D_time_error}, we examine the execution time of the two approaches (\ROneN{training time not included} for the additive approach), comparing it to the relative $L^2$ error. In \cref{tab:case1_3D_timerequired}, we perform a study similar to the one done previously. More precisely, for a given error $e$, we look at the numerical cost of the two approaches that allow us to achieve it. Therefore, for a given $e$, we report the characteristic mesh size $h$, the number of degrees of freedom $N_\text{dofs}$, and the execution time. \ROneN{In this purely online context on the parameter $\bm{\mu}^{(1)}$, the times indicated will include all the steps necessary to solve the EF problem, i.e., mesh construction, system assembly, solving the linear system, and, for the additive approach, the cost of evaluating the derivatives of the prior.}}

\ROneN{However, the offline cost of the additive approach cannot be neglected, as it includes the training time of the prior (in this case: $\num{707.84}$ seconds). This is why, in \cref{fig:case1_3D_time_param}, we are interested in the numerical costs of the two approaches in a parametric framework (since PINN is trained parametrically). To do this, we seek to determine, for a fixed error $e$, the number of parameter instances at which the additive approach is better than the standard approach, including the PINN training time.  More precisely, the offline cost of the standard approach consists simply of constructing the mesh (of the size required to achieve $e$), while the additive approach includes the construction of the (coarser) mesh as well as the training time. For the online phase, both approaches include all the FE steps (assembly and resolution of the system as well as the evaluation of derivatives for the enriched method). It should be noted that the online costs are much higher for the standard approach, due to a finer mesh size.}

\begin{figure}[ht!]
	\centering
	\cvgtimeerror{fig_testcase3D_time_error_FEM_case1_v2_param1_degree1_time_error.csv}{fig_testcase3D_time_error_Corr_case1_v2_param1_degree1_time_error.csv}
	\caption{\RTwo{Considering the \textit{3D Poisson problem} with $\bm{\mu}^{(1)}$, $k=1$ and the PINN prior $u_\theta$. $L^2$ relative error obtained with the standard FEM $e_h^{(1)}$ (solid lines) and the additive approach $e_{h,+}^{(1)}$ (dashed lines) as a function of the \ROneN{online} computation time (in seconds), \ROneN{including mesh construction}.}
	}\label{fig:case1_3D_time_error}
\end{figure}

\begin{table}[ht!]
	\centering
	\timerequired{fig_testcase3D_timerequired_TabTimes_case1_v2_param1_degree1.csv}
	\caption{\RTwo{Considering the \textit{3D Poisson problem} with $\bm{\mu}^{(1)}$, $k=1$ and the PINN prior $u_\theta$. Left -- Characteristic $N$ (associated to the characteristic mesh size $h$) required to reach a fixed error $e$ for standard FEM and the additive approach. Middle -- Number of degrees of freedom $N_\text{dofs}$ associated with each case. Right -- Execution time (in seconds) \ROneN{of the online phase} associated with each case. \ROneN{The results in brackets are extrapolations provided to give an order of magnitude.}}
	}\label{tab:case1_3D_timerequired}
\end{table}

\begin{figure}[ht!]
	\centering
	\begin{subfigure}{0.48\linewidth}
		\centering
		\cvgtimeerrorparam{fig_testcase3D_parametric_framework_output13.csv}
		\vspace{0.4cm}
		\caption{\ROneN{To achieve an error of $10^{-3}$}}
	\end{subfigure}
	\begin{subfigure}{0.48\linewidth}
		\centering
		\cvgtimeerrorparamD{fig_testcase3D_parametric_framework_output54.csv}
		\caption{\ROneN{To achieve an error of $5\cdot 10^{-4}$}}
	\end{subfigure}
	\caption{\ROneN{Computation time (in seconds) of the total process (including the offline phase) as a function of the number of parameter instances $n_p$ to be solved. The vertical dashed line represents the number of parameter instances from which the additive approach becomes more efficient than the standard FEM.}}\label{fig:case1_3D_time_param}
\end{figure}

\RTwo{In \cref{fig:case1_3D_time_error}, we observe that for a given error, the additive approach is significantly faster than the standard FEM. This is confirmed in \cref{tab:case1_3D_timerequired}, where we can see that a coarser mesh can be used to achieve a given accuracy with our enriched FEM compared to the standard FEM, leading to reduced computational time. For example, to reach an error of $10^{-3}$, we need $N=152$ (execution time: \ROneN{$76.5$} seconds) with the standard FEM, whereas with the additive approach, we only need $N=12$ (execution time: \ROneN{$0.052$} seconds). This represents a speed-up of approximately \ROneN{$1471$} times.}

\ROneN{However, the cost of training the network cannot be ignored. In \cref{fig:case1_3D_time_param}, we see that when we include a parametric context (including offline cost), the additive approach becomes more advantageous than the standard approach when using $n_p=19$ parameter sets for an error of $10^{-3}$ (left figure) and $n_p=5$ parameter sets for an error of $5\cdot 10^{-4}$ (right figure). And it is with this framework that approaches enriched by PINNs become interesting.}

\subsubsection{\RTwo{Gains achieved with the additive approach}}\label{sec:Lap3D_gains}

\RTwo{Considering a set $\mathcal{S}$ of $n_p=50$ parameter instances,
we now evaluate the gains $G_{+,\theta}$ and $G_+$ defined in~\eqref{eq:gain_add_num}.
The results are presented in \cref{tab:case1_3D}
for $k=1$ and $N \in \{20, 40\}$, where we observe gains relatively close to the 2D test case (\cref{sec:Lap2Dlow}).}

\begin{table}[ht!]
	\centering
	\GainsTableOnedeg{fig_testcase3D_gains_Tab_stats_case1_v2_degree1.csv}
	\caption{\RTwo{Considering the \textit{3D Poisson problem}, $k=1$ and the PINN prior $u_\theta$. Left -- Gains in $L^2$ relative error of the additive method with respect to PINN. Right -- Gains in $L^2$ relative error of our approach with respect to FEM.}}\label{tab:case1_3D}
\end{table}

\subsubsection{Influence of prior quality on the additive approach}\label{sec:Lap3D_badPINN}

\RTwo{Lastly, we focus on the impact of PINN quality on the results obtained with the additive approach. In \cref{tab:case1_3D_badPINNs}, we look at the average gain obtained with different priors on the same sample of $n_p=50$ parameters as in \cref{sec:Lap3D_gains}. More specifically, we train 9 PINNs (including the PINN from the previous sections) by varying only the number of epochs $n_\text{epochs}$ and the number of collocation points $N_\text{col}$. }

\begin{table}[ht!]
	\centering
	\begin{filecontents*}{fig_testcase3D_badPINNs_Tab_allv_mean_degree1.csv}
Ncol,epoch2500,epoch5000,epoch10000
10000,57.27497823210055,76.38175888666673,77.74352507358225
20000,78.33104088455771,87.78975112868149,88.26037092384931
40000,110.51681213588499,137.34388869761437,138.58717857719333
\end{filecontents*}
\pgfplotstabletypeset[
		col sep=comma,
		every head row/.style={
		before row={\toprule[1.pt]
		& \multicolumn{3}{c}{\textbf{$n_\text{epochs}$}} \\
		\cmidrule(lr){1-1} \cmidrule(lr){2-4}
		},
		after row=\cmidrule(lr){1-1} \cmidrule(lr){2-4}},
		every last row/.style={after row=\bottomrule[1.pt]},
		columns/Ncol/.style={column name=\textbf{$N_\text{col}$}},
		columns/epoch2500/.style={column name=\textbf{$\num{2500}$},fixed},
		columns/epoch5000/.style={column name=\textbf{$\num{5000}$},fixed},
		columns/epoch10000/.style={column name=\textbf{$\num{10000}$},fixed},
		columns={Ncol,epoch2500,epoch5000,epoch10000},
		precision=2
	]{fig_testcase3D_badPINNs_Tab_allv_mean_degree1.csv}

	\caption{\RTwo{Considering the \textit{3D Poisson problem}, $k=1$. Mean gain in $L^2$ relative error of the additive method with respect to FEM, by varying the quality of the PINN prior $u_\theta$ (by changing $n_\text{epochs}$ and $N_\text{col}$).}
	}\label{tab:case1_3D_badPINNs}
\end{table}

\RTwo{In \cref{tab:case1_3D_badPINNs}, we can see that the number of collocation points considered seems to have a significant impact on the gains obtained. On the other hand, the number of epochs appears to be less significant. The previous results, therefore, seem to have been obtained by choosing the training among the nine that gives the best gains, while having a relatively correct number of epochs.}

    \section{Conclusion and future work}\label{sec:conclusion}
    In this work, we explored a new approach combining FEM and predictions from neural networks.
The neural network prediction is used to enhance the FEM prediction, by correcting the FEM approximation space.
Two strategies were investigated: an additive correction and a multiplicative one.
For both approaches, we have proved a priori error estimates for both the $H^1$ semi-norm and the $L^2$ norm.
We have also highlighted a link between these two techniques.
Moreover, the constant appearing in these inequalities is compared with the case of classical FEM.
Numerical simulations on parametric problems in one, two and three dimensions confirm our theoretical analyses. The various numerical test cases have shown that PINNs are good candidates for our enriched methods due to their ability to approximate the derivatives of the solution, which is necessary for the quality of our a priori error estimates. \RBothN{The non-smooth one-dimensional transmission test case, which does not fit into the theoretical framework of the two enriched approaches, also showed a similar behavior when using an appropriate prior. In addition, the ability of PINNs} to approximate the solution of the parametric PDE over a set of parameters also showed that the proposed approaches are much more interesting in terms of numerical cost than the standard method. Solutions to improve the quality of the prior and, thus, the quality of the results have also been highlighted, with Sobolev training in particular. We have also observed that the additive approach offers greater robustness and a more straightforward implementation than the multiplicative one.

The present work opens up several perspectives. For instance, the additive and multiplicative can be easily adapted to non-linear equations. Moreover, the prediction could also be used to build an optimal mesh before the FEM resolution, for instance, via a posteriori error estimates. Furthermore, more complex geometries can be considered through the use of level-sets.
\RBoth{In this paper, we restrict our attention to regular shapes. For geometries with re-entrant corners, the corresponding singularities can be incorporated into the neural network, following the approach of \cite{reconns2024}, while a posteriori error estimators may be employed to drive mesh adaptation in the FEM correction. This direction will be explored in future work.}
\ROne{In this work, the prior is provided by the prediction of a PINN, which offers the advantage of accurately approximating derivatives. Alternatively, a reduced-order model could serve as a prior, provided its derivatives are sufficiently precise, although this option may become computationally demanding.}
\RTwo{Finally, the use of neural operators (e.g.\ Fourier Neural Operators, see~\cite{LiKovAziLiuBhaStuAna2021}) instead of neural networks could be an interesting alternative to explore, allowing for better generalization without being restricted to parameterized functions.}

    \section{Acknowledgements}

    As part of the ``France 2030'' initiative, this work has benefited from a national grant
    managed by the French National Research Agency (Agence Nationale de la Recherche) attributed
    to the Exa-MA project of the NumPEx PEPR program, under the reference ANR-22-EXNU-0002.
    F. F. acknowledges funding by the European Union
    with ERC Project \textsc{Incorwave} -- grant 101116288. This work was also supported by the Agence Nationale de la Recherche, Project PhiFEM, under grant ANR-22- CE46-0003-01.

    \bibliographystyle{alpha}
    \bibliography{bibliography}

\newcommand{\etalchar}[1]{$^{#1}$}
\begin{thebibliography}{BBDVC{\etalchar{+}}06}

\bibitem[Aea15]{AlnBle2015}
M.~Alnæs and J.~Blechta et~al.
\newblock {The FEniCS Project Version 1.5}.
\newblock {\em Archive of Numerical Software}, 3(100), 2015.

\bibitem[AFH{\etalchar{+}}25]{AghFraHilMicVig2025}
J.~Aghili, E.~Franck, R.~Hild, V.~Michel-Dansac, and V.~Vigon.
\newblock Accelerating the convergence of {N}ewton's method for nonlinear
  elliptic {PDE}s using {F}ourier neural operators.
\newblock {\em Commun. Nonlinear Sci.}, 140(2):108434, 2025.

\bibitem[AL{\O}{\etalchar{+}}14]{alnaes_unified_2014}
M.~S. Aln{\ae}s, A.~Logg, K.~B. {\O}lgaard, M.~E. Rognes, and G.~N. Wells.
\newblock Unified form language: {A} domain-specific language for weak
  formulations of partial differential equations.
\newblock {\em ACM Trans. Math. Softw.}, 40(2):1--37, 2014.

\bibitem[AS97]{ALMEIDA1997291}
R.~C. Almeida and R.~S. Silva.
\newblock A stable {P}etrov-{G}alerkin method for convection-dominated
  problems.
\newblock {\em Comput. Method. Appl. M.}, 140(3--4):291--304, 1997.

\bibitem[BBDVC{\etalchar{+}}06]{bazilevs2006isogeometric}
Y.~Bazilevs, L.~Beirão Da~Veiga, J.~A. Cottrell, T.~J.~R. Hughes, and
  G.~Sangalli.
\newblock Isogeometric analysis: approximation, stability and error estimates
  for $h$-refined meshes.
\newblock {\em Math. Models Methods Appl. Sci.}, 16(07):1031--1090, 2006.

\bibitem[BBO04]{babuvska2004generalized}
Ivo Babu{\v{s}}ka, Uday Banerjee, and John~E Osborn.
\newblock Generalized finite element methods—main ideas, results and
  perspective.
\newblock {\em International Journal of Computational Methods}, 1(01):67--103,
  2004.

\bibitem[Bea23]{baratta_dolfinx_2023}
I.~A. Baratta and J.~P.~Dean et~al.
\newblock {DOLFINx: The next generation FEniCS problem solving environment},
  2023.

\bibitem[BGV09]{belytschko2009review}
Ted Belytschko, Robert Gracie, and Giulio Ventura.
\newblock A review of extended/generalized finite element methods for material
  modeling.
\newblock {\em Modelling and Simulation in Materials Science and Engineering},
  17(4):043001, 2009.

\bibitem[BHL93]{berkooz1993proper}
G.~Berkooz, P.~Holmes, and J.~L. Lumley.
\newblock {The Proper Orthogonal Decomposition in the Analysis of Turbulent
  Flows}.
\newblock {\em Annu. Rev. Fluid Mech.}, 25(1):539--575, 1993.

\bibitem[BLM24]{BadLiMar2024}
S.~Badia, W.~Li, and A.~F. Martín.
\newblock Finite element interpolated neural networks for solving forward and
  inverse problems.
\newblock {\em Comput. Method. Appl. M.}, 418:116505, 2024.

\bibitem[BM97]{babuvska1997partition}
I.~Babuška and J.~M. Melenk.
\newblock The partition of unity method.
\newblock {\em Int. J. Numer. Meth. Eng.}, 40(4):727--758, 1997.

\bibitem[BMP{\etalchar{+}}19]{brunet2019physics}
J.-N. Brunet, A.~Mendizabal, A.~Petit, N.~Golse, \'E. Vibert, and S.~Cotin.
\newblock Physics-based deep neural network for augmented reality during liver
  surgery.
\newblock In {\em Medical Image Computing and Computer Assisted Intervention -
  MICCAI 2019}, pages 137--145. Springer International Publishing, 2019.

\bibitem[BS08]{brenner2008mathematical}
S.~C. Brenner and L.~R. Scott.
\newblock {\em {The Mathematical Theory of Finite Element Methods}}.
\newblock Springer New York, 2008.

\bibitem[Caf98]{Caf1998}
R.~E. Caflisch.
\newblock {Monte Carlo and quasi-Monte Carlo methods}.
\newblock {\em Acta Numer.}, 7:1--49, 1998.

\bibitem[CDCG{\etalchar{+}}22]{cuomo2022scientific}
S.~Cuomo, V.~S. Di~Cola, F.~Giampaolo, G.~Rozza, M.~Raissi, and F.~Piccialli.
\newblock {Scientific Machine Learning Through Physics-Informed Neural
  Networks: Where we are and What's Next}.
\newblock {\em J. Sci. Comput.}, 92(3), 2022.

\bibitem[CDG08]{Cockburn2008}
B.~Cockburn, B.~Dong, and J.~Guzmán.
\newblock {A superconvergent LDG-hybridizable Galerkin method for second-order
  elliptic problems}.
\newblock {\em Math. Comput.}, 77(264):1887--1916, 2008.

\bibitem[CDL{\etalchar{+}}23]{cotin2023phi}
S.~Cotin, M.~Duprez, V.~Lleras, A.~Lozinski, and K.~Vuillemot.
\newblock {$\phi$-FEM: An Efficient Simulation Tool Using Simple Meshes for
  Problems in Structure Mechanics and Heat Transfer}.
\newblock In {\em Partition of Unity Methods}, pages 191--216. Wiley Online
  Library, 2023.

\bibitem[Cia02]{ciarlet2002finite}
P.~G. Ciarlet.
\newblock {\em The Finite Element Method for Elliptic Problems}.
\newblock Society for Industrial and Applied Mathematics, 2002.

\bibitem[CLC{\etalchar{+}}16]{CanLeyChiGonCueFeuBerHue2016}
D.~Canales, A.~Leygue, D.~Chinesta, D.~González, E.~Cueto, E.~Feulvarch, J.-M.
  Bergheau, and A.~Huerta.
\newblock {Vademecum-based GFEM (V-GFEM): optimal enrichment for transient
  problems}.
\newblock {\em Int. J. Numer. Meth. Eng.}, 108(9):971--989, 2016.

\bibitem[Dem23]{demkowicz2023mathematical}
L.~F. Demkowicz.
\newblock {\em {Mathematical Theory of Finite Elements}}.
\newblock Society for Industrial and Applied Mathematics, 2023.

\bibitem[DHMM24]{DolHeiMisMos2024}
V.~Dolean, A.~Heinlein, S.~Mishra, and B.~Moseley.
\newblock Multilevel domain decomposition-based architectures for
  physics-informed neural networks.
\newblock {\em Comput. Method. Appl. M.}, 429:117116, 2024.

\bibitem[DL20]{duprez2020phi}
M.~Duprez and A.~Lozinski.
\newblock {$\phi$-FEM: A Finite Element Method on Domains Defined by
  Level-Sets}.
\newblock {\em SIAM J. Numer. Anal.}, 58(2):1008--1028, 2020.

\bibitem[DLL22]{duprez2023new}
M.~Duprez, V.~Lleras, and A.~Lozinski.
\newblock A new $\phi$-{FEM} approach for problems with natural boundary
  conditions.
\newblock {\em Numer. Meth. Part. D. E.}, 39(1):281--303, 2022.

\bibitem[DLL23]{duprez2023phi}
M.~Duprez, V.~Lleras, and A.~Lozinski.
\newblock $\phi$-{FEM}: an optimally convergent and easily implementable
  immersed boundary method for particulate flows and {S}tokes equations.
\newblock {\em ESAIM: Mathematical Modelling and Numerical Analysis},
  57(3):1111--1142, 2023.

\bibitem[DLLV23]{DupLleLozVui2023}
M.~Duprez, V.~Lleras, A.~Lozinski, and K.~Vuillemot.
\newblock $\phi$-{FEM} for the heat equation: optimal convergence on unfitted
  meshes in space.
\newblock {\em C. R. Math.}, 361(G11):1699--1710, 2023.

\bibitem[DRMM24]{DeRMisMol2024}
T.~De~Ryck, S.~Mishra, and R.~Molinaro.
\newblock {wPINNs: Weak Physics Informed Neural Networks for Approximating
  Entropy Solutions of Hyperbolic Conservation Laws}.
\newblock {\em SIAM J. Numer. Anal.}, 62(2):811--841, 2024.

\bibitem[EG04]{Ern2004TheoryAP}
A.~Ern and J.-L. Guermond.
\newblock {\em {Theory and Practice of Finite Elements}}.
\newblock Springer New York, 2004.

\bibitem[ES24]{ern2024convergence}
A.~Ern and M.~Steins.
\newblock {Convergence Analysis for the Wave Equation Discretized with Hybrid
  Methods in Space (HHO, HDG and WG) and the Leapfrog Scheme in Time}.
\newblock {\em J. Sci. Comput.}, 101(1), 2024.

\bibitem[Eva22]{evans2022partial}
L.~C. Evans.
\newblock {\em Partial differential equations}.
\newblock Number~19 in Graduate studies in mathematics. American Mathematical
  Society, Providence, Rhode Island, second edition, 2022.

\bibitem[EY18]{e2017deepritzmethoddeep}
W.~E and B.~Yu.
\newblock {The Deep Ritz Method: A Deep Learning-Based Numerical Algorithm for
  Solving Variational Problems}.
\newblock {\em Commun. Math. Stat.}, 6(1):1--12, 2018.

\bibitem[FB10]{fries2010extended}
T.-P. Fries and T.~Belytschko.
\newblock The extended/generalized finite element method: {A}n overview of the
  method and its applications.
\newblock {\em Int. J. Numer. Meth. Eng.}, 84(3):253--304, 2010.

\bibitem[FBCD22]{Frambati2022practical}
S.~Frambati, H.~Barucq, H.~Calandra, and J.~Diaz.
\newblock Practical unstructured splines: {A}lgorithms, multi-patch spline
  spaces, and some applications to numerical analysis.
\newblock {\em J. Comput. Phys.}, 471:111625, 2022.

\bibitem[FMDN24]{FraMicNav2024}
E.~Franck, V.~Michel-Dansac, and L.~Navoret.
\newblock {Approximately well-balanced Discontinuous Galerkin methods using
  bases enriched with Physics-Informed Neural Networks}.
\newblock {\em J. Comput. Phys.}, 512:113144, 2024.

\bibitem[FST{\etalchar{+}}24]{feng_hybrid_2024}
X.~Feng, H.~Shangguan, T.~Tang, X.~Wan, and T.~Zhou.
\newblock A hybrid {FEM}-{PINN} method for time-dependent partial differential
  equations, 2024.
\newblock arXiv:2409.02810 [math].

\bibitem[GKLS24]{grossmann2023can}
T.~G. Grossmann, U.~J. Komorowska, J.~Latz, and C.-B. Schönlieb.
\newblock Can physics-informed neural networks beat the finite element method?
\newblock {\em IMA J. Appl. Math.}, 89(1):143--174, 2024.

\bibitem[HCB05]{hughes2005isogeometric}
T.~J.~R. Hughes, J.~A. Cottrell, and Y.~Bazilevs.
\newblock {Isogeometric analysis: CAD, finite elements, NURBS, exact geometry
  and mesh refinement}.
\newblock {\em Comput. Method. Appl. M.}, 194(39-41):4135--4195, 2005.

\bibitem[HMP12]{hiptmair2013error}
R.~Hiptmair, A.~Moiola, and I.~Perugia.
\newblock Error analysis of {T}refftz-discontinuous {G}alerkin methods for the
  time-harmonic {M}axwell equations.
\newblock {\em Math. Comput.}, 82(281):247--268, 2012.

\bibitem[HPS17]{hungria2017hdg}
A.~Hungria, D.~Prada, and F.-J. Sayas.
\newblock {HDG methods for elastodynamics}.
\newblock {\em Comput. Math. Appl.}, 74(11):2671--2690, 2017.

\bibitem[IGMS22]{ImbMoiSto2022}
L.-M. Imbert-G{\'{e}}rard, A.~Moiola, and P.~Stocker.
\newblock A space{\textendash}time quasi-{T}refftz {DG} method for the wave
  equation with piecewise-smooth coefficients.
\newblock {\em Math. Comp.}, 92(341):1211--1249, 2022.

\bibitem[JKK20]{JagKhaKar2020}
A.~D. Jagtap, E.~Kharazmi, and G.~E. Karniadakis.
\newblock Conservative physics-informed neural networks on discrete domains for
  conservation laws: {A}pplications to forward and inverse problems.
\newblock {\em Comput. Methods Appl. Mech. Engrg.}, 365:113028, 2020.

\bibitem[JKN18]{JohKnoNov2018}
V.~John, P.~Knobloch, and J.~Novo.
\newblock Finite elements for scalar convection-dominated equations and
  incompressible flow problems: a never ending story?
\newblock {\em Comput. Vis. Sci.}, 19(5--6):47--63, 2018.

\bibitem[Joh02]{faadibruno2002}
Warren Johnson.
\newblock The curious history of faa di bruno's formula.
\newblock {\em American Mathematical Monthly}, 109:217--234, 03 2002.

\bibitem[KB15]{KinBa2015}
D.~Kingma and J.~Ba.
\newblock {Adam: A Method for Stochastic Optimization}.
\newblock In {\em International Conference on Learning Representations (ICLR)},
  San Diego, CA, USA, 2015.

\bibitem[KZK21]{KhaZhaKar2021}
E.~Kharazmi, Z.~Zhang, and G.~E.~M. Karniadakis.
\newblock $hp$-{VPINNs}: {V}ariational physics-informed neural networks with
  domain decomposition.
\newblock {\em Comput. Methods Appl. Mech. Engrg.}, 374:113547, 2021.

\bibitem[LK90]{LeeKan1990}
H.~Lee and I.~S. Kang.
\newblock Neural algorithm for solving differential equations.
\newblock {\em J. Comput. Phys.}, 91(1):110--131, 1990.

\bibitem[LKA{\etalchar{+}}21]{LiKovAziLiuBhaStuAna2021}
Z.~Li, N.~B. Kovachki, K.~Azizzadenesheli, B.~Liu, K.~Bhattacharya, A.~Stuart,
  and A.~Anandkumar.
\newblock {Fourier Neural Operator for Parametric Partial Differential
  Equations}.
\newblock In {\em International Conference on Learning Representations}, 2021.

\bibitem[LLF98]{LagLikFot1998}
I.~E. Lagaris, A.~Likas, and D.~I. Fotiadis.
\newblock Artificial neural networks for solving ordinary and partial
  differential equations.
\newblock {\em IEEE Trans. Neural Netw.}, 9(5):987--1000, 1998.

\bibitem[MF94]{MeaFer1994}
A.~J. Meade and A.~A. Fernandez.
\newblock The numerical solution of linear ordinary differential equations by
  feedforward neural networks.
\newblock {\em Math. Comput. Model.}, 19(12):1--25, 1994.

\bibitem[MJLR24]{MarJenLesRic2024}
N.~Margenberg, R.~Jendersie, C.~Lessig, and T.~Richter.
\newblock {DNN-MG: A hybrid neural network/finite element method with
  applications to 3D simulations of the Navier-Stokes equations}.
\newblock {\em Comput. Method. Appl. M.}, 420:116692, 2024.

\bibitem[MP17]{moiola2018space}
A.~Moiola and I.~Perugia.
\newblock {A space--time Trefftz discontinuous Galerkin method for the acoustic
  wave equation in first-order formulation}.
\newblock {\em Numer. Math.}, 138(2):389--435, 2017.

\bibitem[Pea19]{paszke2019pytorchimperativestylehighperformance}
A.~Paszke and S.~Gross et~al.
\newblock {\em {PyTorch: an imperative style, high-performance deep learning
  library}}, pages 8026--8037.
\newblock Curran Associates Inc., Red Hook, NY, USA, 2019.

\bibitem[PFB24]{Pham2024stabilization}
H.~Pham, F.~Faucher, and H.~Barucq.
\newblock {Numerical investigation of stabilization in the Hybridizable
  Discontinuous Galerkin method for linear anisotropic elastic equation}.
\newblock {\em Comput. Method. Appl. M.}, 428:117080, 2024.

\bibitem[PFS{\etalchar{+}}19]{park2019deepsdflearningcontinuoussigned}
J.~J. Park, P.~Florence, J.~Straub, R.~Newcombe, and S.~Lovegrove.
\newblock {DeepSDF: Learning Continuous Signed Distance Functions for Shape
  Representation}.
\newblock In {\em 2019 IEEE/CVF Conference on Computer Vision and Pattern
  Recognition (CVPR)}, pages 165--174. IEEE, 2019.

\bibitem[PR07]{patera2007reduced}
A.~T. Patera and E.~M. Rønquist.
\newblock Reduced basis approximation and \emph{a posteriori} error estimation
  for a {B}oltzmann model.
\newblock {\em Comput. Method. Appl. M.}, 196(29--30):2925--2942, 2007.

\bibitem[PRV{\etalchar{+}}01]{prud2002reliable}
C.~Prud'homme, D.~V. Rovas, K.~Veroy, L.~Machiels, Y.~Maday, A.~T. Patera, and
  G.~Turinici.
\newblock {Reliable Real-Time Solution of Parametrized Partial Differential
  Equations: Reduced-Basis Output Bound Methods}.
\newblock {\em J. Fluids Eng.}, 124(1):70--80, 2001.

\bibitem[Rav00]{ravindran2000reduced}
S.~S. Ravindran.
\newblock A reduced-order approach for optimal control of fluids using proper
  orthogonal decomposition.
\newblock {\em Int. J. Numer. Meth. Fl.}, 34(5):425--448, 2000.

\bibitem[Red19]{j2005introduction}
J.~N. Reddy.
\newblock {\em Introduction to the Finite Element Method, Fourth Edition}.
\newblock McGraw-Hill Education, New York, N.Y., 4th edition. edition, 2019.

\bibitem[RPK19]{RAISSI2019686}
M.~Raissi, P.~Perdikaris, and G.~E. Karniadakis.
\newblock Physics-informed neural networks: {A} deep learning framework for
  solving forward and inverse problems involving nonlinear partial differential
  equations.
\newblock {\em J. Comput. Phys.}, 378:686--707, 2019.

\bibitem[SBRW22]{scroggs_basix_2022}
M.~W. Scroggs, I.~A. Baratta, C.~N. Richardson, and G.~N. Wells.
\newblock Basix: a runtime finite element basis evaluation library.
\newblock {\em J. Open Source Softw.}, 7(73):3982, 2022.

\bibitem[SCB01]{strouboulis2001generalized}
T.~Strouboulis, K.~Copps, and I.~Babuška.
\newblock The generalized finite element method.
\newblock {\em Comput. Method. Appl. M.}, 190(32-33):4081--4193, 2001.

\bibitem[SDRW22]{scroggs_construction_2022}
M.~W. Scroggs, J.~S. Dokken, C.~N. Richardson, and G.~N. Wells.
\newblock {Construction of Arbitrary Order Finite Element Degree-of-Freedom
  Maps on Polygonal and Polyhedral Cell Meshes}.
\newblock {\em ACM Trans. Math. Softw.}, 48(2):1--23, 2022.

\bibitem[SJHH21]{son2021sobolevtrainingphysicsinformed}
H.~Son, J.~W. Jang, W.~J. Han, and H.~J. Hwang.
\newblock Sobolev training for physics informed neural networks, 2021.

\bibitem[SKPP24]{sikora2024comparison}
M.~Sikora, P.~Krukowski, A.~Paszyńska, and M.~Paszyński.
\newblock {Comparison of Physics Informed Neural Networks and Finite Element
  Method Solvers for advection-dominated diffusion problems}.
\newblock {\em J. Comput. Sci.}, 81:102340, 2024.

\bibitem[SMB{\etalchar{+}}20]{sitzmann2020implicitneuralrepresentationsperiodic}
V.~Sitzmann, J.~Martel, A.~Bergman, D.~Lindell, and G.~Wetzstein.
\newblock Implicit neural representations with periodic activation functions.
\newblock In H.~Larochelle, M.~Ranzato, R.~Hadsell, M.~F. Balcan, and H.~Lin,
  editors, {\em Advances in Neural Information Processing Systems}, volume~33,
  pages 7462--7473. Curran Associates, Inc., 2020.

\bibitem[SMMB00]{sukumar2000extended}
N.~Sukumar, N.~Moës, B.~Moran, and T.~Belytschko.
\newblock Extended finite element method for three-dimensional crack modelling.
\newblock {\em Int. J. Numer. Meth. Eng.}, 48(11):1549--1570, 2000.

\bibitem[SS22]{Sukumar_2022}
N.~Sukumar and A.~Srivastava.
\newblock Exact imposition of boundary conditions with distance functions in
  physics-informed deep neural networks.
\newblock {\em Comput. Method. Appl. M.}, 389:114333, 2022.

\bibitem[SY06]{nocedal_quasi_newton_2006}
W.~Sun and Y.-X. Yuan.
\newblock Quasi-{Newton} {Methods}.
\newblock In {\em Numerical {Optimization}}, New York, NY, 2006. Springer.

\bibitem[Tea20]{TanSri2020}
M.~Tancik and P.~Srinivasan et~al.
\newblock {Fourier Features Let Networks Learn High Frequency Functions in Low
  Dimensional Domains}.
\newblock In {\em Advances in Neural Information Processing Systems},
  volume~33, pages 7537--7547. Curran Associates, Inc., 2020.

\bibitem[TPMM25]{reconns2024}
J.~M. Taylor, D.~Pardo, and J.~Muñoz-Matute.
\newblock {Regularity-conforming neural networks (ReCoNNs) for solving partial
  differential equations}, 2025.

\bibitem[vDSG25]{skardova_finite_2024}
K.~\v{S}kardová, A.~Daby-Seesaram, and M.~Genet.
\newblock {Finite element neural network interpolation: Part I---interpretable
  and adaptive discretization for solving PDEs}.
\newblock {\em Comput. Mech.}, 2025.

\bibitem[Wah95]{wahlbin2006superconvergence}
L.~B. Wahlbin.
\newblock {\em {Superconvergence in Galerkin Finite Element Methods}}.
\newblock Springer Berlin Heidelberg, 1995.

\bibitem[WLZ25]{WanLiZha2025}
D.~Wang, H.~Li, and Q.~Zhang.
\newblock {General enrichments of stable GFEM for interface problems: Theory
  and extreme learning machine construction}.
\newblock {\em Appl. Numer. Math.}, 214:143--159, 2025.

\bibitem[XLBJ25]{XioLonBorJia2025}
W.~Xiong, X.~Long, S.~P.~A. Bordas, and C.~Jiang.
\newblock The deep finite element method: {A} deep learning framework
  integrating the physics-informed neural networks with the finite element
  method.
\newblock {\em Comput. Method. Appl. M.}, 436:117681, 2025.

\bibitem[XPYZ25]{XiaPenYaoZho2025}
Z.~Xiang, W.~Peng, W.~Yao, and W.~Zhou.
\newblock {Hybrid Finite-Difference Physics-Informed Neural Networks Partial
  Differential Equation Solver for Complex Geometries}.
\newblock {\em J. Thermophys. Heat Tr.}, pages 1--18, 2025.

\bibitem[ZAK22]{ZhaAlK2022}
W.~Zhang and M.~Al~Kobaisi.
\newblock {On the Monotonicity and Positivity of Physics-Informed Neural
  Networks for Highly Anisotropic Diffusion Equations}.
\newblock {\em Energies}, 15(18):6823, 2022.

\end{thebibliography}

    \appendix
\section{Notations and definitions}\label[appendix]{app:notations}

The aim of this section is to introduce the notations used throughout the paper. We present the notations related to the parametric PDE (\cref{tab:notations_PDE}), to the neural network (\cref{tab:notations_PINN}), and to the finite element methods (\cref{tab:notations_FEM}).

\renewcommand{\arraystretch}{1.1}  

\begin{table}[ht!]
    \centering
    \begin{tabular}{c|c}
        \textbf{Notation} & \textbf{Definition} \\
        \hline
        $\Omega$ & Spatial domain \\
        $d$ & Spatial dimension \\
        $\bm{x}=(x_1,\dots,x_d)$ & Spatial coordinates \\
        \hline
        $\mathcal{M}$ & Parameter space \\
        $p$ & Number of parameters \\
        $\bm{\mu}=(\mu_1,\ldots,\mu_p)$ & Parameter vector \\
        \hline
        $M$ & Lifting constant \\
        $u$ & Solution of the problem \\
        $u_M$ & Solution of the lifted problem (by $M$) \\
        $f$ & Right-hand side of the problem \\
        $\mathcal{L}$ & Parametric differential operator of the problem \\
        $R$ & Reaction coefficient \\
        $C$ & Convection coefficient \\
        $D$ & Diffusion matrix \\
        Pe & Péclet number \\
    \end{tabular}
    \caption{Notations introduced for the parametric PDE.}
    \label{tab:notations_PDE}
\end{table}

\begin{table}[ht!]
    \centering
    \begin{tabular}{c|c}
        \textbf{Notation} & \textbf{Description} \\
        \hline
        $u_\theta$ & Neural network prediction of $u$ \\
        $u_{\theta,M}$ & Neural network prediction of $u_M$ \\
        \multirow{2}{*}{$\varphi$} & Level-set function used to impose BCs  \\
        & and generate sampling in PINNs \\
        $\theta$ & Trainable parameters of the neural network \\
        $\theta^\star$ & Optimal parameters of $\theta$ \\
        \hline
        $J_r$ & Residual loss \\
        $J_b$ & Boundary loss \\
        $J_\text{data}$ & Data loss \\
        $J_\text{\rm sob} $ & Sobolev loss \\
    \end{tabular}
    \caption{Notations considered for the neural network.}
    \label{tab:notations_PINN}
\end{table}

\begin{table}[ht!]
    \centering
    \begin{tabular}{c|c|c}
        & \textbf{Notation} & \textbf{Description} \\
        \hline
        \multirow{5}{*}{\rotatebox{90}{\small Standard FEM}}
        & $V_h^0$ & Finite element approximation space \\
        & $u_h$ & Finite element approximation of $u$ \\
        & $h$ & Characteristic size of the mesh \\
        & $\mathcal{I}_h$ & Lagrange interpolation operator \\
        & $k$ & Polynomial degree of the finite element approximation \\
        \hline
        \multirow{4}{*}{\rotatebox{90}{\small Additive}}\;\multirow{4}{*}{\rotatebox{90}{\small enrichment}} & $V_h^+$ & Finite element approximation space enriched with additive prior \\
        & $u_h^+$ & Finite element approximation of $u$ in $V_h^+$ \\
        & $p_h^+$ & Finite element approximation of $u-u_\theta$ in $V_h^0$ \\
        & $C_\text{\rm gain}^+$ & Additive gain constant \\
        \hline
        \multirow{6}{*}{\rotatebox{90}{\small Multiplicative}}\;\multirow{6}{*}{\rotatebox{90}{\small enrichment}} & $V_h^\times$ & Finite element approximation space enriched with multiplicative prior \\
        & $u_h^\times$ & Finite element approximation of $u$ in $V_h^\times-M$ \\
        & $p_h^\times$ & Finite element approximation of $u_M/u_{\theta,M}$ in $1+V_h^0$ \\
        & $C_{\text{\rm gain},H^1}^\times$ & Multiplicative gain constant in $H^1$ semi-norm \\
        & $C_{\text{\rm gain},L^2}^\times$ & Multiplicative gain constant in $L^2$ norm \\
        & $\tilde{\mathcal{I}}_h$ & Modified Lagrange interpolation operator \\
    \end{tabular}
    \caption{Notations used in the various finite element methods.}
    \label{tab:notations_FEM}
\end{table}

\end{document}